\documentclass[smallextended]{svjour3}       

\smartqed  

\usepackage{graphicx}
\usepackage{amssymb}
\usepackage[utf8]{inputenc}
\usepackage{dsfont}
\usepackage[T1]{fontenc}
\usepackage{color}
\usepackage{amsmath}
\usepackage{caption}
\usepackage[labelformat=simple]{subcaption}
 
\usepackage{multirow}
\usepackage{indentfirst}
\usepackage{makecell}
\usepackage{algpseudocode,algorithm}
\usepackage{etoolbox} 
\usepackage{fancyhdr}
\usepackage{pdfpages}
\usepackage{amssymb}
\usepackage[all]{xy}
\usepackage{tikz}
\usepackage{array}
\usepackage{overpic}
\usepackage{bm}
\usepackage{comment}
\usepackage{placeins}
\usepackage{flafter}
\usepackage{mathptmx}     
\usepackage{enumerate}
\usepackage[pagewise]{lineno}

\begin{document}

\title{An efficient topology optimization algorithm for large-scale three-dimensional structures}

\titlerunning{TO algorithm for large-scale 3D structures}

\author{A. Vitorino \and F. A. M. Gomes}

\institute{A. Vitorino \at
              Institute of Mathematics, Statistics and Scientific Computing, State University of Campinas, Campinas, SP, Brazil \\
              \email{oalfredovitorino@gmail.com}
           \and
           F. A. M. Gomes \at
              Institute of Mathematics, Statistics and Scientific Computing, State University of Campinas, Campinas, SP, Brazil \\
              \email{chico2@unicamp.br}
}

\date{Received: date / Accepted: date}

\maketitle
 
\begin{abstract}
Designing the topology of three-dimensional structures is a challenging problem due to its memory and time consumption. In this paper, we present a robust and efficient algorithm for solving large-scale 3D topology optimization problems. The robustness of the algorithm is ensured by adopting a globally convergent sequential linear programming method with a stopping criterion based on the first-order optimality conditions of the nonlinear problem. To increase the algorithm’s efficiency, it is combined with a multiresolution scheme that employs different discretizations to deal with displacement, design, and density variables. In addition, the time spent solving the linear equilibrium systems is substantially reduced using multigrid as a preconditioner for the conjugate gradient method. Since multiresolution can lead to the appearance of unwanted artefacts in the structure, we propose an adaptive strategy for increasing the degree of the displacement elements, with a technique for suppressing unnecessary variables that provides accurate solutions with a moderate impact on the algorithm’s performance. We also propose a new thresholding strategy, based on gradient information, to obtain structures composed only by solid or void regions. Computational experiments carried out in Matlab prove that the new algorithm effectively generates high-resolution structures at a low computational cost.
\end{abstract}

% ---------- %

\section{Introduction}

Designing a structure that is both safe and economical requires great skill, experience and creativity from designers or engineers. 
Structural topology optimization is a mathematical methodology that helps in selecting a suitable initial shape for a new structure, facilitating and speeding up its design. 

In the simplest formulation of the problem, we define the domain in which the structure must be built, the external loads, the supports, the amount of material available and, possibly, regions of the domain where material is required or prohibited. The aim of topology optimization is to decide how the material should be distributed in the domain so that the structure has the greatest possible rigidity, supporting the application of external loads without suffering excessive displacements and deformations, remaining in static equilibrium and satisfying a maximum volume constraint. 

Mathematically, we associate to each point in the domain a function that assumes a value of $1$ if the point contains material or $0$ otherwise. However, from a computational point of view, it is extremely difficult and expensive to make each point be represented by just these two possibilities. Therefore, we generally use a continuous function that can take values in the interval $[0,1]$, so that each value represents the density of material at one point in the domain. However, we still need to eliminate these intermediate densities, for which Bends\o e \cite{Bendsoe} suggested using the SIMP model (Solid Isotropic Material with Penalization), in which the density is raised to an exponent greater than 1. In this work, we make use of this approach.
Other ways of dealing with the topology optimization problem include, for example, level set methods \cite{Dijk}, isogeometric analysis \cite{WangY} and moving morphable components \cite{Zhang}.

Since its introduction, topology optimization has gained prominence in academia and industry, with various applications in civil, mechanical, aerospace and automotive engineering, among other areas. Bibliographic reviews of the subject are available at \cite{SigmundMaute,BendsoeSigmund,Deaton,Eschenauer}.
From a computational point of view, the biggest challenge in the area has been the development of efficient algorithms for solving large problems, especially three-dimensional ones. An analysis of the strategies proposed and the challenges encountered in speeding up the topology optimization process are presented in \cite{Mukherjee}. 

In this article, we aim to provide a robust and efficient algorithm for solving three-dimensional problems, combining some existing techniques and proposing new strategies that allow us to obtain better quality solutions. Although we only consider the structural problem in its simple form, the strategies presented here can be applied to more complex problems, such as those involving manufacturing constraints \cite{ZegardPaulino,Zhu2021} or the design of compliant mechanisms \cite{Sigmund1997,Zhu}.  

When the SIMP approach is adopted, the topology optimization problem is converted into a nonlinear programming problem with many variables, so methods that require the calculation of second derivatives are usually avoided. The algorithms most commonly applied to solve the problem are the optimality criteria method \cite[section 1.2]{BendsoeSigmund}, the method of moving asymptotes (MMA) \cite{Svanberg} and sequential linear programming \cite{Sigmund1997,GomesSenne}. In this work, we implement the sequential linear programming algorithm described by Gomes and Senne \cite{GomesSenne}, which is globally convergent, and we adopt a stopping criterion based on the first-order optimality conditions of the problem, similar to the one proposed in \cite{GomesSenne2014}.   

The most expensive step of most topology optimization algorithms is the computation of the objective function, which involves solving a linear system at each iteration. Several strategies have been applied to reduce the cost of this step, such as reanalysis \cite{Senne}, reduced order models \cite{Gogu,Xiao} and multigrid methods \cite{Amir,PeetzElbanna}.

Amir {\it et al.} \cite{Amir} use the geometric multigrid method as a preconditioner of the conjugate gradient method for solving the linear systems. Peetz and Elbanna \cite{PeetzElbanna} showed that the algebraic multigrid preconditioners are more robust than the geometric multigrid preconditioners when combined with the GMRES method. In our implementation, we test both the geometric and the algebraic multigrid method as a preconditioner for the conjugate gradient method. However, our algebraic multigrid version is different from the one used in \cite{PeetzElbanna} and is based on the work of Franceschini {\it et al.} \cite{Franceschini} and Magri {\it et al.} \cite{Paludetto}. Furthermore, to combine multigrid with the higher order finite element method, our geometric multigrid preconditioner is similar to the $h$-multigrid cited in \cite{Sundar}, in which the shape functions of the finite element method are used to construct the prolongation matrix.

In order to make the algorithm more efficient and, at the same time, obtain more detailed and better quality structures, we also adopted the multiresolution technique proposed by Nguyen {\it et al.} \cite{Nguyen2010,Nguyen2012}. Other works using multiresolution include \cite{Groen,Gupta2020,Park,Yoo}. In this technique, meshes with different refinements are used for different steps of the topology optimization process: a coarser mesh for approximating the nodal displacements, an intermediate mesh for solving the optimization problem and a finer mesh to define the distribution of material. As a result, the final resolution of the structure is high, while the computational cost of the most expensive steps of the algorithm is reduced.

The drawback of the multiresolution approach is the fact that the low precision in the approximation of the displacements (that is based on the coarser mesh) can cause inaccuracies in the distribution of material densities (computed using the finer mesh). These inaccuracies are not noticed by the optimization method and end up generating artificially rigid regions in the form of QR-patterns, incomplete bars or ``floating'' parts in the structure, known as artefacts (see Gupta {\it et al.} \cite{Gupta2018}). The appearance of artefacts was also reported in the works of Groen {\it et al.} \cite{Groen}, Gupta {\it et al.} \cite{Gupta2017} and Nguyen {\it et al.} \cite{Nguyen2017}. The results of these studies suggest that using a density filter and increasing the degree of approximation of the displacements can mitigate the appearance of artefacts. 

In our numerical tests, we also observed the appearance of artefacts in the solutions obtained with multiresolution, especially when the radius of the density filter is smaller than 1. To avoid this, we used Lagrange and serendipity finite elements with quadratic or cubic shape functions, in order to improve the approximation of the displacements and, consequently, the quality of the solutions. 

However, using finite elements of a higher degree increases the time taken to solve the linear systems, since the number of nodes in the mesh increases. To alleviate this problem, Gupta {\it et al.} \cite{Gupta2020} proposed an adaptive strategy that modifies the degree of the finite elements and the number of design variables of each element during the optimization process. The idea is to increase these values in regions of the domain that require greater precision, such as those that may contain artefacts, and reduce the values in regions where low precision may suffice. The results obtained with this strategy are promising, although only two-dimensional problems were tested. 

In our work, we propose a new adaptive strategy for increasing the degree of the finite elements, which consists of (i) solving the problem with linear finite elements; (ii) suppressing some variables in regions of the domain filled with solid or void elements; and (iii) solving the problem again with higher order finite elements. This process of fixing variables and solving a problem with a better displacement approximation is repeated until the structure found is adequate, which occurs with quadratic or, at most, cubic elements.

Although the SIMP model leads to solutions with few intermediate densities, the filter used to reduce checkerboard patterns can cause the appearance of intermediate densities in the transition between solid and void, affecting not only the shape of the structure, but also its stiffness. The development of techniques to obtain discrete solutions is therefore a relevant issue in topology optimization \cite{Wang,Xu}. 
In this work, we propose a new thresholding strategy, with the aim of obtaining solutions formed only by solid or void elements. Our strategy is based on the gradient of the Lagrangian  and allows us to obtain good quality solutions that are close to local minimizers. 

This text is organized as follows. In Section \ref{sec: topopt}, we present the structural topology optimization problem and briefly review the ideas behind the multigrid method for solving linear equilibrium systems. In Section \ref{sec: mr}, we explain the multiresolution strategy adopted. The sequential linear programming algorithm is presented in Section \ref{sec: SLP}. In Section \ref{sec: thresholding}, we describe our density projection strategy. The adaptive strategy for increasing the degree of finite elements is presented in Section \ref{sec: fix}. Finally, in Section \ref{sec: results}, the test problems and the computational results obtained with a Matlab implementation of the new algorithm are discussed.

% ---------- %

\section{Structural topology optimization} \label{sec: topopt}

Let us assume that the domain $\Omega \subset \mathds{R}^3$ in which a structure must be contained, its supports, the external loads and an upper limit for the amount of material are known.

After using the finite element method to discretize the domain, we want to determine which elements will contain material. To avoid dealing with an integer optimization problem, we consider that each element has a density $\rho_e \in [0,1]$. 
Moreover, to reduce the occurrence of intermediate densities, we apply the SIMP model \cite{Bendsoe}, in which each density is raised to an exponent greater than $1$. As the stiffness of the finite element is directly related to the properties of the material, this density penalization is introduced in the element stiffness matrix, $\mathrm{K}^{(e)}$, using the formula
\begin{equation} \label{eq: matrixK}
\mathrm{K}^{(e)}(\rho_e) = (E_{min} + {(\rho_e)}^p(E_0 - E_{min}))\, \int_{\Omega_e} \mathrm{B}^T\mathrm{D}\mathrm{B} \; d\Omega_e,
\end{equation}
for $e = 1, 2, ..., n_{el}$, where $n_{el}$ is the total number of elements in the mesh, $\rho_e$ is the density of element $e$, $\mathrm{B}$ is the strain–displacement matrix, $\mathrm{D}$ is the matrix that contains the elastic parameters of the material, $\Omega_e$ is the element domain, $p$ is the penalty parameter of the SIMP model, $E_0$ is the Young's modulus of the solid material, and $E_{min}$ is a very small positive real number, used to prevent the global stiffness matrix $\mathrm{K}(\rho)$ from becoming singular due to the presence of densities very close to zero. 

Assuming that we have a linearly elastic structure, the static equilibrium condition is represented by the linear system $\mathrm{K}(\rho) \mathrm{u} = \mathrm{f}$, where $\mathrm{f}$ is the vector of nodal loads and $\mathrm{u}$ is the nodal displacements vector. Approximating the structure's average compliance by $\mathrm{f}^T\mathrm{u}$, our topology optimization problem may be stated as
\begin{equation} \label{prob: topopt_ori}
	\begin{aligned}
		& \underset{\rho}{\text{Min}} & & \mathrm{f}^T\mathrm{u} \\
		& \text{s.t.} & & \mathrm{K}(\rho) \mathrm{u} = \mathrm{f}, \\
		& & & \sum_{e=1}^{n_{el}} v_e \rho_e \leq V_{max}, \\
		& & & 0 \leq \rho_e \leq 1, \quad e = 1, 2, ..., n_{el},
	\end{aligned}
\end{equation} 
where $n_{el}$ is the total number of elements in the mesh, $V_{max}$ is the maximum volume allowed, and $v_e$ and $\rho_e$ are, respectively, the volume and density of the element $e$.

After imposing the boundary conditions, $\mathrm{K}(\rho)$ becomes positive definite, so we can remove the displacements vector from (\ref{prob: topopt_ori}) writing $\mathrm{u}=~{\mathrm{K}(\rho)}^{-1} \mathrm{f}$, as well as $\mathrm{f}^T\mathrm{u} = \mathrm{f}^T{\mathrm{K}(\rho)}^{-1} \mathrm{f}$. In this case, our nonlinear optimization problem reduces to 
\begin{equation} \label{prob: topopt}
	\begin{aligned}
		& \underset{\rho}{\text{Min}} & & \mathrm{f}^T{\mathrm{K}(\rho)}^{-1} \mathrm{f} \\
		& \text{s.t.} & & \sum_{e=1}^{n_{el}} v_e \rho_e \leq V_{max}, \\
		& & & 0 \leq \rho_e \leq 1, \quad e = 1, 2, ..., n_{el}.
	\end{aligned}
\end{equation} 
It should be clear, however, that the objective function of (\ref{prob: topopt}) is computed solving the system $\mathrm{K}(\rho) \mathrm{u} = \mathrm{f}$, so $\mathrm{K}(\rho)^{-1}$ is never explicitly formed. 

\subsection{Finite element types} \label{sec: elements}

We suppose that the domain $\Omega \subset \mathds{R}^3$ has the shape of a right rectangular prism and we use elements that also have this shape, called rectangular prismatic elements. Each element node has three degrees of freedom, which represent the displacements of the node in the coordinate directions. 

The Lagrange and serendipity finite element families are considered, with shape functions of degree $1$, $2$ and $3$. Linear Lagrange elements have 8 nodes and are very attractive for three-dimensional problems. However, in some cases, the linear approximation of the structure's displacements may not be adequate, so polynomials of higher degree should be used. 
When the 27-node quadratic Lagrange elements are used, the approximation of the displacements becomes more accurate, albeit at a higher computational cost. If even greater precision is desired, it is possible to use a cubic Lagrange element, which has $64$ nodes. However, this element is generally not competitive with the equivalent element from the serendipity family.

Serendipity elements use fewer polynomial terms in the functions that approximate the displacements and, consequently, have fewer nodes ($20$ nodes for quadratic and $32$ nodes for cubic elements). Besides, the degree of approximation obtained is equivalent to that of the corresponding Lagrange element, which is why serendipity elements are widely used in applications. 
For three-dimensional problems, the use of shape functions with a degree greater than $3$ can undermine the efficiency of the optimization algorithm, even for elements from the serendipity family.

\subsection{Density filters} 

Solving the topology optimization problem in discretized form (\ref{prob: topopt}) may present numerical instabilities. Sigmund and Petersson \cite{SigmundPetersson} classify these instabilities into three main categories: local minima, mesh dependence and checkerboard patterns.

The first problem is related to the fact that the objective function has many local minima, so different solutions can be found using different parameters of the algorithm. Mesh dependence refers to the case in which we do not obtain, qualitatively, the same structure for different mesh refinements.

The checkerboard pattern, in turn, is the distribution of the material forming regions where void and solid elements alternate in a periodic fashion. According to Díaz and Sigmund \cite{DiazSigmund}, one of the causes of this pattern, which makes the structure artificially stiffer, is the type of finite element used. 

One way to prevent the occurrence of checkerboard patterns is to apply a mathematical operator, called spatial filter. In this work, we use the weighted average density filter proposed by Bruns and Tortorelli \cite{BrunsTortorelli}, that showed the best results when combined with the SLP algorithm. The idea behind this type of filter is to replace the density of each element $e$ with a weighted average of the densities of the elements in the neighborhood $B(e,r_{min})$, where $r_{min}>0$ is the filter radius. Letting $dist(e,j)$ represent the Euclidean distance between the centers of elements $e$ and $j$, we define the filtered density as 
\begin{equation*}
	\tilde{\rho}_e \equiv \tilde{\rho}_e(\rho) = \sum_{j \in B(e,\, r_{min})} \frac{\omega_{ej}}{\omega_e} \rho_j,
\end{equation*}
where the weighting factor $\omega_{ej}$ is given by
\begin{equation} \label{eq: filter_weights}
	\omega_{ej} = \left. \exp\left(-\frac{dist(e,j)^2}{2\left(r_{min}/3\right)^2} \right)\middle/ 2 \pi \left(r_{min}/3\right) \right.
\end{equation}
if $dist(e,j) \leq r_{min}$ and $0$ otherwise, and $\omega_e$ is the sum of the weight factors of all the neighbors of element $e$.

The original densities are replaced by the filtered densities in the objective function and in the constraints of the topology optimization problem (\ref{prob: topopt}). As a result, the optimality conditions of the problem remain compatible. Other advantages of this filter are the preservation of the linearity of the volume constraint, the low cost of its application when the radius $r_{min}$ is not very large and the fact that the gradients of the objective function and the constraints of the problem can be calculated analytically. 

Our algorithm has a pre-filtering step, which consists of obtaining the neighbors of each element and calculating the weight factors. This step can be carried out just once, in order to make the filter application cheaper.

There are several other types of spatial filters in the literature. Sigmund \cite{Sigmund99} uses a filter that modifies the gradient vector of the objective function, based on a weighted average of the first
 derivatives of the function calculated in a neighborhood of each element. In this case, as the gradient is modified and the objective function is not changed, the optimality conditions of the problem become incompatible, which makes the application of globally convergent optimization methods difficult. Bruns \cite{Bruns} presents a strategy that combines density filtering with an alternative to the SIMP model for penalizing intermediate densities, called the Sinh model, in reference to the use of the hyperbolic sine function. In this case, the penalty is applied to the volume restriction, which becomes nonlinear. An advantage of the Sinh model is obtaining results with few intermediate densities.

\subsection{The static equilibrium system}

Since the global stiffness matrix, $\mathrm{K}$, depends on the densities $\rho$, we need to solve a linear system in the form $\mathrm{K}(\rho) \mathrm{u} = \mathrm{f}$ at each iteration. This is the most time-consuming step in the optimization process, exceeding half of the total time spent by the algorithm (see \cite{Liu,Mukherjee}). It is therefore necessary to adopt efficient methods for solving such linear systems. 
Under the usual boundary conditions, the stiffness matrix is symmetric positive definite, which makes it possible to use the Cholesky factorization to solve the linear system. However, despite being a straightforward and accurate method, the Cholesky factorization can require a lot of time and memory when applied to three-dimensional problems. 

Numerical experiments show that an iterative solving method, such as the conjugate gradient method, is more efficient in cases where the number of elements is very large \cite{Liu}. But for the method to be efficient, it is important to choose a suitable preconditioner. 
Unfortunately, when the usual preconditioners -- such as a diagonal matrix or the incomplete Cholesky factorization -- are adopted, the linear systems still require the largest portion of the time spent solving three-dimensional topology optimization problems. For this reason, in this work, we follow the proposal of Amir {\it et al.} \cite{Amir} and use the multigrid method \cite{Briggs,Trottenberg} as a preconditioner. 

Consider the discretization of the domain $\Omega$ into a mesh $\Omega_h$ and the linear system $\mathrm{A}_h\mathrm{u}_h = \mathrm{f}_h$. Also consider another coarser mesh, $\Omega_H$, in which the system is cheaper to solve. In the multigrid method, we apply a few iterations of a stationary iterative method (Jacobi, Gauss-Seidel, SSOR, etc.) to solve the linear system on the original mesh until the error becomes smooth. Then we switch to a coarser mesh and obtain an approximation to the error by solving the residual system. Finally, we return to the finer mesh and correct the approximation to the solution of the original system.  

What distinguishes the various multigrid implementations is the definition of the coarse grid, as well as the mechanism used to transfer information from one mesh to another. It is a common practice to define $\mathrm{P}$, the prolongation operator from $\Omega_H$ to $\Omega_h$, and to use $\mathrm{R}=\mathrm{P}^T$ as the restriction operator from $\Omega_h$ to $\Omega_H$. A typical two-grid step (cycle) of a multigrid method is given in Algorithm \ref{alg: multigrid}.

\begin{algorithm}[t] 
	\caption{Two-grid cycle of the multigrid method} \label{alg: multigrid}
	\begin{algorithmic}[1]
    \State \hspace{-0.22cm}\begin{tabular}[t]{l}
		          Given the initial approximation $\mathrm{u}_h^{(k)}$, \\
							compute $\mathrm{v}_h$ applying $\eta_1$ iterations of a smoothing method to the system $\mathrm{A}_h\mathrm{v}_h = \mathrm{f}_h$
					 \end{tabular}
		\State Compute the residual $\mathrm{r}_h = \mathrm{f}_h-\mathrm{A}_h\mathrm{v}_h$
		\State Restrict the residual in the coarse grid: $\mathrm{r}_H = \mathrm{R}\mathrm{r}_h$
		\State Solve on $\Omega_H$: $\mathrm{A}_H\mathrm{e}_H = \mathrm{r}_H$
		\State Prolongate the correction: $\mathrm{e}_h = \mathrm{P}\mathrm{e}_H$
		\State Compute the corrected approximation: $\overline{\mathrm{v}}_h = \mathrm{v}_h + \mathrm{e}_h$
    \State \hspace{-0.22cm}\begin{tabular}[t]{l}
		          Using $\overline{\mathrm{v}}_h$ as an initial approximation, \\
							compute $\mathrm{u}_h^{(k+1)}$  applying $\eta_2$ iterations of a smoothing method to the system $\mathrm{A}_h\mathrm{u}_h = \mathrm{f}_h$
					 \end{tabular}
	\end{algorithmic}
\end{algorithm}

The cycle is repeated until some stopping criterion is reached, such as the relative norm of $\mathrm{r}_h$ being less than some pre-established tolerance. We can also apply the multigrid method again to solve the residual system $\mathrm{A}_H\mathrm{e}_H = \mathrm{r}_H$, in a recursive manner, giving rise to cycles with more meshes, such as the well known V and W cycles (see \cite{Trottenberg}). 
In general, the smoothing method is just a classical iterative method such as Jacobi or SSOR, and $\eta_1$ and $\eta_2$ are small. Furthermore, it is also not necessary to solve the residual system exactly.

As stated earlier, constructing the coarse mesh and obtaining the prolongation operator depends on the version of the method being used. Our implementation includes both a geometric and an algebraic version of the multigrid. 
In our geometric multigrid (GMG) algorithm, we consider that each element of the coarse mesh contains $8$ elements of the fine mesh, so that each edge of the element in the coarse mesh is twice the size of the corresponding edge in the fine mesh. 

The extension to a fine mesh node is done through interpolation, using the surrounding coarse mesh nodes as interpolating points. The prolongation matrix $\mathrm{P}$, that can be easily described for a mesh with linear elements, becomes more complicated if we want to use higher degree elements. 
For Lagrange-type elements, we can still use linear interpolation. For serendipity elements, we propose constructing the prolongation operator using interpolation via the shape functions of the finite element method. Letting $\xi_i$, $\eta_i$ and $\zeta_i$ be the local coordinates of a node $i$ and supposing that we know the values $\mathrm{v}_H(\xi_i, \eta_i, \zeta_i)$ at the nodes of the coarse mesh element, the value $\mathrm{v}_h(\xi,\eta,\zeta)$ for a fine mesh node contained in this element is given by
\begin{equation*}
	\mathrm{v}_h(\xi,\eta,\zeta) = \sum_{i=1}^{n_e} \phi_i(\xi,\eta,\zeta) \mathrm{v}_H(\xi_i, \eta_i, \zeta_i),
\end{equation*}
where $n_e$ is the number of nodes of the element in the coarse mesh and $\phi_i$, for $i=1,2,...,n_e$, are the element shape functions. In the case of a quadratic serendipity element, for example, the extension is done by quadratic interpolation, using the $20$ nodes of the element as interpolating points.  

The algebraic multigrid (AMG) is a more robust method, as it can be applied to a wider variety of problems, including those in which the mesh is irregular or unstructured. In this case, only the algebraic equations of the linear system are used to construct the components of the method. 
In our algebraic multigrid implementation, the coarse mesh is generated using the classic process described by St\"{u}ben \cite{Stuben}, but with the connections between nodes calculated following a more recent alternative method, proposed by Magri {\it et al.} \cite{Paludetto} and improved by Franceschini {\it et al.} \cite{Franceschini}, which promises to be more efficient in the case of structural problems. In this alternative approach, the strength of connection between the nodes of the algebraic mesh are obtained based on a test space, constructed using a method based on the minimization of the Rayleigh quotient by conjugate gradients. The prolongation operator is defined using the so called \textit{Dynamic Pattern Least Squares method}, in which the weights of the prolongation are obtained as the solution of a least squares problem. 

% ---------- %

\section{Multiresolution} \label{sec: mr} 

The \textit{resolution} of the structure, i.e. the level of detail it has, depends on the number of elements in the mesh, in a similar way to the resolution of an image, that is usually measured in pixels. 
Obtaining more detailed structures with smoother contours requires a large number of elements in the mesh. Increasing the number of elements, however, causes the number of variables, and therefore the size of the linear equilibrium systems, to grow excessively, especially in the three-dimensional case, resulting in a high consumption of memory and computational time. Most strategies to reduce this time/memory consumption are focused on the solution of the equilibrium system, and include reanalysis techniques \cite{Senne}, order reduction \cite{Gogu,Xiao} or the use of multigrid methods \cite{Amir,PeetzElbanna}.

In this section, we present a multiresolution strategy based on the work of Nguyen {\it et al.} \cite{Nguyen2010,Nguyen2012}. The idea is to work with different resolutions at each step of topology optimization: a coarser mesh for solving the linear systems, an intermediate mesh for solving the optimization problem and a finer mesh on which the distribution of material is defined. As a result, the time spent on the more expensive steps is reduced, while the final resolution of the structure is high.

\subsection{Mesh definition}

Let's consider three meshes for the topology optimization problem:
\begin{enumerate}
	\item the {\it displacement mesh}, that is used to obtain an approximation for the displacements by solving the linear equilibrium system $\mathrm{K}\mathrm{u} = \mathrm{f}$;
	\item the {\it design variable mesh}, where we solve the optimization problem; and
	\item the {\it density mesh}, that is used to represent the material distribution.
\end{enumerate} 
An element belonging to the first mesh is called a \textit{displacement element}, while an element from the third mesh is denoted a \textit{density element}. In the second mesh we have the so called \textit{design variables}, numerical values which may or may not have a physical meaning. 
In the work by Nguyen {\it et al.} \cite{Nguyen2010}, the design variables are defined as the values of the material densities at the centers of the density elements. However, it is not necessary for the design and density meshes to coincide. In \cite{Nguyen2012}, the authors show that the time spent solving the problems can be reduced using three different meshes. In our work, we test both cases.

Taking the displacement mesh as the base, since it is the coarsest one, we define the density mesh dividing each displacement element into ${(n_{mr})}^3$ equal elements. Therefore, each edge of a density element corresponds to $1/n_{mr}$ of the edge of a displacement element. Similarly, we divide each displacement element into ${(d_{mr})}^3$ equal elements, which will form the design variable mesh. Note that $1 < d_{mr} \leq n_{mr}$ and both parameters, $d_{mr}$ and $n_{mr}$, must be integers.

\begin{figure}[h]
  \centering
	\begin{subfigure}{0.3\columnwidth}
		\centering
		\includegraphics[width=0.55\textwidth]{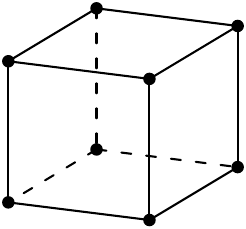}
		\caption{Displacement element.}
		\label{fig: mr_elemA}
	\end{subfigure}%
	\begin{subfigure}{0.3\columnwidth}
		\centering
		\includegraphics[width=0.55\textwidth]{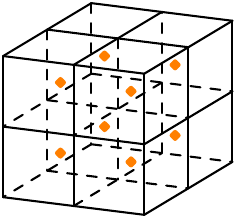}
		\caption{Design variables.}
		\label{fig: mr_elemB}
	\end{subfigure}%
	\begin{subfigure}{0.3\columnwidth}
		\centering
		\includegraphics[width=0.55\textwidth]{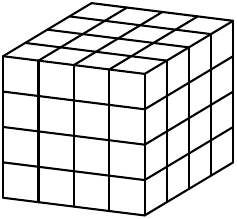}
		\caption{Density elements.}
		\label{fig: mr_elemC}
	\end{subfigure}	
	\caption{Elements in multiresolution meshes, with $n_{mr} = 4$ and $d_{mr} = 2$.}
	\label{fig: mr_elem}
	\vspace{-0.8cm}
\end{figure}

For example, if $d_{mr} = 2$ and $n_{mr} = 4$, each displacement element ``contains'' $8$ design variables and $64$ density elements, as illustrated in Figure \ref{fig: mr_elem}. In Figure \ref{fig: mr_elemA} the black dots represent the displacement element nodes, while in Figure \ref{fig: mr_elemB} the orange dots are the centers of the elements in the design variable mesh. Figure \ref{fig: mr_elemC} shows the density elements. 

\subsection{Projection of design variables to density elements}

At each iteration of our algorithm, after solving an optimization problem in the mesh of design variables, we need to update the density mesh elements.
To make it clear how the densities are calculated, let us denote by $\mathrm{x}$ the vector of design variables and by $\rho$ the vector of material densities. If $d_{mr} = n_{mr}$, we consider that each density is equivalent to a design variable, so we have $\rho = \mathrm{x}$. In the case where $d_{mr} \neq n_{mr}$, we define $\rho = f_p(\mathrm{x})$ where $f_p$ denotes a projection function.

\begin{figure}[h]
	\centering
	\includegraphics[width=0.6\linewidth]{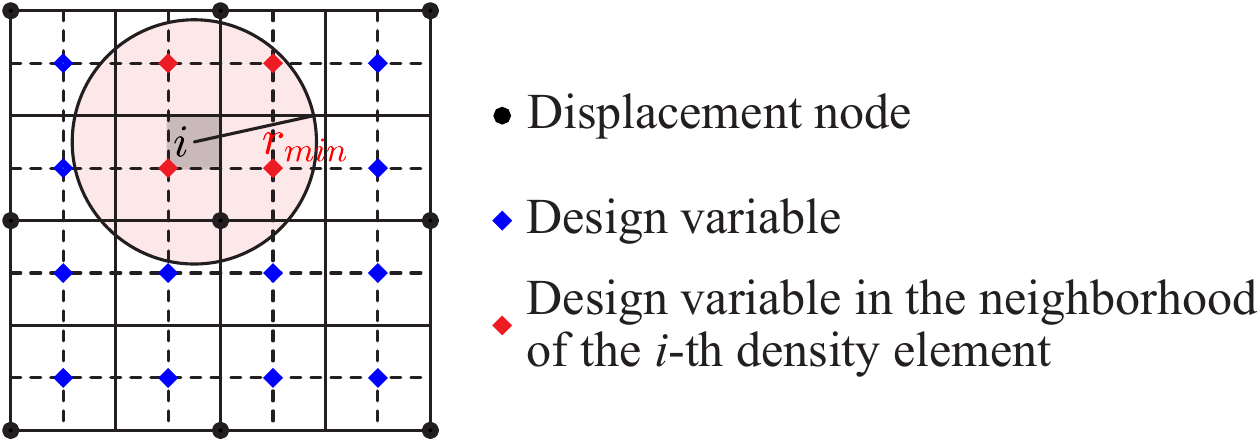}
	\caption{Neighborhood of a two-dimensional density element for $n_{mr} = 4$ and $d_{mr} = 2$.}
	\label{fig: mr_projection}
\end{figure}

We have chosen to use a linear projection, which works in a similar way to applying the density filter. Given a radius $r_{min} > 0$, for each density element $i$, the density $\rho_i$ is calculated as a weighted average of the values of the design variables in a neighborhood $B(i,r_{min})$ of the element: 

\begin{equation} \label{eq: projdsgn}
	\rho_i = \sum_{j \in B(i,\,r_{min})} \frac{\omega_{ij}}{\omega_i} \, x_j,
\end{equation}
where $\omega_{ij}$ is a weight factor and $\omega_i = \sum_{j \in B(i,\,r_{min})} \omega_{ij}$. 
Figure \ref{fig: mr_projection} illustrates the neighborhood of a density element in a two-dimensional problem. 

We kept the weight factors used when applying the density filter, so $\omega_{ij}$ is calculated using (\ref{eq: filter_weights}), where now $dist(i,j)$ denotes the Euclidean distance between the centers of the density element $i$ and element $j$ in the design variable mesh. Note that, when the design variables are equivalent to the densities, the projection is equivalent to applying the filter. 

Considering the distinction between the vector of design variables $\mathrm{x}$ and the vector of densities $\rho$, the topology optimization problem can be rewritten as follows:

\begin{equation} \label{probMR}
	\begin{aligned}
		& \underset{\mathrm{x}}{\text{Min}} & & \mathrm{f}^T {\mathrm{K}(\rho)}^{-1} \mathrm{f} \\
		& \text{s.t.} & & \rho = f_p(\mathrm{x}), \\
		& & & \sum_{i=1}^{n_{\rho}} v_i \rho_i \leq V_{max}, \\
		& & & 0 \leq x_j \leq 1, \;\; j = 1, 2, ..., n_{\mathrm{x}}, 
	\end{aligned}
\end{equation}
where $n_{\mathrm{x}}$ is the number of design variables and $n_{\rho}$ is the number of density elements. If $n_{el}$ is the number of displacement elements, then $n_{\mathrm{x}} = {(d_{mr})}^3 n_{el}$ and $n_{\rho} = {(n_{mr})}^3 n_{el}$.

\subsection{Assembly of the global stiffness matrix} 

The global stiffness matrix $\mathrm{K}(\rho)$ is formed by superimposing the stiffness matrices of the displacement mesh elements. Traditionally, the stiffness matrix $\mathrm{K}^{(e)}(\rho_e)$ of an element $e$ is given by (\ref{eq: matrixK}), where $\rho_e$ is the density of the element.
In the multiresolution approach, each displacement element consists of $N_n = {(n_{mr})}^3$ density elements, as shown in Figure \ref{fig: mr_elemC} for $n_{mr} = 4$, so we need a new formula for computing $\mathrm{K}^{(e)}$.

Let's consider a rectangular prismatic (displacement) element with length $2\ell$, height $2h$ and width $2w$. Furthermore, let us adopt the local coordinate system $(\xi, \eta, \zeta)$ in the interval $[-1,1]$, with the origin at the barycenter of the element, so that $\xi = x/\ell$, $\eta = y/h$ and $\zeta = z/w$. In this case, we have  $\Omega_e = [-1,1]^3$ and $d \Omega_e = dx \; dy \; dz = \ell h w \; d\xi \; d\eta \; d\zeta$. Thus, the integral that appears in (\ref{eq: matrixK}) is given by
\begin{equation*}
	\int_{\Omega_e} \mathrm{B}^T\mathrm{D}\mathrm{B} \; d\Omega_e = \int_{-1}^{1} \int_{-1}^{1} \int_{-1}^{1} \ell h w \; \mathrm{B}^T\mathrm{D}\mathrm{B} \;d\xi\;d\eta\;d\zeta.
\end{equation*} 

When we divide the element into $N_n = {(n_{mr})}^3$ density elements, we also divide the interval $[-1,1]$ into small subintervals of size $1/n_{mr}$. Therefore, the integral can be calculated as a sum of integrals in smaller intervals:
\begin{equation*}
	\int_{\Omega_e} \mathrm{B}^T\mathrm{D}\mathrm{B} \; d\Omega_e =
	\sum_{i=1}^{N_n} \left( \int_{\Omega_{e_i}} \mathrm{B}^T\mathrm{D}\mathrm{B} \; d\Omega_{e} \right) =
	\sum_{i=1}^{N_n} \mathrm{K}^{(e_i)},
\end{equation*}
where $\Omega_{e_i}$ is the domain of a density element $e_i$ contained in the displacement element $e$. 

Let $[a_{\xi}, b_{\xi}]$, $[a_{\eta}, b_{\eta}]$ and $[a_{\zeta}, b_{\zeta}]$ be the intervals of integration of the density element $e_i$ in the $\xi$, $\eta$ and $\zeta$ directions, respectively, so that $\Omega_{e_i} = [a_{\xi}, b_{\xi}]\times[a_{\eta}, b_{\eta}]\times[a_{\zeta},b_{\zeta}]$. The integration intervals are different for each density element $e_i$, but here we omit the index $e_i$ to simplify the notation. In this case, we have
\begin{equation*}
	\mathrm{K}^{(e_i)} = \int_{a_{\zeta}}^{b_{\zeta}} \int_{a_{\eta}}^{b_{\eta}} \int_{a_{\xi}}^{b_{\xi}} \ell h w \; \mathrm{B}^T\mathrm{D}\mathrm{B} \;d\xi\;d\eta\;d\zeta.
\end{equation*}

If the displacement element is divided into $N_n$ density elements, we need to compute $N_n$ matrices $\mathrm{K}^{(e_i)}$, one for each density element $e_i$. Fortunately, if the meshes are regular, these $N_n$ matrices need to be computed once, as they are repeated for each displacement element.
Finally, the stiffness matrix $\mathrm{K}^{(e)}(\rho)$ of a displacement element is defined as the linear combination of the matrices $\mathrm{K}^{(e_i)}$, weighted by the penalty terms of the SIMP method. That is, 
\begin{equation} \label{eq: matrixK_mr}
	\mathrm{K}^{(e)}(\rho) = \sum_{i=1}^{N_n} (E_{min} + {(\rho_{e_i})}^p(E_0 - E_{min})) \; \mathrm{K}^{(e_i)}.
\end{equation} 
The global stiffness matrix, $\mathrm{K}$, is assembled in the same way as in the traditional topology optimization method, by superimposing the matrices $\mathrm{K}^{(e)}$ for all the mesh elements.

\subsection{Gradients}

To apply optimization methods that involve first derivatives, such as the sequential linear programming method, we need to compute the gradients of the objective function and also of the constraints of the problem.

According to (\ref{eq: projdsgn}), for each density element, the density $\rho_i$ is a weighted average of the values of the design variables in its neighborhood. Therefore, the derivative of the volume constraint in (\ref{probMR}) with respect to $x_k$ may be written as 
\begin{equation*}
	\frac{\partial c(\rho)}{\partial x_k} = \sum_{i \in \bar{B}(k,\, r_{min})} v_i \frac{\omega_{ik}}{\omega_i},
\end{equation*}
where $c(\rho) = \sum_{i=1}^{n_{\rho}} v_i \rho_i - V_{max}$, and $\bar{B}(k,r_{min})$ denotes the set of density elements whose centers are at a distance less than or equal to $r_{min}$ from the center of element $k$ in the design variable mesh. If we consider that all of the elements have a fixed volume, it is possible to compute the Jacobian matrix only once. 

Rewriting the objective function of the problem as 
\begin{equation*}
  f(\rho) = \mathrm{f}^{T}\mathrm{u}(\rho) = \mathrm{u}(\rho)^T\mathrm{K}(\rho)\mathrm{u}(\rho),
\end{equation*}
its first derivatives with respect to the design variables are given by
\begin{equation*}
	\frac{\partial f(\rho)}{\partial x_k} = - \sum_{i \in \bar{B}(k,r_{min})} \mathrm{u}^{T}\left[p(\rho_i)^{(p-1)}(E_0 - E_{min})  \mathrm{K}_i \right] \mathrm{u}\, \frac{\omega_{ik}}{\omega_i},
\end{equation*}
where $\mathrm{K}_i$ is a matrix with the same dimension as $\mathrm{K}$, but formed only by the components of $\mathrm{K}^{(e_i)}$, the stiffness matrix of the density element $i$. In practice, instead of expanding the matrix $\mathrm{K}^{(e_i)}$ to form $\mathrm{K}_i$, we can just extract from the vector $\mathrm{u}$ the displacements of the nodes corresponding to the element $e$ that contains the density element $i$. 

\subsection{Occurrence of artefacts} \label{sec: artefacts} 

\begin{figure}[h]
	\centering
	\includegraphics[width=0.45\linewidth]{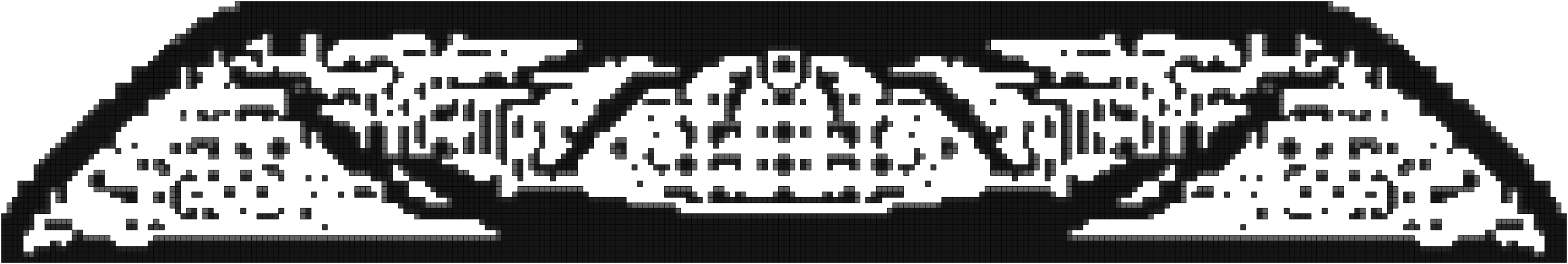}
	\caption{Two-dimensional structure with artefacts.}
	\label{fig: artefacts}
\end{figure}

Unfortunately, one drawback of the multiresolution methods is the occurrence of regions with disconnected or loosely connected material parts, as illustrated in Figure \ref{fig: artefacts}. Gupta {\it et al.} \cite{Gupta2018} studied the emergence of these regions, which they called \textit{numerical artefacts} or \textit{QR-patterns} (due to their pattern in two-dimensional problems).

QR-patterns are artificially stiff regions that occur when the approximation of the displacements is done with shape functions composed by low order polynomials. As the density mesh is finer than the displacement mesh, the errors increase as we pass the information from one mesh to the other, making the artefacts artificially rigid. The occurrence of artefacts was also reported in \cite{Groen,Nguyen2017}. 

The presence of artefacts may also result from the non-uniqueness of the design field. Studying the limitations of multiresolution, Gupta {\it et al.} \cite{Gupta2017}  conclude that, in order to preserve the uniqueness of the design field, we must satisfy $m \leq \mathrm{rank}(\mathrm{K}^{(e)})$, where $\mathrm{K}^{(e)}$ is the element stiffness matrix and $m$ is the number of design variables used to define the material distribution within each finite element.

For the three-dimensional case, $\mathrm{rank}(\mathrm{K}^{(e)}) = 3n_e - 6$, where $n_e$ is the number of nodes in the element. In addition, in this work we consider that the number of design variables per element is given by $m = (d_{mr})^3$. Table \ref{tab: bounds} contains the upper bounds for $m$ and $d_{mr}$ for Lagrange and serendipity rectangular prismatic elements. It should be noted, however, that the quality of the results also depends on the radius of the filter and on the parameter $n_{mr}$, so it is possible to obtain reasonable results even for $d_{mr}$ slightly above the specified limit. 
\begin{table}[h]
\centering
\caption{Maximum number of design points for some rectangular prismatic finite elements.}
\def\arraystretch{1.2}
\setlength{\tabcolsep}{8pt}
\small
\begin{tabular}{c|ccc|cc}
	\hline 
	Element type & \multicolumn{3}{c|}{Lagrange} & \multicolumn{2}{c}{Serendipity} \\ \cline{2-6}
	and degree & 1 & 2 & 3 & 2 & 3 \\ \hline 
	$m$ max. & $18$ & $75$ & $186$ & $54$ & $90$ \\
	$d_{mr}$ max. & $2$ & $4$ & $5$ & $3$ & $4$ \\ \hline 
\end{tabular}
\label{tab: bounds}
\end{table}
%\begin{table}[h]
%\centering
%\caption{Maximum number of design points for some rectangular prismatic finite elements.}
%\def\arraystretch{1.2}
%\setlength{\tabcolsep}{10pt}
%\small
%\begin{tabular}{c|c|c}
	%\hline 
	%{\bf Element type} & {\bf $\bm{m}$ max.} & {\bf $\bm{d_{mr}}$ max.} \\ \hline 
	%Lagrange, degree 1 & $18$ & $2$ \\
	%Lagrange, degree 2 & $75$ & $4$ \\
	%Lagrange, degree 3 & $186$ & $5$ \\
	%Serendipity, degree 2 & $54$ & $3$ \\
	%Serendipity, degree 3 & $90$ & $4$ \\ \hline 
%\end{tabular}
%\label{tab: bounds}
%\end{table}

The detection of artefacts is still an open problem and, usually, we only notice them once the final solution has been obtained. Furthermore, the artefacts may not be as apparent as shown in Figure \ref{fig: artefacts}. 
Strategies currently used to mitigate the appearance of artefacts include increasing the degree of the finite elements and applying a density filter. As in the case of checkerboard patterns, larger filter radii produce a more homogeneous material distribution, preventing the formation of artefacts.

On the other hand, a large filter radius may reduce the resolution of the structure, which is one of the purposes of multiresolution. In order to use smaller radii and thus obtain more detailed structures, it is necessary to increase the degree of the displacement elements, which, in turn, increases the computational cost. For three-dimensional problems, it is necessary to find a balance between filter radius size and the degree of the displacement elements, which must not exceed 3. 

In this work, we obtained good results using multiresolution with linear elements and the application of a filter with a moderate radius size. However, increasing the number of density elements did not lead to an increase in the resolution of the structure. Therefore, to obtain more detailed structures without the presence of artefacts, we combined a small filter radius with an adaptive strategy of increasing the degree of the elements, suppressing some variables in order to increase the efficiency of the algorithm. This strategy will be presented in Section \ref{sec: fix}.

% ---------- %

\section{Sequential linear programming} \label{sec: SLP}

The topology optimization problem can be solved using various nonlinear programming methods known in the literature, the most common being the optimality criteria method \cite{Andreassen,Liu,FerrariSigmund,Amir} and the MMA \cite{Svanberg}, which are easy to implement. When these methods are adopted, the stopping criterion usually does not depend on the optimality conditions of problem (\ref{prob: topopt}), but only on the step size or the number of iterations.

However, to solve large and more challenging topology optimization problems, it is convenient to use a more robust optimization method. In order to obtain better quality solutions, in this work we use the globally convergent sequential linear programming algorithm proposed by Gomes and Senne \cite{GomesSenne}, which has proved to be very efficient for solving topology optimization problems. 

\subsection{Description of the method}

Let us consider a general nonlinear optimization problem written in the form
\begin{equation} 
	\begin{aligned}
		& \underset{\mathrm{x}}{\text{Min}} & & f(\mathrm{x}) \\
		& \text{s.t.} & & c(\mathrm{x}) = 0, \\
		& & & \mathrm{x}_l \leq \mathrm{x} \leq \mathrm{x}_u,
	\end{aligned} \label{prob: P} 
\end{equation} 
where the objective function, $f: \mathds{R}^n \rightarrow \mathds{R}$, and the constraint functions, $c: \mathds{R}^n \rightarrow \mathds{R}^m$, have Lipschitz continuous first derivatives, $\mathrm{x} \in \mathds{R}^n$ contains the design variables and $\mathrm{x}_l, \mathrm{x}_u \in \mathds{R}^n$ define the lower and upper bounds for the components of $\mathrm{x}$.  

The linear approximations for $f$ and $c$ at $\mathrm{x}^{(k)} \in \mathds{R}^n$ are given by
\begin{equation*}
	f\left(\mathrm{x}^{(k)}+\mathrm{s}\right) \approx f\left(\mathrm{x}^{(k)}\right) + \nabla f\left(\mathrm{x}^{(k)}\right)^T\mathrm{s} \ \ \ \mbox{and} \ \ \ 
	c\left(\mathrm{x}^{(k)}+\mathrm{s}\right) \approx c\left(\mathrm{x}^{(k)}\right) + \mathrm{A}\left(\mathrm{x}^{(k)}\right)\mathrm{s}, 
\end{equation*}
where $\smash{\nabla f(\mathrm{x}^{(k)})}$ and $\smash{\mathrm{A}(\mathrm{x}^{(k)}) = [ \nabla c_1 (\mathrm{x}^{(k)}) \cdots \nabla c_m (\mathrm{x}^{(k)}) ]^{T}}$ 
denote the gradient of the objective function and the Jacobian matrix of the constraints, respectively. Therefore, given a point $\mathrm{x}^{(k)} \in \mathds{R}^n$, problem (\ref{prob: P}) can be approximated by the linear programming (LP) problem
\begin{equation} 
	\begin{aligned}
		& \underset{\mathrm{s}}{\text{Min}} & & \nabla f\left(\mathrm{x}^{(k)}\right)^T\mathrm{s} \\
		& \text{s.t.} & & \mathrm{A}\left(\mathrm{x}^{(k)}\right)\mathrm{s} + c\left(\mathrm{x}^{(k)}\right) = 0, \\
		& & & \mathrm{x}_l \leq \mathrm{x}^{(k)}+\mathrm{s} \leq \mathrm{x}_u.
	\end{aligned} \label{prob: Ptil} 
\end{equation} 

SLP is an iterative method that solves a sequence of linear problems in the form (\ref{prob: Ptil}). At each iteration $k$ of the algorithm, a previously computed point $\mathrm{x}^{(k)}$ is used to generate the linear programming problem, for which we obtain a solution $\tilde{\mathrm{s}}$. If this solution is considered acceptable, the approximate solution for (\ref{prob: P}) is then updated according to $\mathrm{x}^{(k+1)} = \mathrm{x}^{(k)} + \tilde{\mathrm{s}}$. Otherwise, problem (\ref{prob: Ptil}) is solved again with more stringent bounds $\mathrm{x}_l$ and $\mathrm{x}_u$, so the step $\mathrm{s}$ is reduced. This process is repeated until some stopping criterion is satisfied.

\subsubsection{Trust regions}

To control the length of the step $\mathrm{s}$, we introduce a \textit{trust region}, in which we believe that the linear approximations of the functions are good enough. Thus, we require that $\mathrm{s}$ also satisfies $\left\| \mathrm{s} \right\|_{\infty} \leq \delta,$ where $\delta > 0$ is the radius of the trust region. Combining the trust region with the constraints in (\ref{prob: Ptil}),  we rewrite the  problem as
\begin{equation} 
	\begin{aligned}
		& \underset{\mathrm{s}}{\text{Min}} & & \nabla f\left(\mathrm{x}^{(k)}\right)^T\mathrm{s} \\
		& \text{s.t.} & & \mathrm{A}\left(\mathrm{x}^{(k)}\right)\mathrm{s} + c\left(\mathrm{x}^{(k)}\right) = 0, \\
		& & & \mathrm{s}_l \leq \mathrm{s} \leq \mathrm{s}_u,
	\end{aligned} \label{prob: Pstar} 
\end{equation} 
where $(\mathrm{s}_l)_i = \max\{ -\delta, ( \mathrm{x}_l)_i - \mathrm{x}_i^{(k)} \}$ and $(\mathrm{s}_u)_i = \min \{ \delta, (\mathrm{x}_u)_i - \mathrm{x}_i^{(k)} \}$, for $i = 1,2,...,n$.

\subsubsection{Feasibility of the linear programming problem} 

The box constraints in (\ref{prob: Pstar}) ensure that the problem is bounded. Even so, it can be infeasible, unless $\mathrm{x}^{(k)}$ is feasible for problem (\ref{prob: P}). To obtain a feasible solution, we need to solve the problem
\begin{equation} 
	\begin{aligned}
		& \underset{(\mathrm{s},\, \mathrm{z})}{\text{Min}} & & \mathrm{e}^T\mathrm{z} \\
		& \text{s.t.} & & \mathrm{A}\left(\mathrm{x}^{(k)}\right)\mathrm{s} + c\left(\mathrm{x}^{(k)}\right) + \mathrm{E}\left(\mathrm{x}^{(k)}\right)\mathrm{z} = 0, \\
		& & & \bar{\mathrm{s}}_l \leq \mathrm{s} \leq \bar{\mathrm{s}}_u, \\
		& & & \mathrm{z} \geq 0, 
	\end{aligned} \label{prob: Pstar2} 
\end{equation}
where $\smash{\mathrm{e}^T = \left[ 1 \ 1 \ \cdots \ 1 \right] \in \mathds{R}^q}$, $\mathrm{z} \in \mathds{R}^q$ is the vector of artificial slack variables corresponding to the infeasible constraints of (\ref{prob: P}) at $\mathrm{x}^{(k)}$, $\smash{\mathrm{E}(\mathrm{x}^{(k)}) \in \mathds{R}^{m\times q}}$ is a matrix formed by columns of the identity matrix $\mathrm{I} \in \mathds{R}^{m \times m}$, or $-\mathrm{I}$, and
\begin{equation*}
	\begin{aligned}
		(\bar{\mathrm{s}}_l)_i &= \max\left\{ -0.8\,\delta; \, {( \mathrm{x}_l)}_i - \mathrm{x}_i^{(k)} \right\}, \qquad		(\bar{\mathrm{s}}_u)_i &= \max\left\{ 0.8\,\delta; \, {( \mathrm{x}_u )}_i - \mathrm{x}_i^{(k)} \right\},
	\end{aligned}
\end{equation*}
for $i = 1,2,...,n$. 
Problem (\ref{prob: Pstar2}) is similar to the artificial problem solved to obtain a feasible solution for the Simplex method. Note that the trust region of this problem has been slightly reduced, to allow problem (\ref{prob: Pstar}) to have a sufficiently large feasible region, so the objective function can be reduced.

When the topology optimization problem has linear constraints (as in (\ref{prob: topopt})) and we choose a feasible starting point $\mathrm{x}^{(0)}$, the remaining points $\mathrm{x}^{(k)}$ generally remain feasible. Therefore, we always try to solve the problem (\ref{prob: Pstar}) to obtain the step $\tilde{\mathrm{s}}$. If the LP solver detects that it has no solution, then we compute $\mathrm{s}_n$, the solution of the artificial problem (\ref{prob: Pstar2}), and define $\tilde{\mathrm{s}} = \mathrm{s}_n$. Fortunately, in our experiments, we were always able to obtain $\tilde{\mathrm{s}}$ by solving (\ref{prob: Pstar}) directly.

\subsubsection{Merit function} 

After obtaining the step $\tilde{\mathrm{s}}$, we need to check if $\mathrm{x}^{(k+1)} = \mathrm{x}^{(k)} + \tilde{ \mathrm{s}}$ is a better approximation than $\mathrm{x}^{(k)}$ to the solution of the original problem (\ref{prob: P}).  We want the new approximation to sufficiently reduce the value of the objective function, but we also need to control the infeasibility of the problem.   
One way of dealing with this duplicity of objectives is to define a merit function $\psi: \mathds{R}^n \times \mathds{R} \rightarrow \mathds{R}$ that combines the two measures (optimality and feasibility). In our SLP method, we use 
\begin{equation*}
	\psi(\mathrm{x},\theta) = \theta f(\mathrm{x}) + (1-\theta)\varphi(\mathrm{x}),
\end{equation*}
where $\varphi(\mathrm{x}) = \frac{1}{2} \left\| c(\mathrm{x}) \right\|_2^2$ is a measure of the infeasibility at $\mathrm{x} \in \mathds{R}^n$ with respect to problem (\ref{prob: P}) and $\theta \in (0,1]$ is a penalty parameter whose purpose is to establish a balance between optimality and feasibility. 

The step acceptance is based on a comparison between the merit function reduction predicted by the linear model (\ref{prob: Pstar}) and the actual reduction obtained for the original problem (\ref{prob: P}). The actual reduction of $\psi$ is defined by
\begin{equation*} 
	A_{red} = \theta \left( f(\mathrm{x}^{(k)}) - f(\mathrm{x}^{(k)} + \tilde{ \mathrm{s}}) \right) + (1-\theta)\left( \varphi(\mathrm{x}^{(k)}) - \varphi(\mathrm{x}^{(k)} + \tilde{ \mathrm{s}}) \right),
\end{equation*}
while the predicted reduction of $\psi$ is given by
\begin{equation*}
	P_{red} = - \theta \, \nabla f\left(\mathrm{x}^{(k)}\right)^T\!\tilde{\mathrm{s}} + (1-\theta)\left( M\left(\mathrm{x}^{(k)},0\right) - M\left(\mathrm{x}^{(k)},\tilde{\mathrm{s}}\right)\right),
\end{equation*}
where $M(\mathrm{x},\mathrm{s}) = \frac{1}{2} \left\| \mathrm{A}(\mathrm{x})\mathrm{s} + c(\mathrm{x}) \right\|_2^2$ is an approximation to $\varphi$.

If the actual reduction of the merit function is less than $10\%$ of the predicted reduction, i.e. if $A_{red} < 0.1 P_{red}$, then the step $\tilde{\mathrm{s}}$ is rejected and we decrease the radius of the trust region, taking $\delta_{k+1} = \min\left\{ 0.25\left\|\tilde{\mathrm{s}}\right\|_{\infty}; \, 0.1\,\delta_k \right\}$.

On the other hand, the step $\tilde{\mathrm{s}}$ is accepted if $A_{red} \geq 0.1 P_{red}$. Furthermore, if $A_{red} \geq 0.5 P_{red}$, i.e. the reduction of the merit function is good enough, we increase the radius of the trust region by defining $\delta_{k+1} = \max \left\{\min\left\{ 2.0 \,\delta_k; \, \left\|\mathrm{x}_u-\mathrm{x}_l\right\|_{\infty} \right\}, \delta_{min}\right\}$,
where $\delta_{min}$ is a small positive real number. The parameters used to accept the step and to increase or decrease the trust region radius were obtained experimentally. Here we have kept the numerical values used in the implementation to simplify the notation. 

The penalty parameter $\theta$ is also updated at each iteration. 
To do this, starting with $\theta_0 = \theta_{max} = 1$, at each iteration $k$ we calculate $\theta_k = \min\left\{\theta_k^{large}, \theta_k^{sup}, \theta_{max} \right\}$, where
\begin{align} 
	\theta_k^{large} &= \left[1 + \frac{N}{(k+1)^{1.1}}\right] \min \left\{\theta_0,...,\theta_{k-1}\right\}, \label{eq: theta1} \\
	\theta_k^{sup} &=
	\left\{ \begin{array}{ll}
		0.5 \left( \dfrac{P_{red}^{fsb}}{P_{red}^{fsb} - P_{red}^{opt}} \right) & \mbox{if} \;\; P_{red}^{opt} \leq 0.5 P_{red}^{fsb}, \\
		1 & \mbox{otherwise}, 
	\end{array} \right. \label{eq: theta2}
\end{align} 
with $P_{red}^{fsb} = M (\mathrm{x}^{(k)},0 ) - M (\mathrm{x}^{(k)},\tilde{\mathrm{s}} )$ and $P_{red}^{opt} = \nabla f(\mathrm{x}^{(k)})^T\!\tilde{\mathrm{s}}$.
When the step is rejected, we set $\theta_{max} = \theta_k$, and when the step is accepted, we take $\theta_{max} = 1$. The constant $N \geq 0$, used to calculate $\theta_k^{large}$, must be adjusted to allow $\theta$ to have a non-monotone decrease over the iterations. In this work, we consider $N = 10^6$.

\subsection{Algorithm}

Given a starting point $\mathrm{x}^{(0)} \in \mathds{R}^n$ satisfying $\mathrm{x}_l \leq \mathrm{x}^{(0)} \leq \mathrm{x}_u$, an initial trust region radius $\delta_0 \geq \delta_{min} > 0$, $\theta_0 = \theta_{max} = 1$ and $k = 0$, Algorithm \ref{alg: SLP} solves the nonlinear optimization problem (\ref{prob: P}) using the SLP method.

\begin{algorithm}[h] 
	\caption{Sequential Linear Programming (SLP)} \label{alg: SLP}
	\begin{algorithmic}[1]
		
		\While{none of the stopping criteria are satisfied}
		
		\State find $\tilde{\mathrm{s}}$, the solution of problem (\ref{prob: Pstar})
		
		\If{$\tilde{\mathrm{s}}$ is infeasible}
		
		\State find $\mathrm{s}_n$, the solution of problem (\ref{prob: Pstar2})
		
		\State $\tilde{\mathrm{s}} \longleftarrow \mathrm{s}_n$
		
		\EndIf
		
		\State update $\theta_k = \min\left\{\theta_k^{large}, \theta_k^{sup}, \theta_{max} \right\}$ using (\ref{eq: theta1}) and (\ref{eq: theta2})
		
		\State compute $A_{red}$ and $P_{red}$
		
		\If{$A_{red} < 0.1 P_{red}$} \hfill {\it (reject the step)} \
		
		\State $\delta_{k} \longleftarrow \min\left\{ 0.25\left\|\tilde{\mathrm{s}}\right\|_{\infty}, \, 0.1\,\delta_k \right\}$
		
		\State $\theta_{max} \longleftarrow \theta_k$
		
		\Else \hfill {\it(accept the step)} \
		
		\State $\mathrm{x}^{(k)} \longleftarrow \mathrm{x}^{(k)} + \tilde{\mathrm{s}}$
		
		\If{$A_{red} \geq 0.5 P_{red}$}
		
		\State $\delta_{k} \longleftarrow \min\left\{ 2.0\,\delta_k, \, \left\|\mathrm{x}_u-\mathrm{x}_l\right\|_{\infty} \right\}$
		
		\EndIf
		
		\State $\delta_{k} \longleftarrow \max\left\{\delta_{k}, \delta_{min}\right\}$
		
		\State update $\mathrm{A}(\mathrm{x}^{(k)})$ and $\nabla f(\mathrm{x}^{(k)})$
		
		\State $\theta_{max} \longleftarrow 1$
		
		\State $k \longleftarrow k+1$
		
		\EndIf
		
		\EndWhile
		
	\end{algorithmic}
\end{algorithm}

Note that an outer iteration starts only when step $\tilde{\mathrm{s}}$ is accepted. However, even when the step is rejected, we have to compute the objective function, which involves solving a linear system to find the displacements. The convergence results of the algorithm are presented in \cite{GomesSenne}.

\subsubsection{Stopping criteria} \label{sec: stopping}

One important feature of our algorithm is the fact that we adopt a stopping criterion based on the Karush-Kuhn-Tucker (KKT) conditions of problem (\ref{prob: P}), ensuring that the solution is close to a stationary point for this problem. 

Let $P_X(\mathrm{y})$ be the orthogonal projection of $\mathrm{y}$ onto the set  $X = \{\mathrm{x} \in \mathds{R}^n \| \mathrm{x}_l \leq \mathrm{x} \leq \mathrm{x}_u \}$ and let $\nabla \mathcal{L} (\mathrm{x}^{(k)},\mathrm{\lambda}^{(k)})$ be the gradient of the Lagrangian function associated to the equality constraints, where $\mathrm{\lambda}^{(k)}$ is the corresponding vector of Lagrange multipliers (that is the dual solution of the LP problem (\ref{prob: Pstar})). At each iteration $k$ we compute the vector $ g_{P}(\mathrm{x}^{(k)}) = P_{X}(\mathrm{x}^{(k)} - \nabla \mathcal{L} (\mathrm{x}^{(k)},\mathrm{\lambda}^{(k)})) - \mathrm{x}^{(k)}$, and require the solution $\mathrm{x}^{(k)}$ to satisfy $\smash{\|g_{P}(\mathrm{x}^{(k)})\|_{\infty} < \varepsilon_g}$. 

As mentioned, it is a common practice in topology optimization to stop the algorithm whenever the step is small enough, i.e. $\left\|\tilde{\mathrm{s}}\right\|_{\infty} < \varepsilon_s$, where $\varepsilon_s > 0$ is a tolerance, or when the decrease in the objective function is small, i.e. $| f(\mathrm{x}^{(k)}) - f(\mathrm{x}^{(k-1)}) | < \varepsilon_f$. 

In our algorithm, these three criteria are combined, and we stop the algorithm when the tolerances $\varepsilon_f$ and $\varepsilon_g$ are reached at the same time for three consecutive iterations or when the $\varepsilon_s$ tolerance is reached for three consecutive iterations. 

Furthermore, we set a limit of $500$ iterations as an additional stopping criterion. In practice, using $\varepsilon_s = 10^{-4}$, $\varepsilon_f = 5\times10^{-2}$ and $\varepsilon_g = 10^{-3}$, the criterion based on the projected gradient together with the decrease of the objective function was the first to be reached in all of our numerical experiments.

% ---------- %

\section{Density projection strategy} \label{sec: thresholding}

Let $V$ be the total domain volume and $v_{frac}$ be the predetermined volume fraction that the structure may occupy. Supposing that the domain has been divided into a mesh with $n_{\rho}$ density elements, the volume constraint of the structural topology optimization (\ref{probMR}) problem can be written as
\begin{equation} \label{eq: rvol} 
	\sum_{i=1}^{n_{\rho}} v_i \rho_i \leq v_{frac} \, V, 
\end{equation} 

Although the SIMP model approach usually leads to discrete topologies, the density filter application can cause intermediate densities to occur more frequently, typically forming transition layers between solid and void regions. In addition to raising manufacturing concerns, intermediate densities affect the objective function value, so the development of techniques to obtain fully discrete solutions is  relevant. 

The goal of the density projection strategy, or simply the thresholding strategy, is to eliminate intermediate densities from the solution vector, $\rho$, without violating the volume constraint (\ref{eq: rvol}). Thus, after applying the strategy, we expect the vector of projected densities, $\tilde{\rho}$, to have exactly $\lfloor v_{frac}\,n_{\rho} \rfloor$ elements with density $1$ and the remaining elements to have density $0$.

We adopt two thresholding strategies, one that is based on the density values and another that takes into account the gradient of the Lagrangian of the  problem, $\nabla \mathcal{L}$. Regardless of the strategy used, the first step is to round the variables that are very close to $0$ or $1$. In our tests, we round $\rho_i$ to $1$ if $\rho_i \geq 0.95$ and round it to $0$ if $\rho_i \leq 0.05$. 

\subsection{Strategy based on the density values}

The first thresholding strategy adopted by the algorithm consists of performing the following steps: 
\begin{enumerate}[1)]
	\item Arrange the density values in descending order; 
	\item Round the $\lfloor v_{frac}\,n_{\rho} \rfloor$ higher densities to $1$; 
	\item Assign $0$ to all of the remaining densities.  
\end{enumerate}
After this attempt, we verify if the vector $d = \tilde{\rho} - \rho$ is a descent direction to the Lagrangian, that is, we verify the angle condition
\begin{equation} \label{eq: angle}
	- \frac{{\nabla \mathcal{L}}^T d}{\|\nabla \mathcal{L}\| \|d\|} > \cos(\vartheta_{max}), 
\end{equation} 
where $\vartheta_{max}$ is the largest angle allowed between the vectors $-\nabla \mathcal{L}$ and $d$ (we adopt $\vartheta_{max} = 89.9^{\circ}$).
If condition (\ref{eq: angle}) is satisfied, we accept the vector $\tilde{\rho}$ as the new density vector. Otherwise, we resort to the strategy based on the gradient of the Lagrangian. 

\subsection{Strategy based on the gradient of the Lagrangian} 

The second thresholding strategy consists of projecting the vector $\hat{\rho}(\alpha) = \rho - \alpha \nabla \mathcal{L}$ onto the box $0 \leq \hat{\rho}(\alpha) \leq 1$. In this case, we need to determine the highest value of $\alpha$, say $\tilde{\alpha}$, such that the projected solution, $\tilde{\rho} = \rho - \tilde{\alpha} \nabla \mathcal{L}$, satisfies the angle condition (\ref{eq: angle}), as well as the two safeguards defined bellow: 
\begin{itemize}
	\item $\tilde{\rho}_i$ can only be rounded to $1$ if it satisfies $\tilde{\rho}_i \geq v^u_{min}$. 
	\item $\tilde{\rho}_i$ can only be rounded to $0$ if it satisfies $\tilde{\rho}_i \leq v^l_{max}$. 
\end{itemize}
In our experiments, we adopt the thresholds $v^u_{min} = 0.3$ and $v^l_{max} = 0.7$. These safeguards prevent the densities close to $0$ from being rounded to $1$ and vice versa, which may cause algorithm instabilities due to abrupt changes in the solution.

Let $\nabla_i \mathcal{L}$ be the component of the gradient of the Lagrangian related to the density $\rho_i$. If $\nabla_i \mathcal{L}>0$, then $\rho_i - \alpha \nabla_i \mathcal{L} < \rho_i$ for $\alpha >0$, indicating that the density should be rounded to $0$. On the other hand, if $\nabla_i \mathcal{L}<0$, the density  should be rounded to $1$.

In practice, we determine, for each component $i$, the value $\alpha_i$ that makes $\rho_i - \alpha_i \nabla_i \mathcal{L} = 0$ (if $\nabla_i \mathcal{L} > 0$) or $\rho_i - \alpha_i \nabla_i \mathcal{L} = 1$ (if $\nabla_i \mathcal{L}< 0$).  After that, we arrange the values $\alpha_i$, $i = 1,2,...,n_{\rho}$, in ascending order and define $\tilde{\alpha}$ as the highest $\alpha_i$ such that $\tilde{\rho} =~\rho - \tilde{\alpha} \nabla \mathcal{L}$ satisfies the safeguards defined above and the angle condition (\ref{eq: angle}). 

\subsection{Quality control of the thresholded solution}

If the thresholding strategy based on the density values is adopted and the angle condition is satisfied, the projected solution vector $\tilde{\rho}$ is discrete. On the other hand, if the strategy based on the gradient of the Lagrangian is used, there is no guarantee that all densities are rounded. Some densities change value but remain intermediate, while those associate to zero Lagrangian gradient entries do not change at all.
Furthermore, the thresholded solution, $\tilde{\rho}$, may not correspond to a local minimizer of the topology optimization problem and may not even satisfy the volume constraint (\ref{eq: rvol}). 

To control the quality of the thresholded solution obtained, we apply again the sequential linear programming (SLP) algorithm, using $\tilde{\rho}$ as the initial density vector. This process of solving the problem and projecting the densities is repeated until one of the following criteria is satisfied. 
\begin{itemize}
	\item The difference between two consecutive thresholded solutions is relatively small and the volume constraint is satisfied within a prescribed tolerance, that is, 
	\begin{equation*}
		\| \tilde{\rho}_k - \tilde{\rho}_{k-1} \|_1 < \varepsilon_N \| \tilde{\rho}_{k-1} \|_1 \ \ \ \mbox{ and } \ \ \ \sum_{i=1}^{n_{\rho}} v_i \rho_i \leq v_{frac} \, V + \varepsilon_V.
	\end{equation*}
	\item The limit number of thresholding attempts is reached. 
\end{itemize}
In our experiments, we define $\varepsilon_N = 0.01$ and $\varepsilon_V = 0.005$, and we limit the number of thresholding attempts to $10$. 

In order to accelerate the reduction of intermediate densities, some parameters of the SLP method may be changed between the first and the remaining calls of the algorithm. A simple but effective idea is to reduce the filter radius when it is much larger than one. In our algorithm, just before the second call to the SLP method, we set $r_{min}^{prj} = \min\left(r_{min}, 1.1\right)$. Another promising but risky idea is to use only linear elements from the second call of the SLP algorithm. This strategy is effective when the densities are very close to $0$ or $1$ in the original solution, but may be dangerous otherwise, so we decided not to adopt it, leaving its analysis for a future study.

\subsection{Projection based on the Heaviside filter}

When the gradient of the Lagrangian is used to project the densities, there may be a large number of intermediate variables and there is no warranty that the volume constraint will be preserved when rounding them to 0 or 1. 

The aim of the smooth Heaviside filter is to bring closer to $1$ the densities greater than a value $\eta \in (0,1)$, as well as to bring closer to $0$ the densities smaller than $\eta$. To cope with the volume constraint, the value $\eta$ must be selected in such a way that $\lfloor v_{frac}\, n_{\rho} \rfloor$ densities approach $1$.
In this work, we use the version of the Heaviside filter based on the hyperbolic tangent, proposed by Wang {\it et al.} \cite{Wang}. In this case, the projected density of the $i$-th element is given by
\begin{equation*}
	h(\rho_i, \beta, \eta) = \frac{\tanh(\beta\eta) + \tanh(\beta(\rho_i-\eta))}{\tanh(\beta\eta) + \tanh(\beta(1-\eta))},
\end{equation*}
where $\rho_i$ is the original density of the element and $\beta$ is a projection parameter that indicates how bold the thresholding strategy will be. 
We choose the value of $\eta$ using the strategy proposed by Xu {\it et al.} \cite{Xu}. Given the parameter $\beta$ and supposing that $\rho$ is the vector of original densities, the value of $\eta$ that preserves the current volume corresponds to the zero of the function $f$ given by $f(\eta) = \sum_{i=1}^{n_{\rho}} v_i h(\rho_i, \beta, \eta) - \sum_{i=1}^{n_{\rho}} v_i \rho_i$.

We remark that $h(\rho_i,\beta,0) \geq \rho_i$, so that $f(0) > 0$ if there is at least one element with density greater than $0$. Analogously, $h(\rho_i,\beta,1) \leq \rho_i$, which implies that $f(1) < 0$. Therefore, by the intermediate value theorem, there is always a value $\eta^{*} \in (0,1)$ such that $f(\eta^{*}) = 0$. We compute $\eta^{*}$ applying Newton's method.

In our implementation, the Heaviside projection procedure is called immediately before and immediately after applying the thresolding strategies. The aim of the first call is to push the densities closer to $0$ or $1$, while the second call is made to preserve the desired volume. 
To speed up the convergence of the densities to $0$ or $1$, the parameter $\beta$ is updated after each density projection attempt using the formula $\beta_{k+1} \longleftarrow \min( \delta_{\beta}\, \beta_k, \, \beta_{max})$, where $\beta_{max}$ is the maximum value admitted for $\beta$. In our tests, we start from $\beta_0 = 1.0$ and adopt $\delta_{\beta} = 2.0$ and $\beta_{max} = 100$.

% ---------- %

\section{Adaptive strategy to increase the degree of the elements} \label{sec: fix}

As stated in Subsection \ref{sec: artefacts}, the solutions obtained with the multiresolution technique can present regions with artificial stiffness, forming undesired artefacts in the structures. The best ways to mitigate this drawback are to increase the degree of the finite elements or to apply a density filter. 

In this work, we try to use small filter radii to obtain more detailed structures and we propose an adaptive strategy to increase the degree of the finite elements in order to avoid artefacts in the structure. This strategy comprises the following steps: 
\begin{enumerate}[1)]
	\item Obtain the solution using elements with degree $1$; 
	\item Choose that variables that can be fixed and remove them from the problem;
	\item Solve the problem again using elements with degree $2$. If necessary, repeat this step using elements with degree $3$ or higher. 
\end{enumerate} 

Initially, we solve the problem using linear elements. If the filter radius is small, the material distribution obtained is likely to contain artefacts. However, this solution already has some relevant information about the topology of the structure, so that it can be a good initial guess for the solution of the problem using elements with a higher degree. Since the initial guess do not need to be accurate, we accelerate the adaptive process relaxing the tolerance used in the stopping criteria of the SLP method, described in Subsection \ref{sec: stopping}. Instead of using $\varepsilon_g$, we adopt $10\varepsilon_g$ as the threshold for the norm of the projected gradient. The original tolerance is used only when we solve the problem with finite elements of the highest degree.

Based on the solid and void regions of the structure obtained with linear elements, we set some variables to 0 or 1 in step (2), so they can be removed from the problem. This procedure, that will be explained in detail in the next subsection, may be repeated in step (3).

Since we use the solution obtained in step (1) as an initial guess and remove some variables, we hope to reduce the time required to obtain the solution in step (3). Our aim is to eliminate the artefacts present in the initial solution in a cheaper way than solving the problem with quadratic elements from the beginning. Nonetheless, if the new solution is still inaccurate, we can solve the problem again increasing the degree of the elements to $3$. This process of fixing variables and solving the problem with a better displacement approximation is repeated until the structure found is adequate. 

After the adaptive process, we apply the density projection strategy presented in Section \ref{sec: thresholding} to obtain the final structure. 

\subsection{Choice of the variables that will be fixed}

\begin{figure}[h]
	\begin{subfigure}{0.5\linewidth}
	  \centering
	  \includegraphics[width=0.85\linewidth]{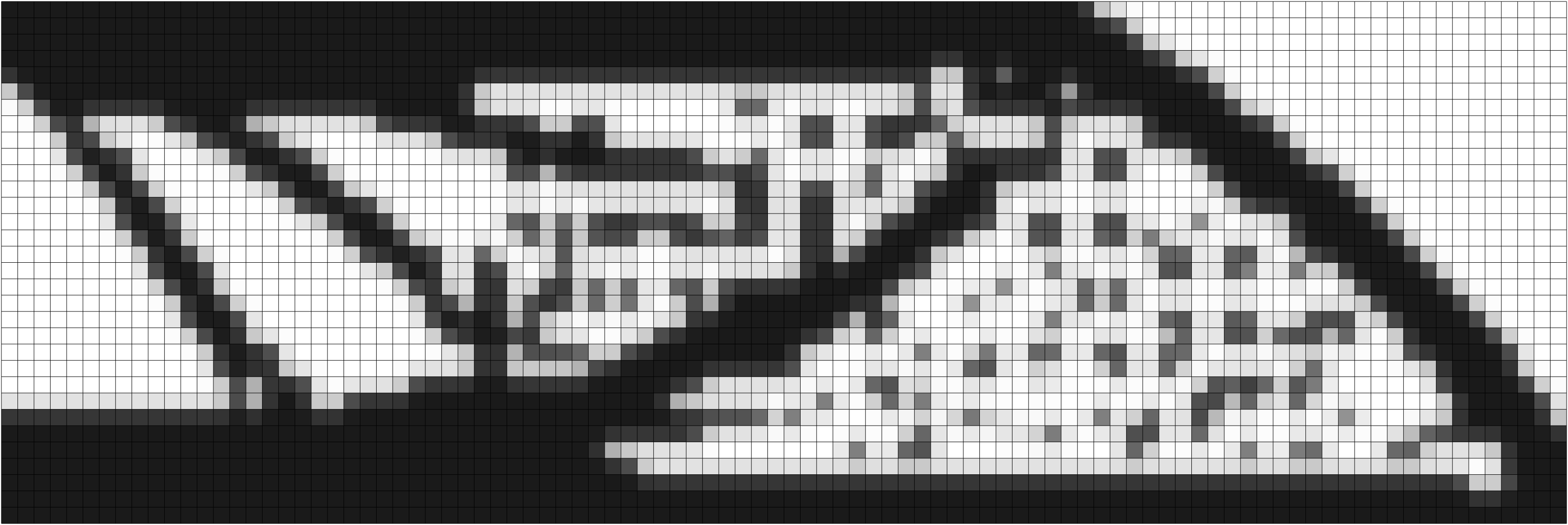}	
	  \caption{Structure obtained \\ with linear finite elements.\centering}
	  \label{fig: mbb2Ddeg1}
	\end{subfigure}
	\begin{subfigure}{0.5\linewidth}
	  \centering
	  \includegraphics[width=0.85\linewidth]{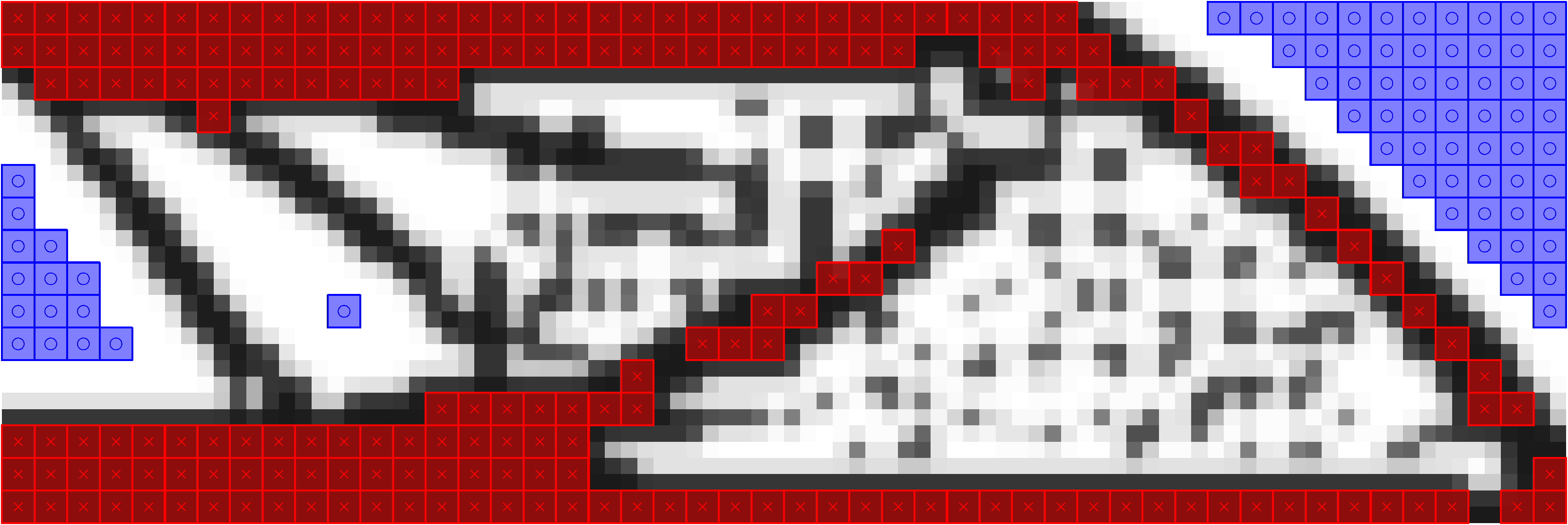}		
    \caption{Structure with void displacement elements highlighted \\in blue and solid displacement elements highlighted in red.\centering}
	  \label{fig: mbb2Dfullelem}
	\end{subfigure}
	\begin{subfigure}{0.5\linewidth}
  	\centering
	  \includegraphics[width=0.85\linewidth]{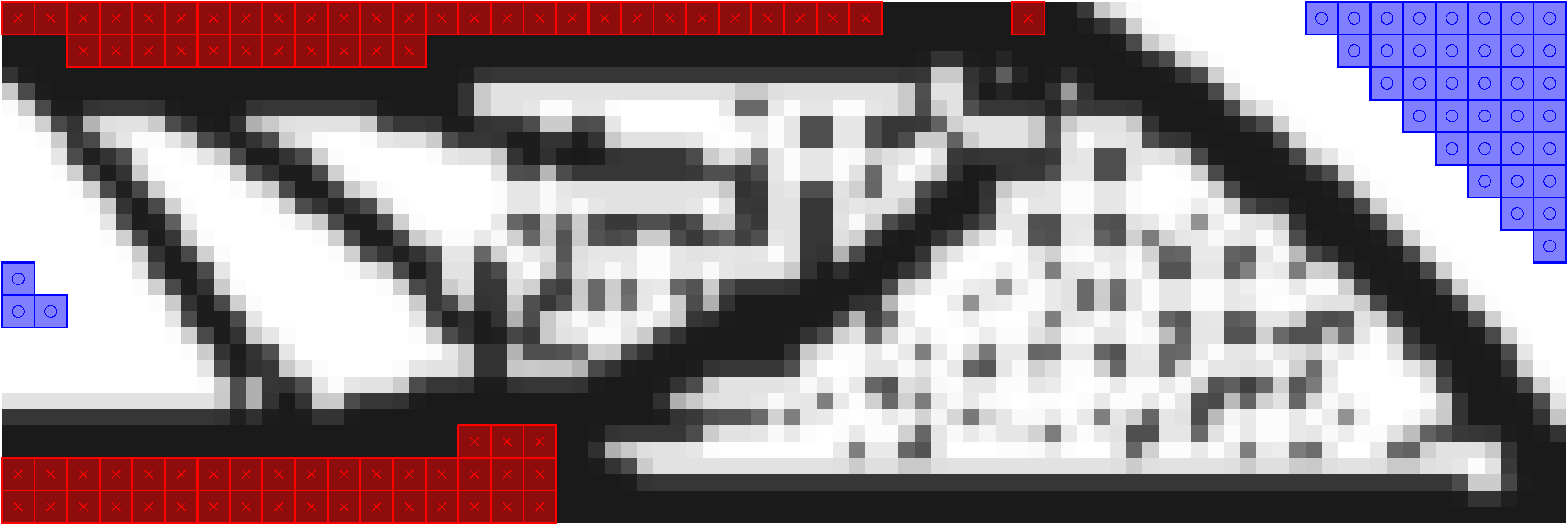}		
	  \caption{Structure with void elements surrounded by void \\ elements highlighted in blue and solid elements \\ surrounded by solid elements highlighted in red.\hspace{0.4cm} \centering}
	  \label{fig: mbb2Dfixelem}
	\end{subfigure}
	\begin{subfigure}{0.5\linewidth}
	  \centering
	  \includegraphics[width=0.85\linewidth]{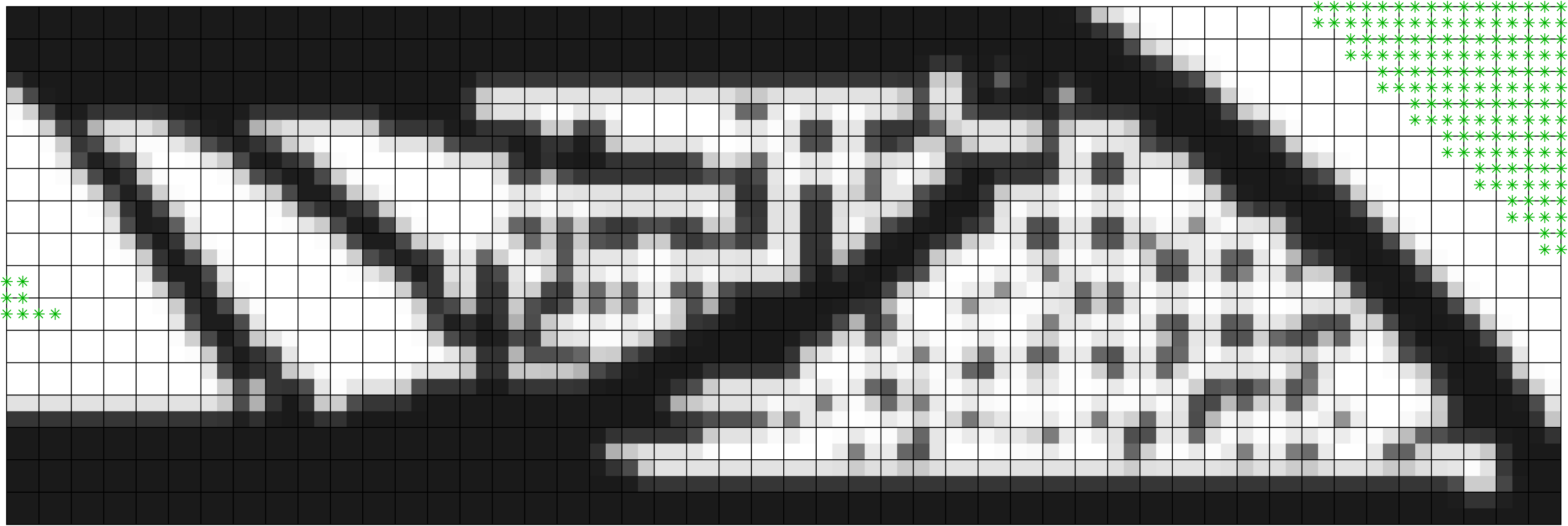}		
	  \caption{Structure with nodes whose degrees \\ of freedom will be suppressed highlighted \\ in green.\ \hspace{0.4cm} \centering}
	  \label{fig: mbb2Dfixnodes}
	\end{subfigure}
	\caption{Two-dimensional MBB beam obtained using multiresolution with $n_{mr} = d_{mr} = 2$ and $r_{min} = 0.8$.}
	\label{fig: mbb2Dx4}
\end{figure}

With the aim of reducing the time spent to solve the problem using elements with degree greater than $1$, we use the solution obtained in step (1) to remove some variables from the problem. To make it easier to understand the proposed idea, we consider a two-dimensional MBB beam example (see Section \ref{sec: results}). 
Figure \ref{fig: mbb2Ddeg1} shows the initial structure obtained using the multiresolution with $n_{mr} = d_{mr} = 2$, $r_{min} = 0.8$ and linear elements. In this figure, we display the entire density mesh, so that void elements are shown in white. We notice that the structure contains many artefacts.

In our multiresolution scheme for three-dimensional problems, we consider that each displacement element is composed by $N_n = (n_{mr})^3$ density elements and contains $N_d = (d_{mr})^3$ design variables. We suppose that a density element $e_i$ contained in the displacement element $e$ has density $\rho_{e_i}$ and we denote by $\nabla_{e_j}f$ the component of the objective function gradient vector corresponding to the design variable $j$ contained in the displacement element $e$. We say that a displacement element is {\it void} if it satisfies the conditions
\begin{equation*}
	\rho_{e_i} \leq \rho_{fix}^l, \mbox{ for all }i = 1,2,...,N_n, \ \ \ \mbox{and} \ \ \ 
	\nabla_{e_j}f \geq -\varepsilon_{grad}, \mbox{ for all }j = 1,2,...,N_d, 
\end{equation*} 
where $\rho_{fix}^l$ and $\varepsilon_{grad}$ are positive constants close to zero. In other words, an element is considered void if all densities inside it are close to $0$ and if the design variables in this element do not increase too much when moving in the direction of $-\nabla f$. 

Similarly, a displacement element is {\it solid} if it satisfies the conditions 
\begin{equation*}
	\rho_{e_i} \geq \rho_{fix}^u, \mbox{ for all }i = 1,2,...,N_n, \ \ \ \mbox{and} \ \ \ 
	\nabla_{e_j}f \leq \varepsilon_{grad}, \mbox{ for all }j = 1,2,...,N_d, 
\end{equation*} 
where $\rho_{fix}^u$ is a real constant close to $1$. Therefore, an element is called solid if all densities inside it are close to $1$ and if the design variables in this element do not decrease too much when moving in the direction of $-\nabla f$.

The safeguards related to the objective function gradient guarantee that, at least in the beginning of the optimization process, the design variables related to void elements will not increase and the design variables inside solid elements will not decrease, so we can fix the values of these variables. 
In our experiments, we use $\rho_{fix}^l = 10^{-6}$, $\rho_{fix}^u = 0.9$ and $\varepsilon_{grad} = 10^{-6}$. We noticed that the threshold $\rho_{fix}^l$, used to find densities close to $0$, must be small to avoid numerical errors, because we also remove from the gradient vector those components related to void elements. 

In our example, each displacement element contains $N_n = 2^2 = 4$ density elements. Moreover, since $n_{mr} = d_{mr}$, the densities correspond to the design variables. Figure \ref{fig: mbb2Dfullelem} shows the structure with the void displacement elements highlighted in blue with a ``$\circ$'' marker and the solid elements highlighted in red with a ``$\times$'' marker. 

Specially at the final iterations of the topology optimization process, the change in the material distribution usually occurs in the boundary of the structure, so it is reasonable to remove from the problem those variables that belong to completely void or solid elements, keeping only the variables inside elements that have mixed or intermediate densities. 
To ensure that the SLP algorithm will have enough space to find a better solution, i.e. to avoid removing from the problem some variables that would otherwise be modified by the algorithm, in our implementation we only fix the design variables that belong to void displacement elements surrounded by void elements or to solid displacement elements surrounded by solid elements. Figure \ref{fig: mbb2Dfixelem} shows these elements for our two-dimensional example.

The components of vector $\mathrm{s}$ corresponding to fixed design variables can be suppressed from the linear programming subproblems solved in steps 2 and 4 of the SLP method (see Algorithm \ref{alg: SLP}). As an alternative, we can assign zero to both the lower and upper bounds of these components, so they will be zero at the solution and the design variables will not be changed.
Note, however, that the densities are obtained from the design variables, so it is important to consider the variables inside solid displacement elements when computing the densities around these elements using equation (\ref{eq: projdsgn}). In practice, since we adopt a conservative strategy for fixing the variables and use a filter radius smaller than the size of the displacement element, the values of the densities inside the fixed elements do not vary much. 

In addition to assigning the value $0$ to the components of the gradient vector corresponding to fixed void elements, we can set to zero the displacements at the internal nodes of the void elements, since the contribution of these elements to the global stiffness matrix, defined in equation (\ref{eq: matrixK_mr}), is negligible. On the other hand, we keep in the problem the nodal displacements of the solid elements. 
Figure \ref{fig: mbb2Dfixnodes} highlights in green with a ``$*$'' marker the nodes whose degrees of freedom will be suppressed from the problem, for our 2D MBB problem with quadratic Lagrange elements. We note that these nodes must be chosen after increasing the degree of the elements. 

In practice, we remove the rows and columns corresponding to the degrees of freedom of these nodes from the global stiffness matrix, and the related components from the nodal load vector, reducing the size of the linear equilibrium systems. We also set the nodal displacements to zero, in a similar way to that adopted for the nodes that are restricted by supports. 
However, it is important to ensure that the entries of the global stiffness matrix related to void element nodes kept in the problem are very small, so there is no abrupt change in the displacements between a node that was removed from the problem and an adjacent node kept in the matrix. Since these entries have the same order of magnitude of $E_{min}$ (the Young's modulus of the void), to avoid numerical errors, we must use a value $E_{min}$ close to zero, but still large enough so that the matrix is not singular. In our tests, we have obtained good results using $E_{min} = 10^{-9}E_0$, being $E_0$ the Young's modulus of the solid material.

If we are using the conjugate gradient method with the multigrid as preconditioner to solve the linear systems, we also need to redo the setup phase of the multigrid method after removing the rows and columns from the global stiffness matrix.

\begin{figure}[h]
	\centering
	\includegraphics[width=0.425\linewidth]{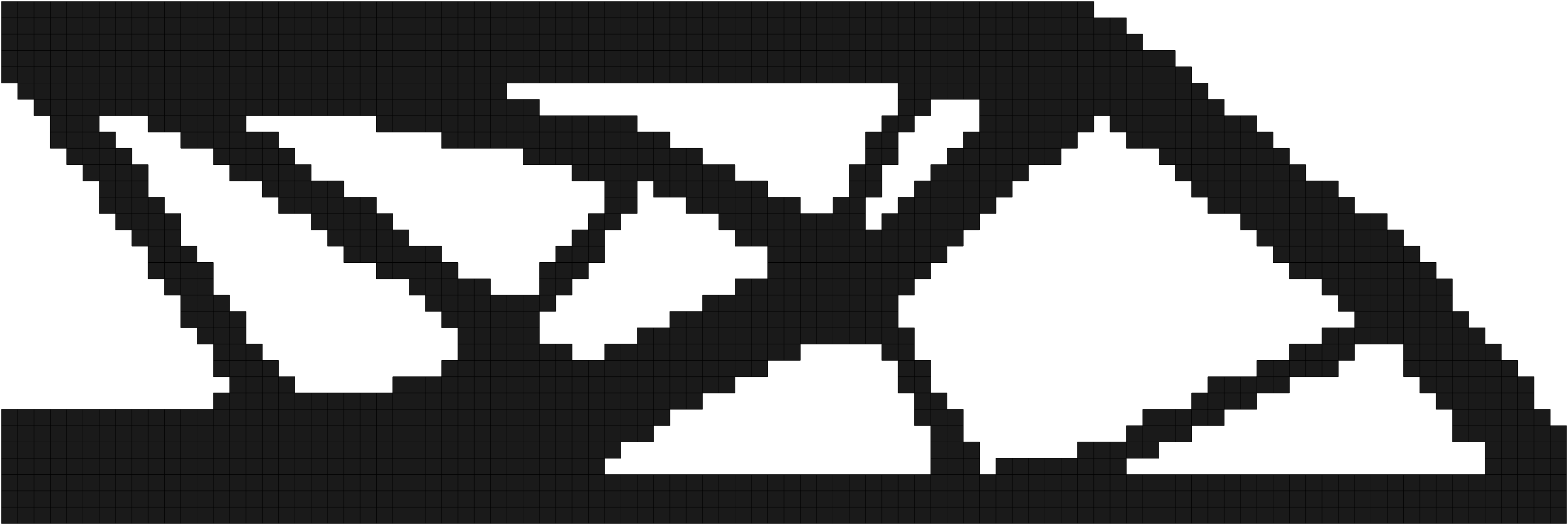}		
	\caption{Final structure obtained after applying the adaptive strategy with quadratic Lagrange elements and projecting the densities.}
	\label{fig: mbb2Ddeg2}
	\vspace{-0.7cm}
\end{figure}

Figure \ref{fig: mbb2Ddeg2} presents the final structure obtained for the two-dimensional example, after applying the adaptive strategy with quadratic Lagrange elements and using the thresholding strategy to project the densities. We notice that this structure does not have visible artefacts, which suggests that increasing the degree of the elements to $2$ was enough to improve the quality of the solution. 

\subsection{Fixed variables update} \label{sec: adapt_fix}

In the process described above, we choose the variables that will be fixed only once, after solving the problem using linear elements. However, during the topology optimization process, new void or solid regions may emerge in the structure, so the fixed variables should change throughout the iterations. 

In fact, choosing the variables that will be fixed based only on the solution obtained with linear elements can be dangerous. For instance, consider the case in which the structure contains an incomplete member due to the presence of a void region in its middle. If we set to zero the design variables in this region, then the algorithm may not be able to converge to a solution that contains this member. 

Taking this into consideration, we propose some strategies to define the variables that should be fixed throughout the optimization of the problem with elements of higher degree. 
\begin{itemize}
	\item {\bf E1:} choose the variables to be fixed once, right after the solution of the problem with linear elements. 
	\item {\bf E2:} choose the variables to be fixed in all SLP iterations for problems with elements of degree greater than one. 
	\item {\bf E3:} choose the variables to be fixed at every $n$ SLP iterations.
	\item {\bf E4:} choose the variables twice. After the solution of the problem with linear elements and at the beginning of the $n$-th SLP iteration.
\end{itemize}

The fixed variables are only updated in outer iterations of the SLP algorithm. In the computational results, we compare each of these strategies, in order to determine which is the safest and more efficient.

% ---------- %

\section{Computational results} \label{sec: results}

In this section, we present the numerical experiments carried out with the algorithm previously described. The implementation was made in MATLAB\textsuperscript{\textcopyright}. The results were obtained in a computer with an AMD Ryzen\texttrademark{} Threadripper\texttrademark{} 1950X processor, with 64GB of RAM memory, using the version R2021b of Matlab. 

To solve the linear equilibrium systems, we use a routine of the conjugate gradient method that employs the geometric or the algebraic multigrid method as preconditioner. Multigrid was applied with 4 mesh levels, W cycle, the Jacobi method with relaxation parameter $\omega = 0.5$ as smoother and $\eta_1 = \eta_2 = 1$ smoothing iterations. These parameters were adjusted experimentally, based on several numerical tests. 

The linear programming (LP) subproblems were solved using the dual simplex method, provided by the internal routine {\tt linprog} of Matlab. Following the notation presented in Section \ref{sec: SLP}, 
we use $\delta_{0} = 0.1$ as the initial trust region radius of the SLP algorithm and $\delta_{min} = 10^{-4}$ as the minimum radius. The parameters of the merit function are $\theta_0 = 1$ and $\theta_{max} = 1$. The numerical values of other parameters, such as those used for accepting or rejecting the step, as well as for reducing or increasing the trust region, have been presented in Algorithm \ref{alg: SLP}.

For all of the problems, the initial design variables vector $\mathrm{x}^{(0)}$ was defined by assigning  the upper limit for the volume fraction to all entries, so that the initial point is feasible for the optimization problem. For instance, if the structure can occupy at most $20\%$ of the domain volume, the initial design variables are all equal to $0.2$.

In addition, we explore the symmetries of the structures whenever possible. Although the shape of the structure may not be known in advance, it will be symmetric with respect to a plane that passes through the domain if the supports and external loads are also symmetric with respect to this plane. Considering that the densities of symmetric elements are equal, we can reduce the number of variables in the problem.

\subsection{Tested structures} \label{sec: tests}

We present below the three-dimensional structures that were used in our tests: the cantilever beam, the MBB beam, the L-shaped beam, and the bridge. We have considered different discretizations for each problem, i.e. meshes with different numbers of elements. In all cases, the three dimensions of each finite element measure $1$ unit of length ($cm$, for example), the Young's modulus of the solid material is $E_0 = 1\;N/{cm}^2$, the Poisson's ratio is $\nu = 0.3$, the intensity of the external loads is $1\;N$ (or $1\;N/{cm}^2$ in the case of distributed loads). The Young's modulus of the void is $E_{min} = 10^{-6}E_0$, with the exception of the tests with the adaptive strategy. in which we use $E_{min} = 10^{-9}E_0$. The penalty parameter of the SIMP model is always $p=3$.

\begin{figure}[h]
	\begin{subfigure}{0.5\linewidth}
	  \centering
	  \includegraphics[width=0.7\linewidth]{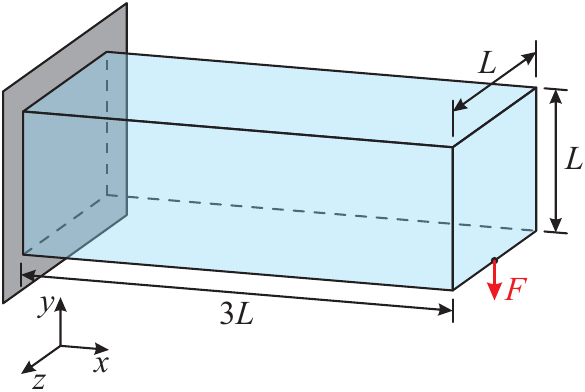}
	  \caption{Cantilever beam.}
	  \label{fig: domainCb}
	\end{subfigure}
	\begin{subfigure}{0.5\linewidth}
	  \centering
	  \includegraphics[width=0.9\linewidth]{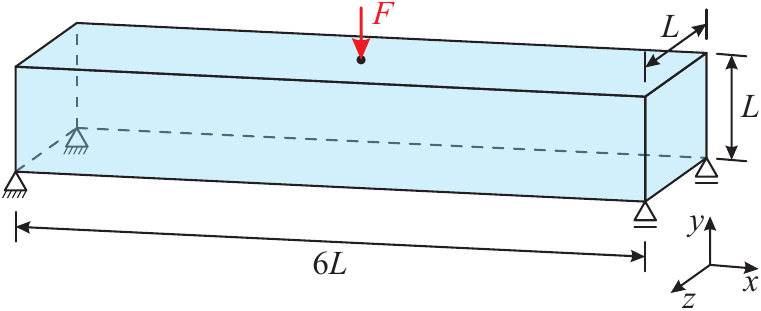}
	  \caption{MBB beam.}
	  \label{fig: domainMbb}
	\end{subfigure}
	\begin{subfigure}{0.5\linewidth}
	  \centering
	  \includegraphics[width=0.656\linewidth]{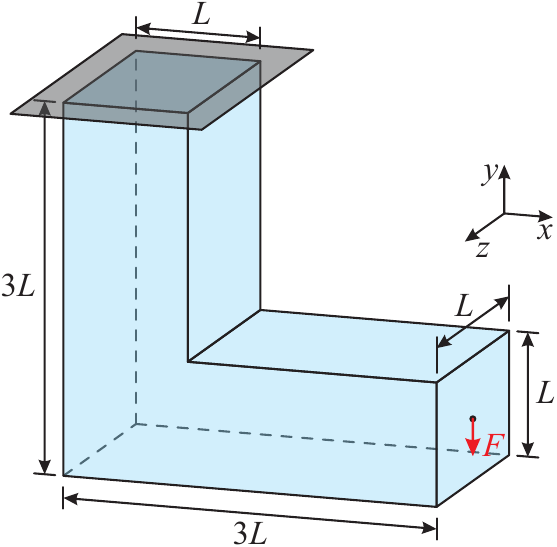}
	  \caption{L-shaped beam.}
	  \label{fig: domainL}
	\end{subfigure}
	\begin{subfigure}{0.5\linewidth}
	  \centering
	  \includegraphics[width=0.9\linewidth]{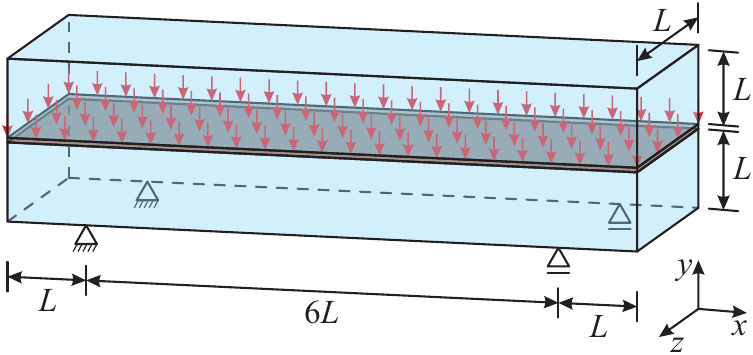}
	  \caption{Bridge.}
	  \label{fig: domainBd}
	\end{subfigure}
	\caption{Initial design domain of the tested structures.}
	\label{fig: domains}
\end{figure}

\begin{description}
	\item[\textbf{Cantilever beam}] \ 

This structure is frequently found in the literature. The left end of the beam is supported (preventing the movement in all directions), as shown in Figure~\ref{fig: domainCb}. A concentrated load is applied downwards at the midpoint of the bottom right edge of the domain. We consider that the final structure can occupy at most $20\%$ of the domain volume ($v_{frac}=0.2$) and that it is symmetric with respect to the $xy$ plane. 

  \item[\textbf{MBB beam}] \ 

Originally produced by the Messerschmitt-Bolkow-Blohm company, this beam is another classic problem. As illustrated in Figure \ref{fig: domainMbb}, it has supports in the four bottom corners and a concentrated load applied downwards at the center of the top face of the domain. The optimal structure can occupy at most $20\%$ of the domain volume and is symmetric with respect to both the $xy$ and the $yz$ planes. The supports are applied to all nodes of the bottom face of each corner element. 

  \item[\textbf{L-shaped beam}] \

The domain of this problem is formed by the union of two right rectangular prisms, that come together in the shape of the letter ``L'', as illustrated in Figure \ref{fig: domainL}. We can also consider that the domain is just one right rectangular prism and that the top right part of the domain is void, i.e. the densities in this region must be zero. The top face is supported and a concentrated load is applied downwards at the center of the right side face. The structure can occupy at most $18\%$ of the domain volume and is symmetric with respect to the $xy$ plane. 

  \item[\textbf{Bridge}] \

Figure \ref{fig: domainBd} illustrates the design domain of this problem. There are four supports in the bottom face, at a distance $L$ from the vertices. There is also a layer that must contain material in the central part of the domain, parallel to the $xz$ plane, where should be the bridge deck. A distributed load is applied downwards over the entire top surface of this layer. The structure can occupy at most $15\%$ of the domain volume and is symmetric with respect to both the $xy$ and the $yz$ planes.

\end{description}

We represent the cantilever beam by the letters {\tt cb}, the MBB beam by {\tt mbb}, the L-shaped beam by {\tt ls} and the bridge by {\tt bd}, followed by the number of elements in each direction $x$, $y$ and $z$ of the mesh, in this order. For instance, {\tt cb60x20x20} represents the cantilever beam with a mesh containing 60 layers of elements in the $x$-direction, 20 layers in the $y$-direction and 20 layers in the $z$-direction.

\subsection{Multigrid tests} 

As previously stated, we solve the linear system $\mathrm{K}(\rho) \mathrm{u} = \mathrm{f}$ using the preconditioned conjugate gradient method (PCG). The available preconditioners include the algebraic multigrid (AMG) and the geometric multigrid (GMG). The incomplete Cholesky factorization (ichol) preconditioner is also available, although it is used mainly for comparison purposes. 

Since multigrid may have a high setup cost, a setup phase is carried out once, before the first SLP iteration, so we just need to recalculate the iteration matrices of the smoother method, the system matrices corresponding to each mesh level and the Cholesky factorization of the matrix on the coarsest mesh at each iteration. 

\begin{figure}[h]
	\centering
	\includegraphics[width=0.64\linewidth]{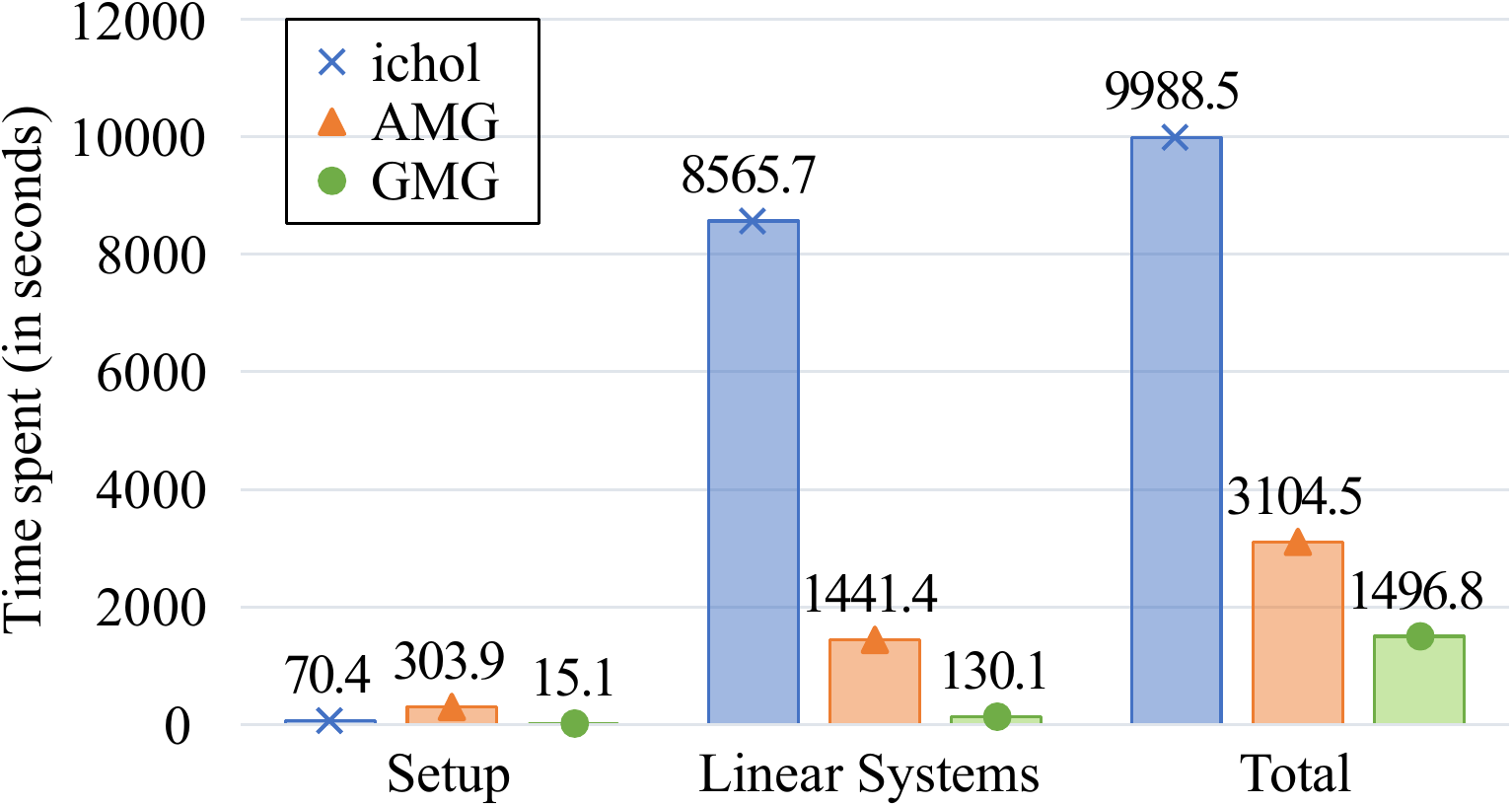}
	\caption{Comparison of the time spent to solve the {\tt mbb288x48x48} problem, using different preconditioners for solving the linear systems. The total time includes all of the remaining steps of the algorithm.}
	\label{fig: MGtimes}
\end{figure}

Figure \ref{fig: MGtimes} presents the time spent by each preconditioner for the {\tt mbb288x48x48} problem. We notice that the setup time is considerable for the AMG preconditioner. On the other hand, the system solution time is substantially reduced when one of the multigrid alternatives is adopted. Using the ichol preconditioner, the time spent solving the linear systems corresponds to $85.8\%$ of the total time, percentage that is reduced to $46.4\%$ and $8.7\%$ for the AMG and the GMG methods, respectively.

When the multigrid method is used, the solution of the linear systems may no longer be the most expensive step of the topology optimization process, especially for large problems. In fact, for the {\tt mbb288x48x48} problem, around 1131 seconds are spent solving linear programming problems, which corresponds to 75\% of the total time for the GMG method but only to 11\% of the time when the ichol preconditioner is used. This result suggests that it may be advantageous to adopt a coarser design variable mesh when using multigrid. We were able to reduce the time spent on this step using the multiresolution approach, as we will see in Subsection \ref{sec: mr_results}.

In general, the GMG preconditioner is the most efficient method, with the only disadvantage that it can only be applied to problems with a regular mesh and some specific discretizations, while the AMG preconditioner is more robust and can be used in a wider variety of problems.

\subsubsection{Number of PCG iterations}

\begin{figure}[h]
	\centering
	\includegraphics[width=0.45\linewidth]{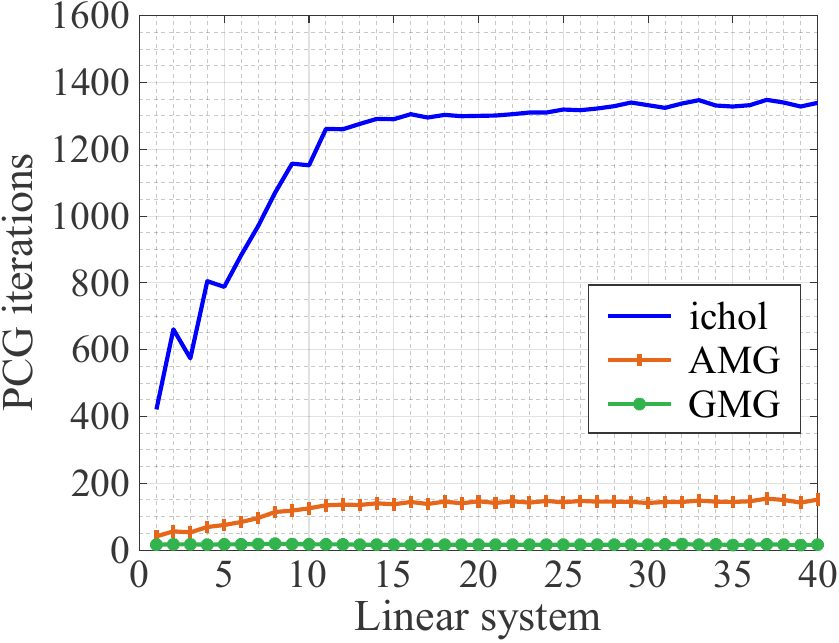}
	\caption{Number of PCG iterations for the {\tt mbb288x48x48} problem.}
	\label{fig: PCGmbb}
\end{figure}

The time spent to solve the linear systems is related to the number of iterations required for the convergence of the PCG method, that is attained when the relative norm of the residual vector is less than $10^{-8}$. Figure \ref{fig: PCGmbb} shows the number of PCG iterations spent solving the linear systems for the {\tt mbb288x48x48} problem.

The PCG method requires a lot more iterations when the ichol preconditioner is used, and less iterations with the GMG preconditioner. We also notice that, for the ichol preconditioner, the number of iterations grows rapidly in the first iterations, probably due to the emergence of void regions in the structure, which may raise the condition number of the stiffness matrix. Since we are using the solution of one linear system as the initial approximation for the next system, the number of PCG iterations becomes more stable when the changes in the shape of the structure become smaller.

\subsubsection{Time growth as a function of mesh refinement} 

\begin{figure}[h]
	\centering
	\includegraphics[width=0.6\linewidth]{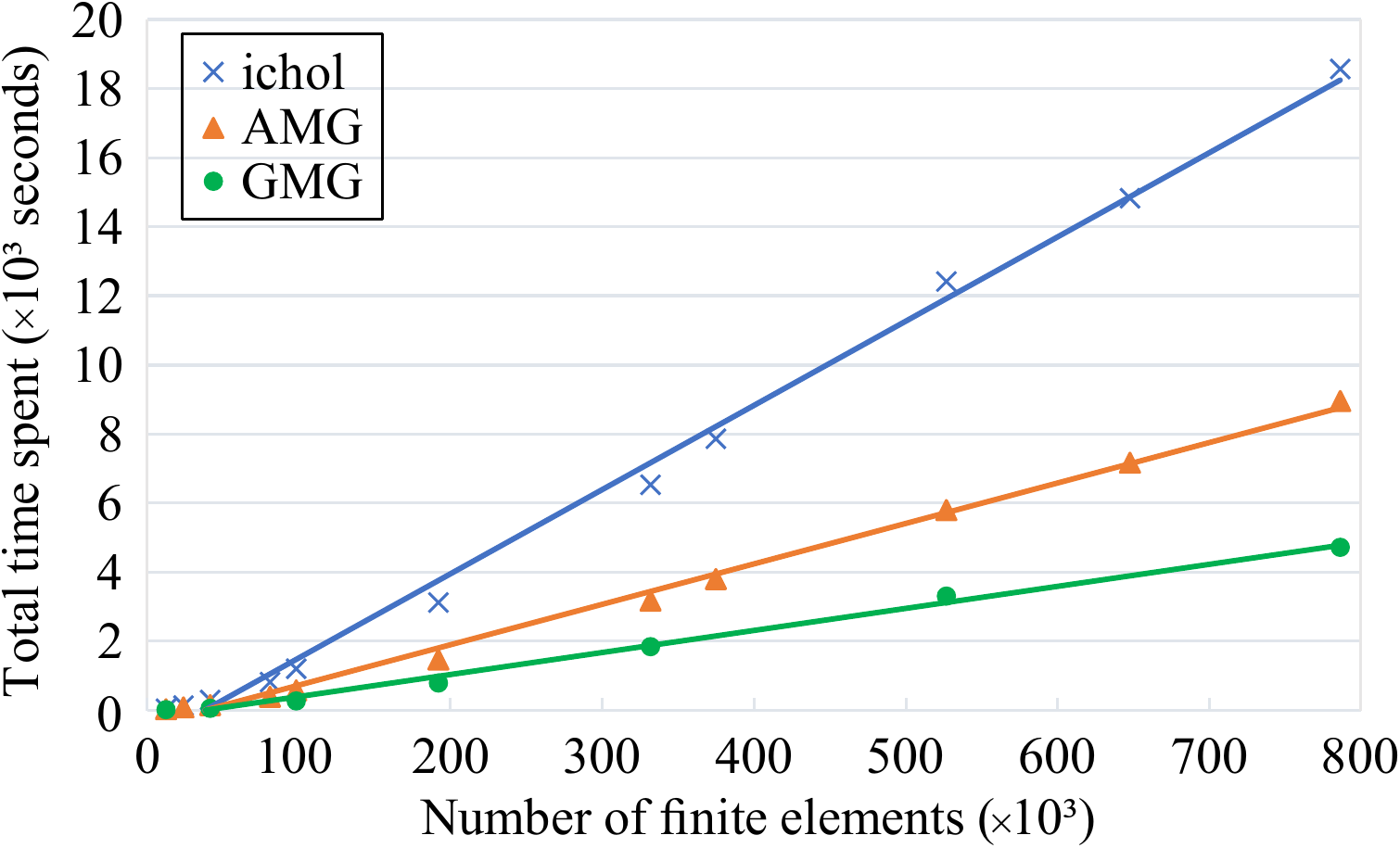}
	\caption{Total time spent by the algorithm as a function of the mesh refinement for the cantilever beam.}
	\label{fig: cbtotaltime}
\end{figure}

The total time spent by our algorithm grows almost linearly with the size of the problem. Figure \ref{fig: cbtotaltime} shows the total time as a function of the number of elements in the mesh for the cantilever beam problem. In this test, we do not consider the symmetries of the structure, so the problems have almost twice the number of variables.

The average rate of growth of the total time with respect to the number of mesh elements is equal to 0.0059 seconds per element for the GMG method, 0.0109 for the AMG method and 0.0226 for the ichol preconditioner. This means that the rate of growth is 84\% higher for AMG method and 280\% higher for the ichol preconditioner compared to the GMG method. 
Due to these results, from now on, we will only use the geometric multigrid preconditioner when solving the linear equilibrium systems.

\subsection{Thresholding strategy tests} 

According to Sigmund \cite{Sigmund2022}, for a more effective comparison among different methods, we must use the same algorithm parameters, consider the same mesh refinement and use boundary conditions that are mesh independent. It is also advisable to recalculate the objective function after projecting the densities, in order to eliminate the influence of intermediate densities in the compliance, and to use the thresholded 0-1 solutions to compare the shapes of the optimal structures. Another good practice is to present the compliance of the fully solid structure, i.e. the structure that covers $100\%$ of the domain, which has the highest possible stiffness. This value gives us information of how good is the solution obtained considering a lower volume fraction.  

In \cite{Sigmund2022}, the author uses a simple strategy to project the variables, that consists of assigning the value $1$ for the highest densities until the maximum desired volume is attained, and considering the remaining density elements as void. In this work, we propose a new thresholding strategy, described in Section \ref{sec: thresholding}, that also takes into account the gradient of the Lagrangian of the problem. 

To evaluate the efficiency of our new thresholding strategy, we compare the compliance of the optimal 0-1 structure with the solution with intermediate densities, the solution given by the simple projection strategy, and the fully solid structure. 
Table \ref{tab: thresholding_compliance} presents the optimal objective function (compliance) values calculated in each case, for the cantilever beam and the MBB beam problems. 

\begin{table}[h]
	\centering
	\caption{Compliance values obtained with the thresholding strategy.}
	\def\arraystretch{1.2}
	\setlength{\tabcolsep}{8pt}
	\small
	\begin{tabular}{c|rrrr}
		\hline 
		Problem & \multicolumn{1}{c}{$F$} & \multicolumn{1}{c}{$F_{rnd}$} & \multicolumn{1}{c}{$F_{prj}$} & \multicolumn{1}{c}{$F_{sol}$} \\ \hline 
		{\tt cb192x64x64} & 21.566 & 18.549 & 18.546 & 11.108 \\ 
		\hline 
		{\tt mbb240x40x40} & 39.597 & 26.800 & 26.762 & 13.285 \\ 
		\hline 
	\end{tabular} 
	\label{tab: thresholding_compliance}
\end{table}

For these examples, as expected, the objective function value for the structure with intermediate densities ($F$) was always the highest, and the compliance of the fully solid structure ($F_{sol}$) was always the smallest. Futhermore, the compliance obtained with the new thresholding strategy ($F_{prj}$) is slightly smaller than the one given by the simple thresholding strategy ($F_{rnd}$). 

We control the quality of the solution by calling the SLP algorithm again using the thresholded solution as the initial density vector. This process is repeated until the maximum of $10$ thresholding attempts is reached or until the difference between two consecutive solutions is relatively small and the violation in the volume constraint is within a prescribed tolerance. In the tests presented here, one attempt was sufficient for the algorithm to reach these stopping criteria. Table \ref{tab: thresholding} shows the number of outer iterations required for the convergence of the SLP algorithm ($N_{it}$), the number of rejected steps ($N_{rs}$), the time spent (in seconds) with the thresholding strategy ($T_{prj}$) and the total time spent (in seconds) to solve the problem ($T_{total}$).
\begin{table}[h]
	\centering
	\caption{Results obtained with the thresholding strategy.}
	\def\arraystretch{1.2}
	\setlength{\tabcolsep}{8pt}
	\small
	\begin{tabular}{c|c|rrr}
		\hline 
		Problem & SLP call & \multicolumn{1}{c}{$N_{it}$ ($N_{rs}$)} & \multicolumn{1}{c}{$T_{prj}$} & \multicolumn{1}{c}{$T_{total}$} \\ \hline 
		\multirow{2}{*}{\tt cb192x64x64} & Initial solution & 25 (+1) & 0.01 & 3992.9 \\
		& After thresholding & 3 (+0) & 0.8 & 169.3 \\
		\hline 
		\multirow{2}{*}{\tt mbb240x40x40} & Initial solution & 36 (+4) & 0.01 & 681.0 \\
		& After thresholding & 7 (+0) & 0.2 & 51.5 \\
		\hline 
	\end{tabular} 
	\label{tab: thresholding}
\end{table} 

The SLP calls after the thresholding consume just a few iterations, since the thresholded solution is already close to a local minimizer. Moreover, the time spent to apply the thresholding strategy is insignificant in relation to the total time. $T_{prj}$ is slightly higher after thresholding because it includes the time spent with the Heaviside projection. 
Table \ref{tab: densities} presents the numbers of void elements ($\rho = 0$), intermediate elements ($\rho \in (0,1)$) and solid elements ($\rho = 1$) in the mesh, for the original and the thresholded optimal solutions, considering the symmetry of the structures.   

\begin{table}[h]
	\centering
	\caption{Material densities distribution in the mesh.}
	\def\arraystretch{1.2}
	\setlength{\tabcolsep}{8pt}
	\small
	\begin{tabular}{c|c|ccc}
		\hline 
		Problem & Solution & $\rho = 0$ & $\rho \in (0,1)$ & $\rho = 1$ \\ \hline 
		\multirow{2}{*}{\tt cb192x64x64} & Original & 218982 & 153938 & 20296 \\
		& Thresholded & 314545 & 54 & 78617 \\
		\hline 
		\multirow{2}{*}{\tt mbb240x40x40} & Original & 43972 & 51765 & 263 \\
		& Thresholded & 76800 & 0 & 19200 \\
		\hline 
	\end{tabular} 
	\label{tab: densities}
\end{table}

We notice that the original solutions have a significant number of intermediate densities, justifying the difference between the values $F$ and $F_{prj}$ presented in Table \ref{tab: thresholding_compliance}. Our thresholding strategy was able to completely eliminate the intermediate densities in the {\tt mbb240x40x40} problem, keeping the structure volume at $20\%$ of the domain volume. In the {\tt cb192x64x64} problem, $54$ intermediate densities remained after thresholding, but this represents only $0.014\%$ of the total number of elements. Also in this case, the structure volume was kept at $20\%$ of the domain volume.

\subsection{Multiresolution tests} \label{sec: mr_results}

We now present the results obtained with the combination of the SLP algorithm and the multiresolution technique described in Section \ref{sec: mr}. The results refer to the thresholded 0-1 solutions, so the methods are compared in a more effective way.

\subsubsection{Comparisons with the traditional method} 

First, we compare the solutions obtained using the multiresolution approach and the traditional topology optimization method. In multiresolution, the density and the design variable meshes are constructed from the displacement mesh, dividing each displacement element into ${(n_{mr})}^3$ density elements and ${(d_{mr})}^3$ design variables. Some combinations of the parameters $n_{mr}$ and $d_{mr}$ were tested.

Since all the edges of the displacement elements have length equal to $1$, the edges of the density elements measure $1/n_{mr}$, so the filter radius (used in the projection of design variables to density elements) can be smaller than the one adopted in the traditional method. In our tests, we choose a radius $r_{min}$ for multiresolution, and define $n_{mr} r_{min}$ as the radius for the traditional method.
Furthermore, we try to keep the boundary conditions compatible applying the supports to the same region for both methods. Nevertheless, we remark that the scenarios in which we solve the problem are still different and it is well known that topology optimization problems have many local minima, so the results obtained by the two methods may not be the same.

\begin{figure}[h]
	\begin{subfigure}{0.5\linewidth}
		\centering
		\includegraphics[width=0.9\linewidth]{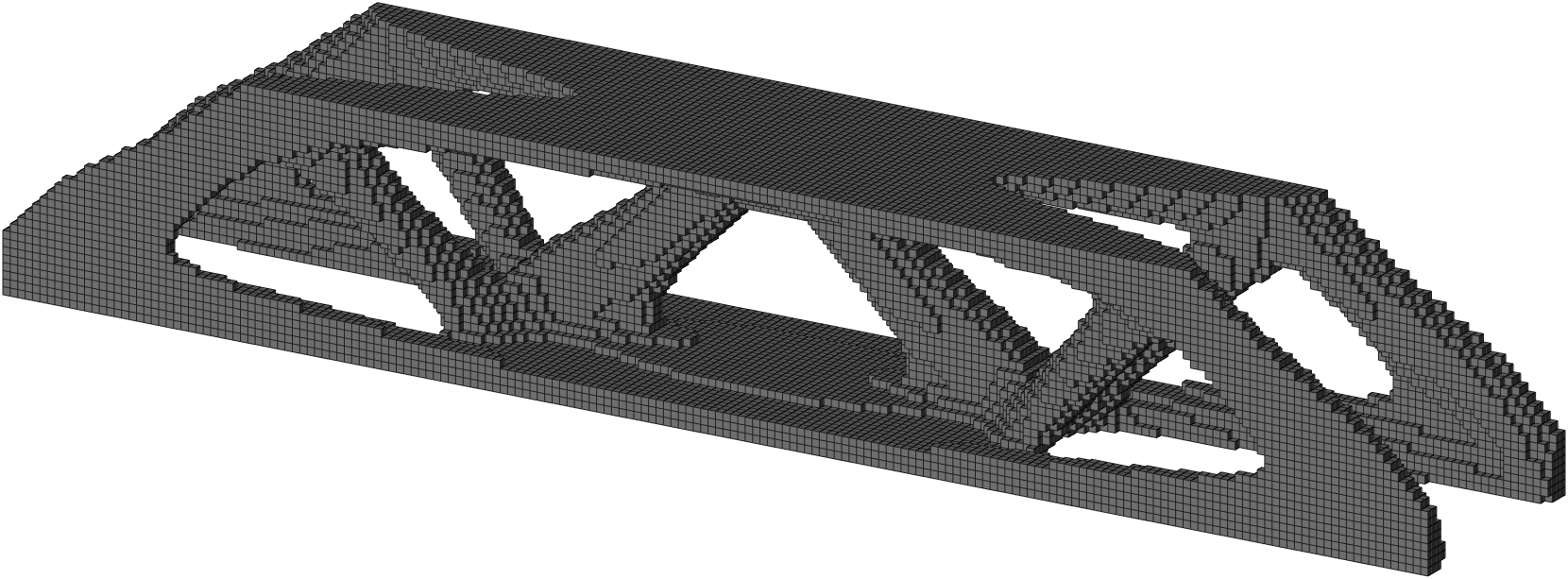}
		\caption{Traditional method.}
	\end{subfigure}
	\begin{subfigure}{0.5\linewidth}
		\centering
		\includegraphics[width=0.9\linewidth]{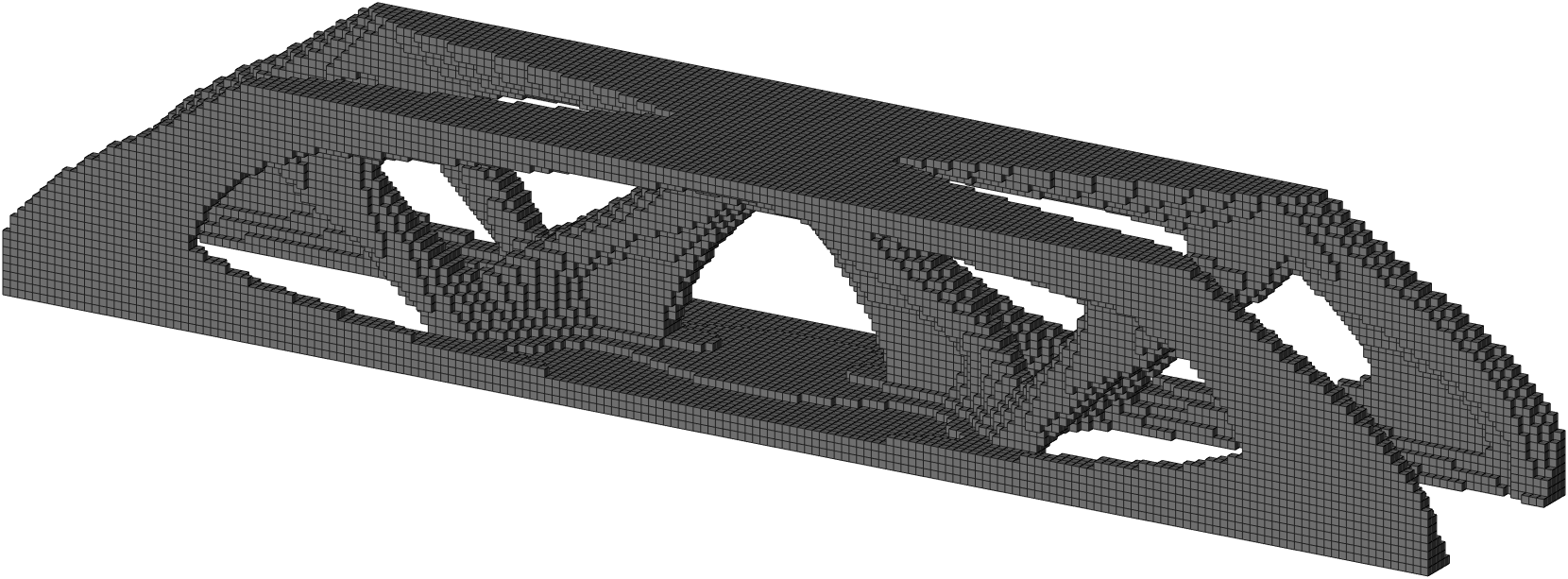}
		\caption{Multiresolution with $n_{mr} = 2$ and $d_{mr} = 2$.}
	\end{subfigure}
	\caption{Optimal topologies for the {\tt mbb192x32x32} problem ($r_{min} = 2.0$).}
	\label{fig: mbb192x32x32_mr}
\end{figure}

\begin{figure}[h]
  \centering
	\begin{subfigure}{0.45\linewidth}
		\centering
		\includegraphics[width=\linewidth]{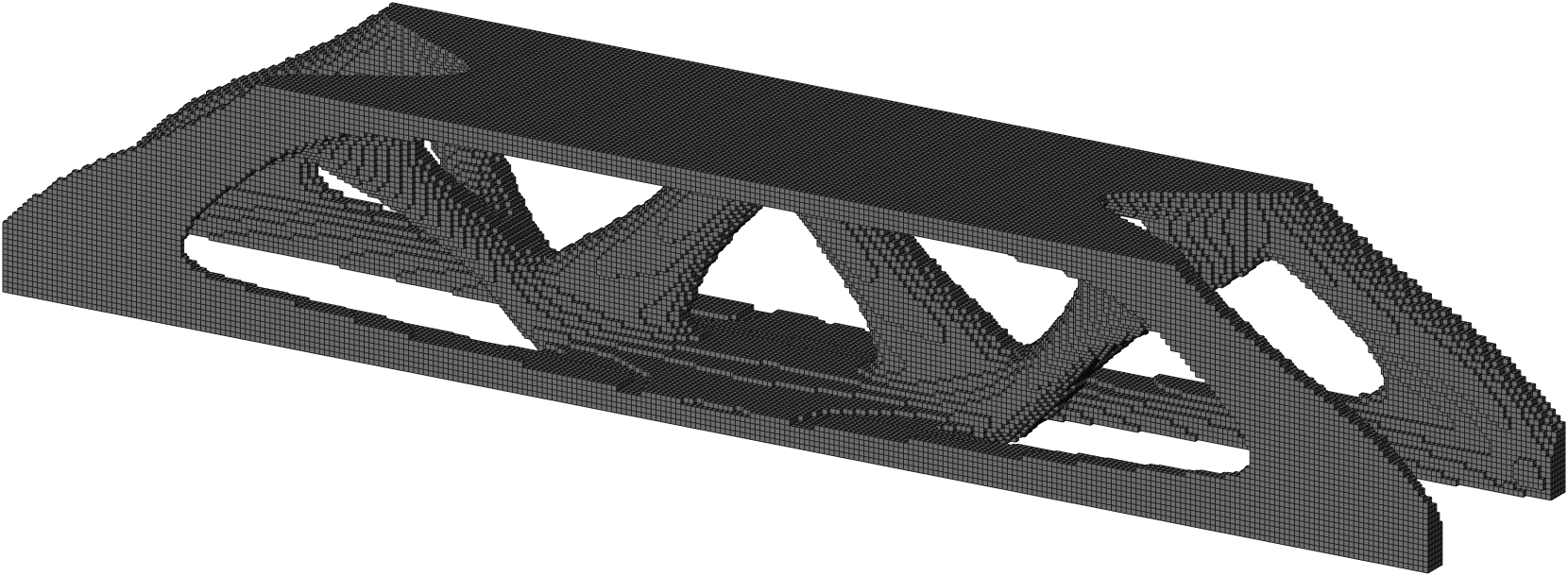}
		\caption{Traditional method.}
	\end{subfigure}\hspace{0.5cm}
	\begin{subfigure}{0.45\linewidth}
		\centering
		\includegraphics[width=\linewidth]{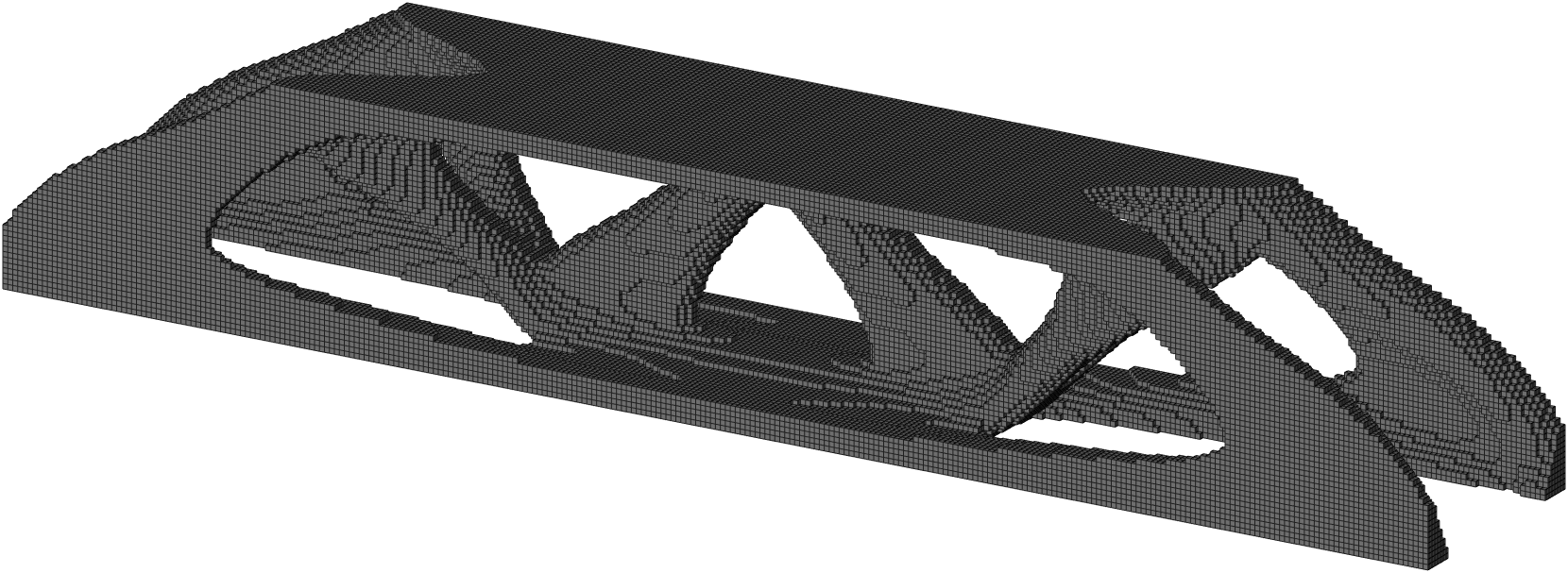}
		\caption{Multiresolution with $n_{mr} = 3$ and $d_{mr} = 3$.}
	\end{subfigure}
	\begin{subfigure}{0.45\linewidth}
		\centering
		\includegraphics[width=\linewidth]{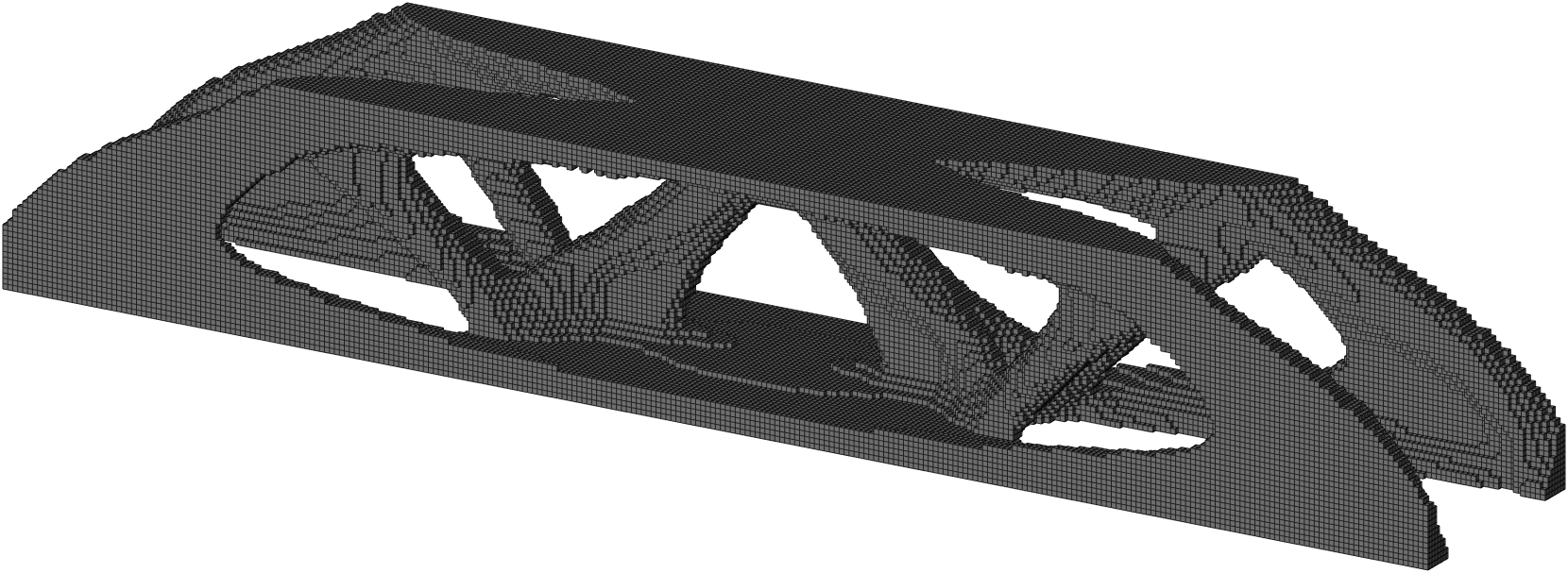}
		\caption{Multiresolution with $n_{mr} = 3$ and $d_{mr} = 2$.}
	\end{subfigure}
	\caption{Optimal topologies for the {\tt mbb288x48x48} problem ($r_{min} = 2.0$).}
	\label{fig: mbb288x48x48_mr}
\end{figure}

In our first test, we divide the MBB beam into 96x16x16 displacement elements and try to find structures with 192x32x32 and 288x48x48 density elements, using $n_{mr}$ equal to $2$ and $3$, respectively. For the larger structure, we consider two alternatives for $d_{mr}$. In all cases, we use $r_{min} = 2.0$. A comparison of the results obtained for both the traditional method and the multiresolution (MR) is presented in Table \ref{tab: mbb_mr}. Figures \ref{fig: mbb192x32x32_mr} and \ref{fig: mbb288x48x48_mr} show the structures obtained after applying the thresholding strategy. 

\begin{table*}[h]
	\centering
	\caption{Comparison between the traditional method and multiresolution for two MBB beam problems.}
	\def\arraystretch{1.2}
	\setlength{\tabcolsep}{8pt}
	\small
	\begin{tabular}{cc|crrr}
		\hline 
		Problem & Method & \multicolumn{1}{c}{$N_{it}$ ($N_{rs}$)} & \multicolumn{1}{c}{$F$} & \multicolumn{1}{c}{$F_{prj}$} & \multicolumn{1}{c}{$T_{total}$} \\ \hline 
		\multirow{2}{*}{\tt mbb192x32x32} & Traditional & 149 (+49) & 42.208 & 29.584 & 1007.9 \\ 
		& MR ($n_{mr} = 2$, $d_{mr} = 2$) & \phantom{1}99 (+28) & 41.884 & 29.097 & 404.6 \\ 
		\hline 
		\multirow{3}{*}{\tt mbb288x48x48} & Traditional & \phantom{1}31 \phantom{1}(+1) & 30.603 & 21.857 & 1310.4 \\ 
		& MR ($n_{mr} = 3$, $d_{mr} = 3$) & \phantom{1}54 (+16) & 30.731 & 21.557 & 2038.4 \\ 
		& MR ($n_{mr} = 3$, $d_{mr} = 2$) & 119 (+35) & 30.053 & 21.527 & 508.0 \\ 
		\hline 
	\end{tabular} 
	\label{tab: mbb_mr}
\end{table*}

For $n_{mr} = 2$, the number of SLP iterations ($N_{it}$) and the total solution time ($T_{total}$) clearly decreased with multiresolution. For $n_{mr} = 3$, the number of iterations increased when multiresolution was adopted. Moreover, for $d_{mr} = 3$, the time spent to solve the problem was greater than that consumed by the traditional method. As will become clear later in this section, this happened because the solution of the LP problems is the most time consuming step of the multiresolution algorithm and, when $n_{mr} = d_{mr}$, the size of the LP problems are the same for the multiresolution and the traditional method, so an increase in the number of iterations jeopardizes the time saved in solving the linear systems on the displacement mesh. This problem can be overcame reducing the value of $d_{mr}$ to 2. In this case, even with a considerably increase in the number of iterations, we notice a significant reduction in the total time. 

\begin{figure}[h]
	\centering
	\begin{subfigure}{0.45\linewidth}
		\centering
		\includegraphics[width=\linewidth]{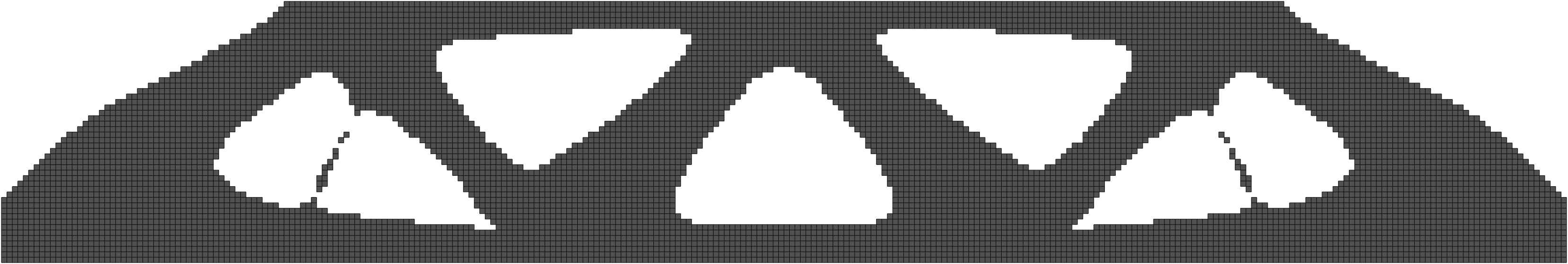}
		\caption{Multiresolution with $n_{mr} = 3$ and $d_{mr} = 3$ \\ (simple thresholding strategy). \centering}
	\end{subfigure} \hspace{0.5cm}
	\begin{subfigure}{0.45\linewidth}
		\centering
		\includegraphics[width=\linewidth]{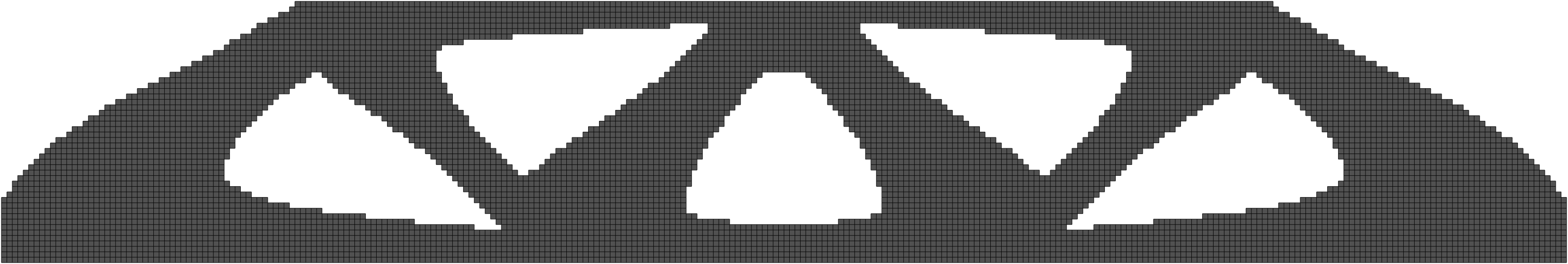}
		\caption{Multiresolution with $n_{mr} = 3$ and $d_{mr} = 3$ \\ (new thresholding strategy). \centering}
	\end{subfigure}
	\caption{Side views of the MBB beams.}
	\label{fig: mbb_side}
\end{figure}

In general, the structures obtained with multiresolution have similar shapes to those obtained with the traditional method. However, applying the simple thresholding strategy to the optimal solution with intermediate densities may result in a structure with some incomplete members, as shown in Figure \ref{fig: mbb_side} for the {\tt mbb288x48x48} problem with $n_{mr} = d_{mr} = 3$.  The side views of the structures make it clear that these artefacts can be completely eliminated using our new thresholding strategy. 

It should be noted that, when linear finite elements are used, a value greater than $2$ for the parameter $d_{mr}$ violates the upper limit imposed in Table \ref{tab: bounds}. According to Gupta {\it et al.} \cite{Gupta2017}, this may lead to non-unique optimal solutions and favor the occurrence of artefacts. However, we see that good results can still be obtained in some cases, if a careful choice of the density filter and thresholding strategy is made. 

As mentioned before, the problems solved with multiresolution and with the traditional method are not exactly the same, so it is difficult to properly compare the compliance values achieved by both methods. To circumvent this problem, we take the optimal density vector found with multiresolution and use it to compute the compliance in the same mesh as the traditional method. This correction was used to compute the function values presented in Table \ref{tab: mbb_mr} for both the solution with intermediate densities ($F$) and the thresholded solution ($F_{prj}$). As we see, the function values are slightly smaller for multiresolution with $d_{mr} = 2$.

In order to obtain more conclusive results about the accuracy of the multiresolution solutions, we carried out a new test in which the optimal design variables obtained with multiresolution are used as an initial guess and the problem is solved again with the traditional method. Our purpose is to find out if the objective function value given by multiresolution ($F_{mr}$) is close to the one found by the traditional method ($F_{pos}$). 
Table \ref{tab: mr_pos} shows both function values, as well as the number of SLP iterations ($N_{pos}$) needed to solve the problem again, for MBB beam problems with $n_{mr} = d_{mr}$. The figures in the table show that the number of iterations $N_{pos}$ is small, and that the ratio $F_{mr}/F_{pos}$ is approximately equal to $1$ in all cases, indicating that the solutions obtained with the multiresolution scheme are considerably good.

\begin{table}[h]
	\centering
	\caption{Results obtained applying the traditional method using the MR solution as the initial guess.}
	\def\arraystretch{1.2}
	\setlength{\tabcolsep}{8pt}
	\small
	\begin{tabular}{c|crrc}
		\hline 
		Problem & $N_{pos}$ & \multicolumn{1}{c}{$F_{mr}$} & \multicolumn{1}{c}{$F_{pos}$} & \multicolumn{1}{c}{$F_{mr}/F_{pos}$} \\ \hline 
		{\tt mbb192x32x32} & 7 (+2) & 41.884 & 41.852 & 1.001 \\
		\hline 
		{\tt mbb288x48x48} & 3 (+0) & 30.731 & 30.704 & 1.001 \\
		\hline
	\end{tabular} 
	\label{tab: mr_pos}
\end{table}

Let us now investigate the time spent on each step of the multiresolution algorithm. In our analysis, we consider only the first $20$ outer iterations of the SLP algorithm, in order to allow a fair comparison between the scenarios. The results obtained for the {\tt mbb288x48x48} problem are shown in Table \ref{tab: timesmbb_mr}, where the notation MR ($n_{mr}$, $d_{mr}$) represents the multiresolution with parameters $n_{mr}$ and $d_{mr}$. 

\begin{table}[h]
	\centering
	\caption{Time spent (in seconds) on each step of the algorithm for the {\tt mbb288x48x48} problem, considering the first 20 SLP outer iterations.}
	\def\arraystretch{1.2}
	\setlength{\tabcolsep}{8pt}
	\small
	\begin{tabular}{l|ccc}
		\hline 
		\multicolumn{1}{c|}{Step} & \multicolumn{1}{c}{Traditional} & \multicolumn{1}{c}{MR (3, 3)} & \multicolumn{1}{c}{MR (3, 2)} \\ \hline 
		Pre-filtering & \phantom{7}63.7 & \phantom{7}64.5 & \phantom{7}6.8 \\
		Filtering & \phantom{75}2.7  & \phantom{75}3.0 & \phantom{7}0.8 \\
		Assembly of matrix $\mathrm{K}$ & \phantom{7}79.7 & \phantom{75}3.6 & \phantom{7}3.4 \\
		Preconditioner setup & \phantom{75}8.5 & \phantom{75}0.4 & \phantom{7}0.4 \\
		Linear systems & \phantom{7}68.1 & \phantom{75}7.5 & \phantom{7}7.3 \\
		Gradients & \phantom{75}6.5 & \phantom{75}2.7 & \phantom{7}2.4 \\
		LP subproblems & 752.0 & 830.4 & 77.9 \\
		Other & \phantom{75}0.3  & \phantom{75}0.1 & \phantom{7}0.1 \\
		\hline 
		Total & \multicolumn{1}{c}{981.7} & \multicolumn{1}{c}{912.2} & \multicolumn{1}{c}{98.9} \\ \hline 
	\end{tabular} 
	\label{tab: timesmbb_mr}
	\vspace{-0.3cm}
\end{table}

From the table, we notice a reduction in the time spent assembling the stiffness matrix, setting up the preconditioner and solving the linear systems with the multiresolution, due to the fact that the equilibrium systems are solved in a coarsest displacement mesh. When $n_{mr} = d_{mr} = 3$, the number of design variables is equal in both methods, so that the time spent on the pre-filtering and filtering steps were similar. In contrast, the time spent solving the LP subproblems was greater, suggesting that these subproblems become more complicated. Altogether, the time spent to solve the problem was reduced by only 7.1\% in this case. 

A much greater reduction was achieved using $d_{mr} = 2$. In this case, the design variable space is significantly small, so the time spent on the filtering process and especially on solving LP problems has been significantly reduced.

\subsubsection{The artefacts drawback} 

Refining the density mesh (i.e. increasing the value of $n_{mr}$) leads to a reduction in the size of the elements. If the refinement is done keeping the filter radius, the region of application of the filter remains the same for all resolutions, so the structures obtained are mesh independent, as we can see in Figures \ref{fig: mbb192x32x32_mr} and \ref{fig: mbb288x48x48_mr}. However, mesh refinement is generally used with the aim of obtaining more detailed structures, i.e. structures with new holes or members that cannot be formed when using a low resolution. To achieve this objective, we also need to reduce the filter radius.  

\begin{figure}[h]
	\centering 
	\includegraphics[width=0.45\linewidth]{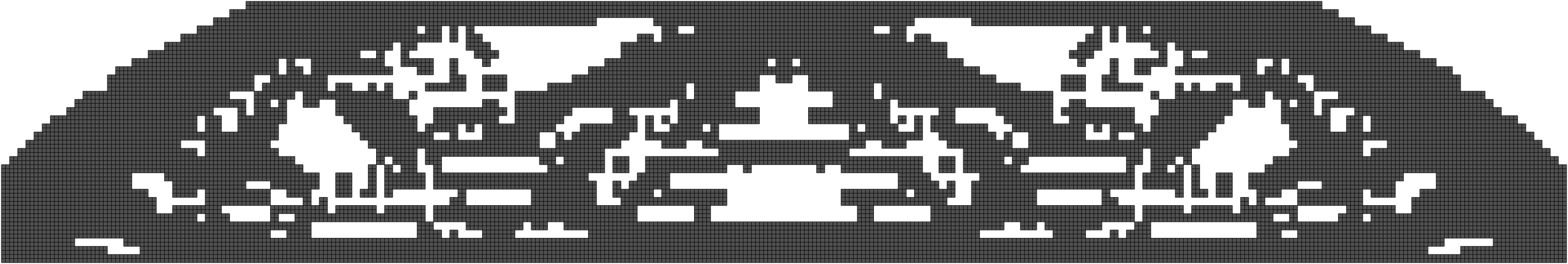}
	\caption{MBB beam obtained using MR with $n_{mr} = 4$, $d_{mr} = 2$, $r_{min} = 0.6$ and linear finite elements.}
	\label{fig: mbb_artefacts}
\end{figure}

Unfortunately, reducing the filter radius may lead to the emergence of artefacts, as shown in Figure \ref{fig: mbb_artefacts} for the {\tt mbb384x64x64} problem with $r_{min} = 0.6$. One way to alleviate this drawback is to apply a density filter. However, to completely eliminate the artefacts we need to use a large density filter radius, which brings us back to the low resolution trouble. Therefore, if we want to obtain structures with more details, we need to rely on other strategies, and increasing the degree of the polynomials used to approximate the displacements is one of the most effective alternatives.

Next, we present the results obtained increasing the degree of the elements used to approximate the displacements. The disadvantage of using finite elements with a higher degree is that the size of the stiffness matrices also increase, so the solution of the linear systems can be again the most expensive step of the algorithm, even with multigrid methods. To deal with this problem, we suggest the adoption of an adaptive strategy to increase the degree of the elements.

\subsection{High-order element tests}

To measure the impact of using higher order finite elements, we once again solve the MBB beam problem using multiresolution with $n_{mr} = 4$, $d_{mr} = 2$, and $r_{min} = 0.6$. We compare the use of elements of the Lagrange and serendipity families with degree $1$, $2$ and $3$, presented in Subsection \ref{sec: elements}. In our tests, we adopt the SSOR algorithm as the smoother of the multigrid method, since it performs better than the Jacobi method for quadratic and cubic elements.

\subsubsection{Comparison between Lagrange and serendipity elements} 

\begin{figure*}[h]
	\begin{subfigure}{0.5\linewidth}
		\centering
		\includegraphics[width=0.9\linewidth]{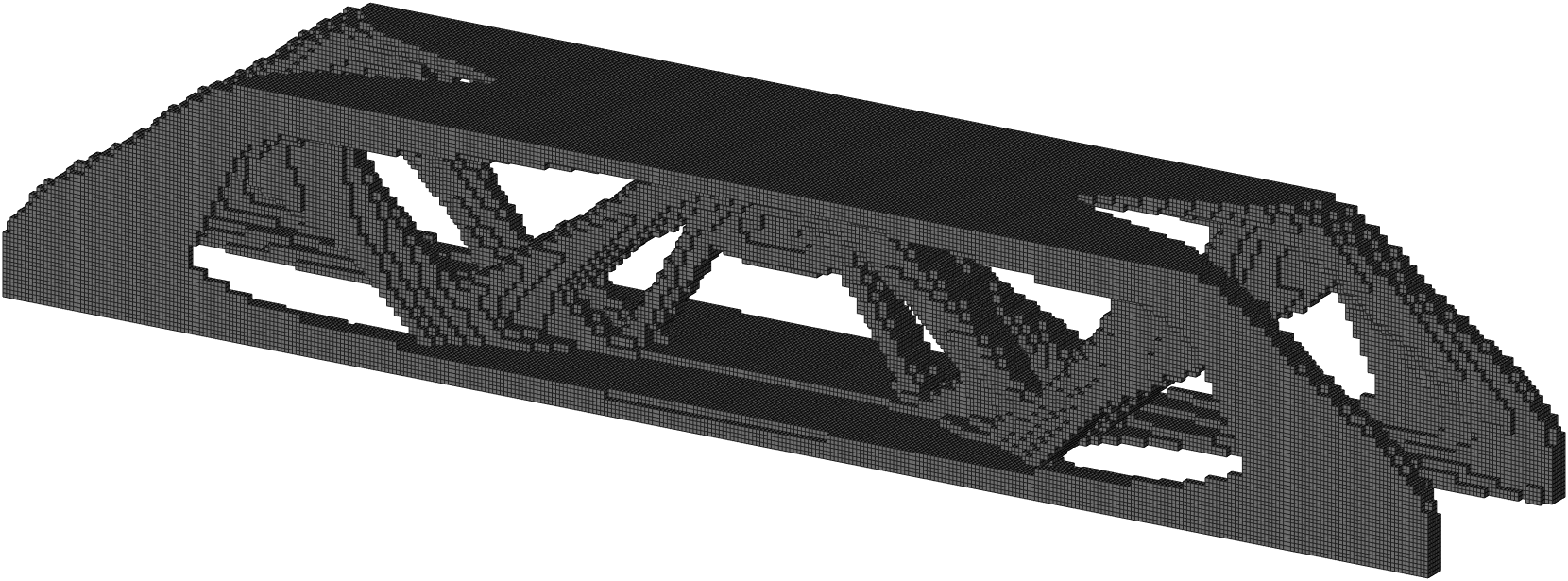}
	\end{subfigure}%
	\begin{subfigure}{0.5\linewidth}
		\centering
		\includegraphics[width=0.9\linewidth]{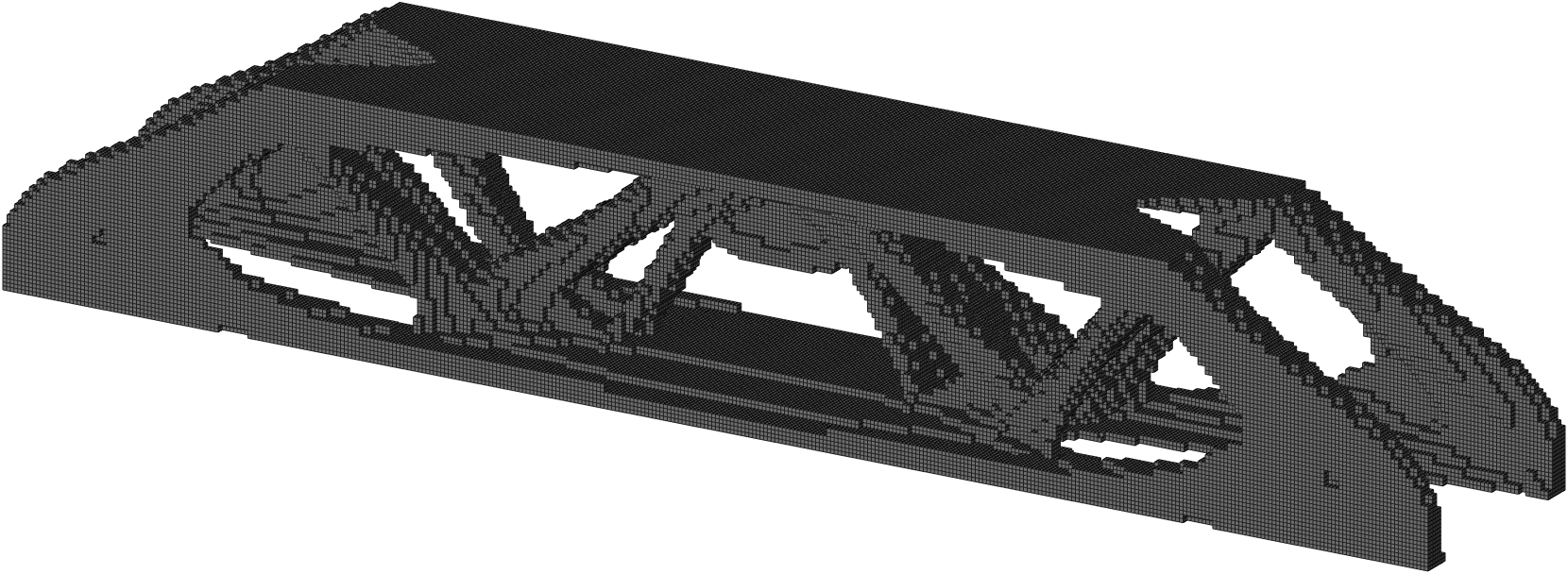}
	\end{subfigure}
	
	\begin{subfigure}{0.5\linewidth}
		\centering
		\includegraphics[width=0.9\linewidth]{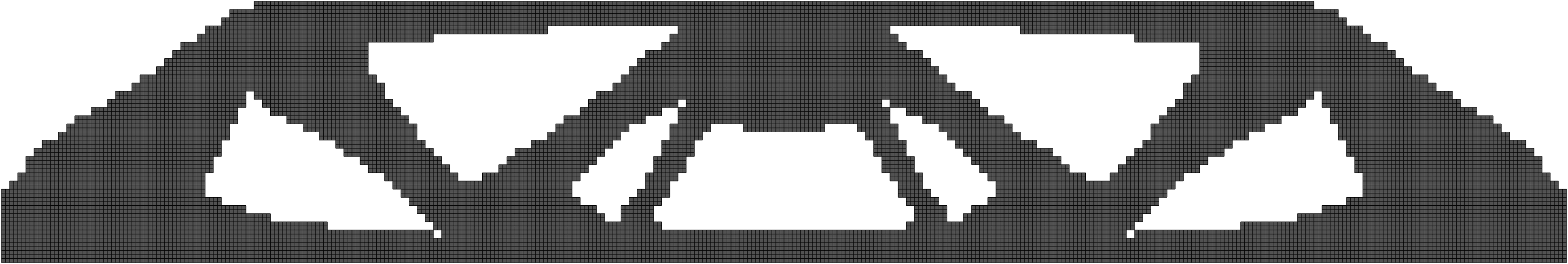}
		\caption{Lagrange elements, degree $2$.}
	\end{subfigure}%
	\begin{subfigure}{0.5\linewidth}
		\centering
		\includegraphics[width=0.9\linewidth]{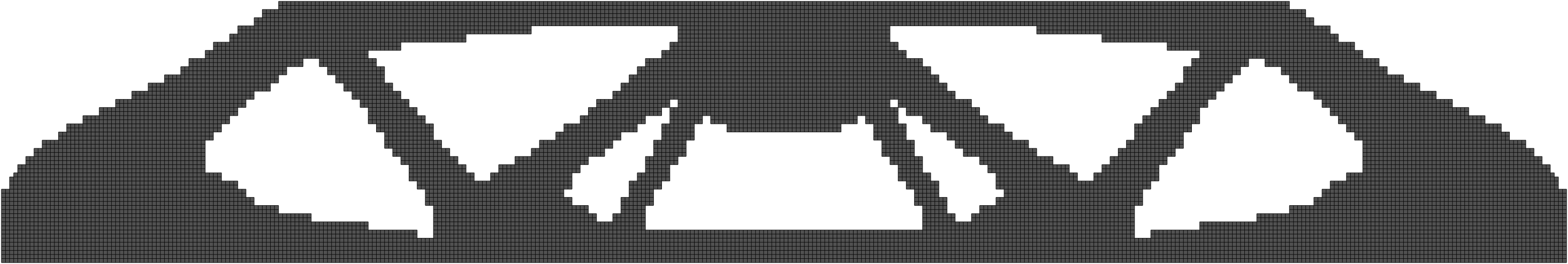}
		\caption{Serendipity elements, degree $2$.}
	\end{subfigure}
	
	\begin{subfigure}{0.5\linewidth}
		\centering
		\includegraphics[width=0.9\linewidth]{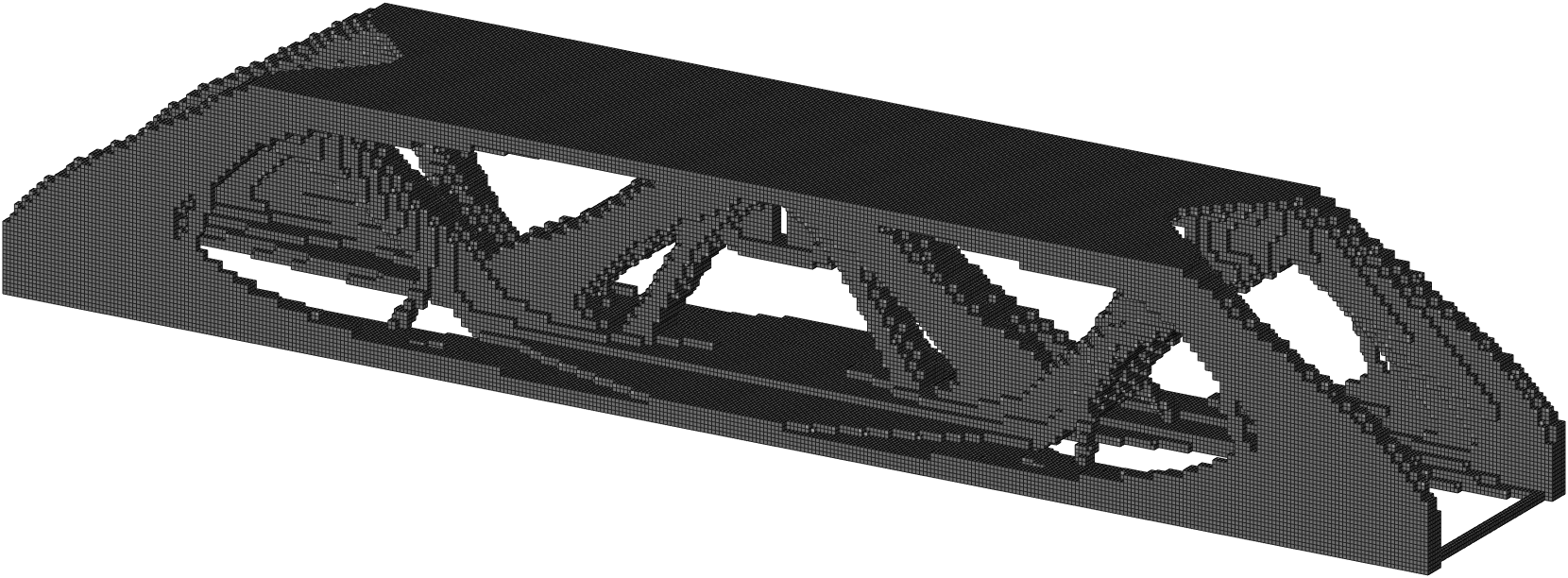}
	\end{subfigure}%
	\begin{subfigure}{0.5\linewidth}
		\centering
		\includegraphics[width=0.9\linewidth]{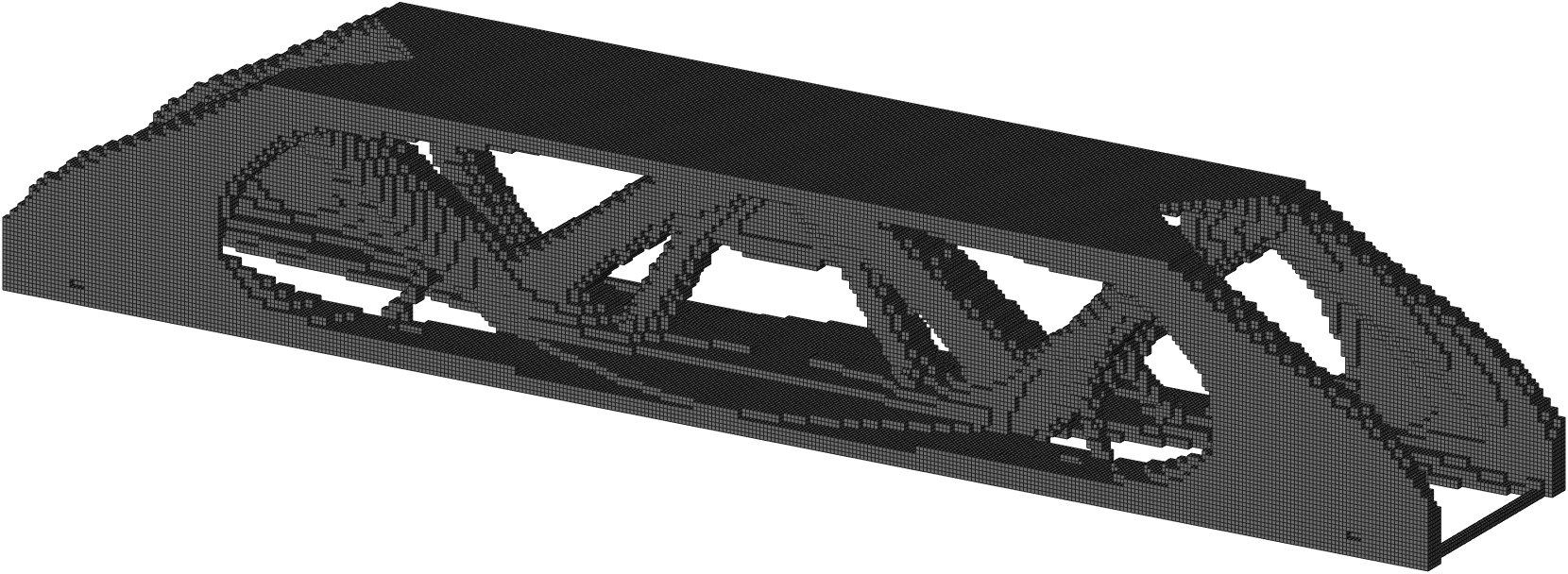}
	\end{subfigure}
	
	\begin{subfigure}{0.5\linewidth}
		\centering
		\includegraphics[width=0.9\linewidth]{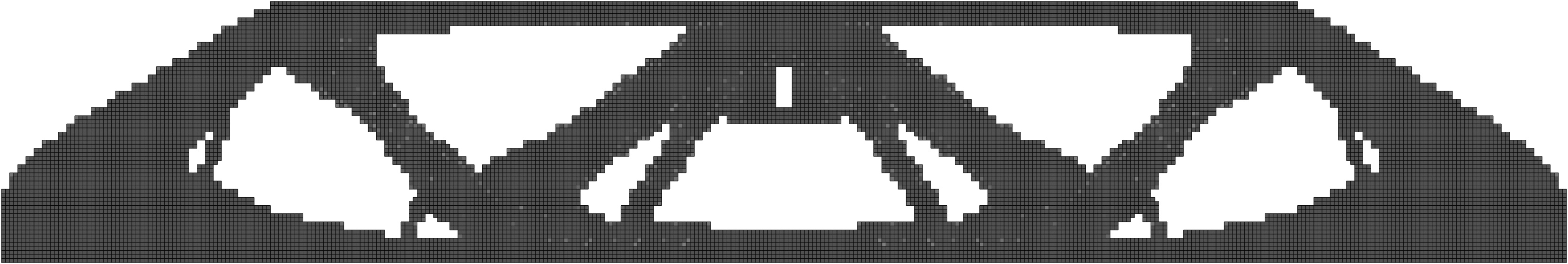}
		\caption{Lagrange elements, degree $3$.}
	\end{subfigure}%
	\begin{subfigure}{0.5\linewidth}
		\centering
		\includegraphics[width=0.9\linewidth]{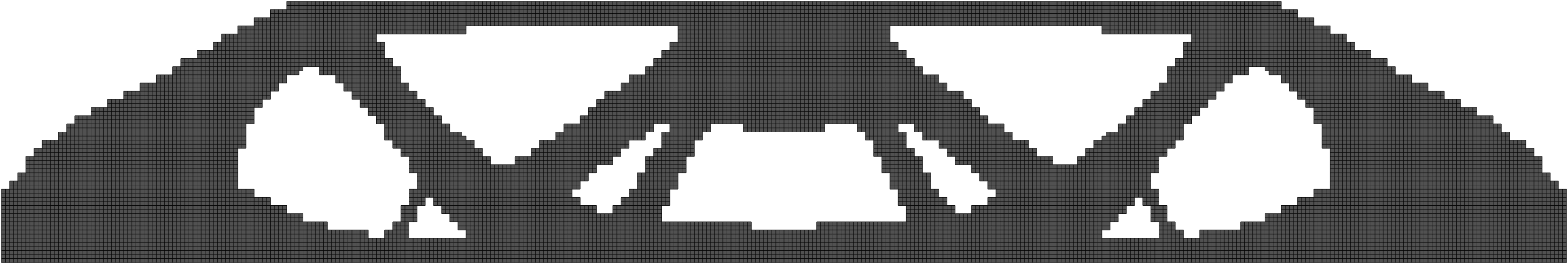}
		\caption{Serendipity elements, degree $3$.}
	\end{subfigure}
	
	\caption{Structures obtained for the {\tt mbb384x64x64} problem, using MR with $n_{mr} = 4$, $d_{mr} = 2$, $r_{min} = 0.6$.}
	\label{fig: mbb_deg}
\end{figure*}

As shown in Figure \ref{fig: mbb_artefacts}, the material distribution obtained for the MBB beam applying multiresolution with linear elements is full of artefacts, making it difficult to interpret the shape of the structure. Figure \ref{fig: mbb_deg} shows the topologies obtained for the same problem using elements with degree greater than $1$. In this figure, we display the structures after applying the thresholding strategy, with oblique and side views to allow a better visualization of the details.

From the figure, we clearly notice that, using a better approximation of the displacements, the members of the structure became better defined and the artefacts disappeared. Moreover, the topologies obtained with elements of degree $3$ are more detailed than those obtained with degree $2$. For three-dimensional problems, it seems to be enough to use elements with degree $2$ or $3$ to prevent the occurrence of artefacts. 

Table \ref{tab: increase_deg} shows the numerical results obtained using different types of elements. We use the notation L$d$ (S$d$) to represent an element of the Lagrange (serendipity) family with degree $d$. $\tilde{F}_{prj}$ is the optimal objective function value for the thresholded solution and $F_{prj}$ is the ``corrected'' compliance calculated after computing the nodal displacements in the density mesh with linear elements.   

\begin{table}[h]
	\centering
	\caption{Results obtained for the {\tt mbb384x64x64} problem using MR with $n_{mr} = 4$, $d_{mr} = 2$, $r_{min} = 0.6$ and different types of elements.}
	\def\arraystretch{1.2}
	\setlength{\tabcolsep}{8pt}
	\small
	\begin{tabular}{c|crrr}
		\hline 
		Element & $N_{it}$ ($N_{rs}$) & \multicolumn{1}{c}{$\tilde{F}_{prj}$} & \multicolumn{1}{c}{$F_{prj}$} & \multicolumn{1}{c}{$T_{total}$} \\ \hline 
		L1 & 177 (+0) & 49.590 & 88.658 & 371.8 \\
		L2 & 41 (+1) & 62.241 & 17.811 & 713.7 \\
		S2 & 43 (+1) & 57.039 & 18.039 & 1136.1 \\
		L3 & 86 (+27) & 78.364 & 17.978 & 18518.5 \\
		S3 & 90 (+30) & 58.652 & 18.050 & 14199.9 \\
		\hline 
	\end{tabular} 
	\label{tab: increase_deg}
\end{table}

As we see, $\tilde{F}_{prj}$ is much smaller with linear elements, indicating that the structure with artefacts has a greater (but artificial) stiffness. Furthermore, the objective function value varies a lot depending on the type and the degree of the elements, which makes the comparison difficult and inaccurate. After computing the corrected function value $F_{prj}$, not only do we realize that the compliance is much bigger with linear elements, but we also notice that the values became closer and more comparable with quadratic and cubic elements. In general, $F_{prj}$ is slightly smaller for the Lagrange elements than for the serendipity ones, as well as slightly smaller for elements with degree $2$. However, we emphasize that this comparison is not perfect, since the problems solved are not identical and they have many local minima. 

Although the difference in the objective function value ($F_{prj}$) is relatively small when we change from degree $2$ to degree $3$, the increase in the total time ($T_{total}$) is significant, as the size of the linear equilibrium systems also increases. Since the results shown here suggest that the solutions found with quadratic and cubic elements are considerably good, the use of elements with degree greater than $3$ should be avoided, as they are disadvantageous from the point of view of computational cost.

For a better comparison of the time spent with each type of element, we also solved the problems with a fixed budget of $30$ outer iterations of the SLP method. Table \ref{tab: increase_deg_times} presents the time spent (in seconds) on the algorithm steps affected by the type of finite element, as well as the total time ($T_{total}$). 

\begin{table}[h]
	\centering
	\caption{Time spent (in seconds) on each step for the {\tt mbb384x64x64} problem using MR with different types of elements.}
	\def\arraystretch{1.2}
	\setlength{\tabcolsep}{8pt}
	\small
	\begin{tabular}{l|rrrrr}
		\hline 
		\multicolumn{1}{c|}{Step} & \multicolumn{1}{c}{L1} & \multicolumn{1}{c}{L2} & \multicolumn{1}{c}{S2} & \multicolumn{1}{c}{L3} & \multicolumn{1}{c}{S3}\\ \hline 
		Assembly of $\mathrm{K}$ &  5.5 &  65.7 &  30.0 &  642.9 &  103.2 \\
		Precond. setup          &  0.6 &  19.8 &  12.7 &  148.9 &   38.4 \\
		Linear systems          & 16.2 & 393.9 & 736.0 & 5669.2 & 3945.7 \\
		Gradients               &  8.1 &  22.4 &  14.3 &  355.5 &   30.5 \\
		\hline 
		Total                   & 98.7 & 563.8 & 857.5 & 6885.2 & 4182.0 \\
		\hline 
	\end{tabular} 
	\label{tab: increase_deg_times}
\end{table}

Comparing the elements of the Lagrange and serendipity families with the same degree, we notice that serendipity elements require less time to assemble matrix $\mathrm{K}$, set up the preconditioner and calculate the gradient. This was expected given that the serendipity elements have less nodes than the Lagrange elements with the same degree. On the other hand, the time spent on the solution of the linear systems and, consequently, the total time are greater with quadratic serendipity elements and smaller only with cubic elements.

\begin{figure}[h]
	\centering
	\includegraphics[width=0.45\linewidth]{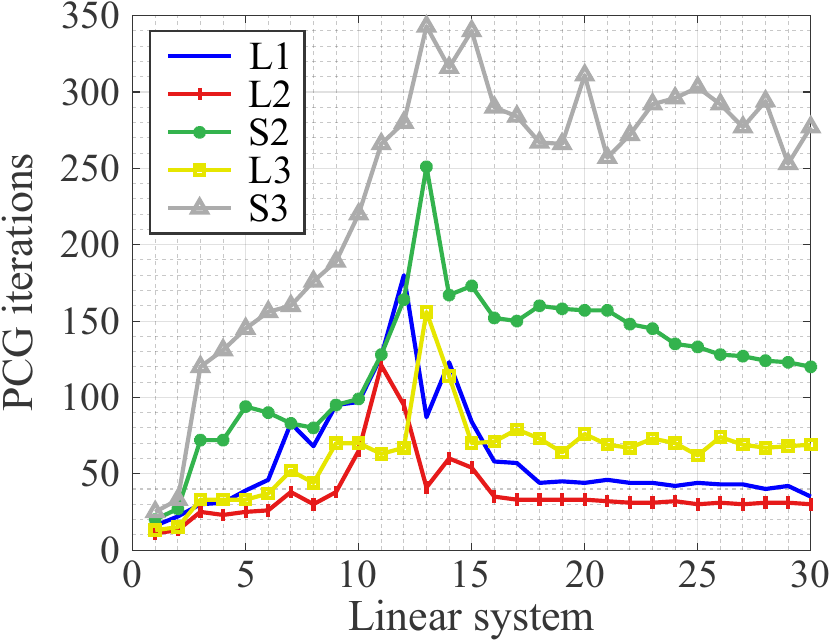}
	\caption{Number of PCG iterations required to solve the linear systems for the {\tt mbb384x64x64} problem, with different types of elements.}
	\label{fig: PCG_deg}
\end{figure}

To better understand how the solution of linear systems is affected by the element type, we show in Figure \ref{fig: PCG_deg} the number of PCG iterations required to solve the linear system at each iteration of the SLP algorithm, for the MBB beam problem.

We observe that elements of the Lagrange family require less PCG iterations than the serendipity ones. We believe that this occurs because the geometric multigrid method is more effective with Lagrange elements. The nodes in this type of element are equally spaced, so that the prolongation operator can be constructed in the usual way, through a linear interpolation, using the surrounding coarse mesh nodes as interpolating points. For serendipity elements, the position of the nodes makes the prolongation more complicated. In this case, the prolongation operator is constructed using interpolation via the shape functions of the finite element method. Although this strategy worked well, the multigrid method did not achieve the same performance as for the Lagrange elements. 

For serendipity elements of degree 2, despite the size of the linear systems being smaller, the increase in the number of PCG iterations made the time taken to solve them to be greater than with Lagrange elements. For cubic serendipity elements, on the other hand, the reduction on the size of the systems is so significant that the time is reduced even with the increase in the number of iterations. Based on these results, we conclude that, when the geometric multigrid method is employed, it is more advantageous to use elements of the Lagrange family with degree at most $2$ and elements of the serendipity family with degree greater than or equal to $3$.

\subsubsection{Adaptive strategy to increase the degree of the elements}

As noted in the previous subsection, increasing the degree of the finite elements improves the quality of the solutions of the topology optimization problem, but also increases the computational cost. Now, we test the adaptive strategy proposed in Section \ref{sec: fix}, whose purpose is to obtain more detailed and accurate solutions without harming too much the algorithm efficiency.
To avoid numerical errors that may arise when we suppress the nodal displacements in void regions, we use $E_{min} = 10^{-9}E_0$ as the Young's modulus of the void. Moreover, we use the thresholds $\rho_{fix}^l = 10^{-6}$ and $\rho_{fix}^u = 0.9$ and $\varepsilon_{grad} = 10^{-6}$ as described in Section \ref{sec: fix}.

We decided to apply our adaptive strategy to the solution of more challenging tests, such as the MBB beam and the cantilever beam problems with a volume limit of $15\%$ of the total domain volume, as well as the L-shaped beam and the bridge problems. For the first three problems, we use multiresolution with $n_{mr} = 4$ and $d_{mr} = 2$. For the bridge, we set $n_{mr} = 3$ and $d_{mr} = 2$. All these structures can be considered large-scale problems, since the density meshes have more than one million elements.

\begin{figure*}[h]
	\begin{subfigure}{0.5\linewidth}
		\centering
		\includegraphics[width=0.9\linewidth]{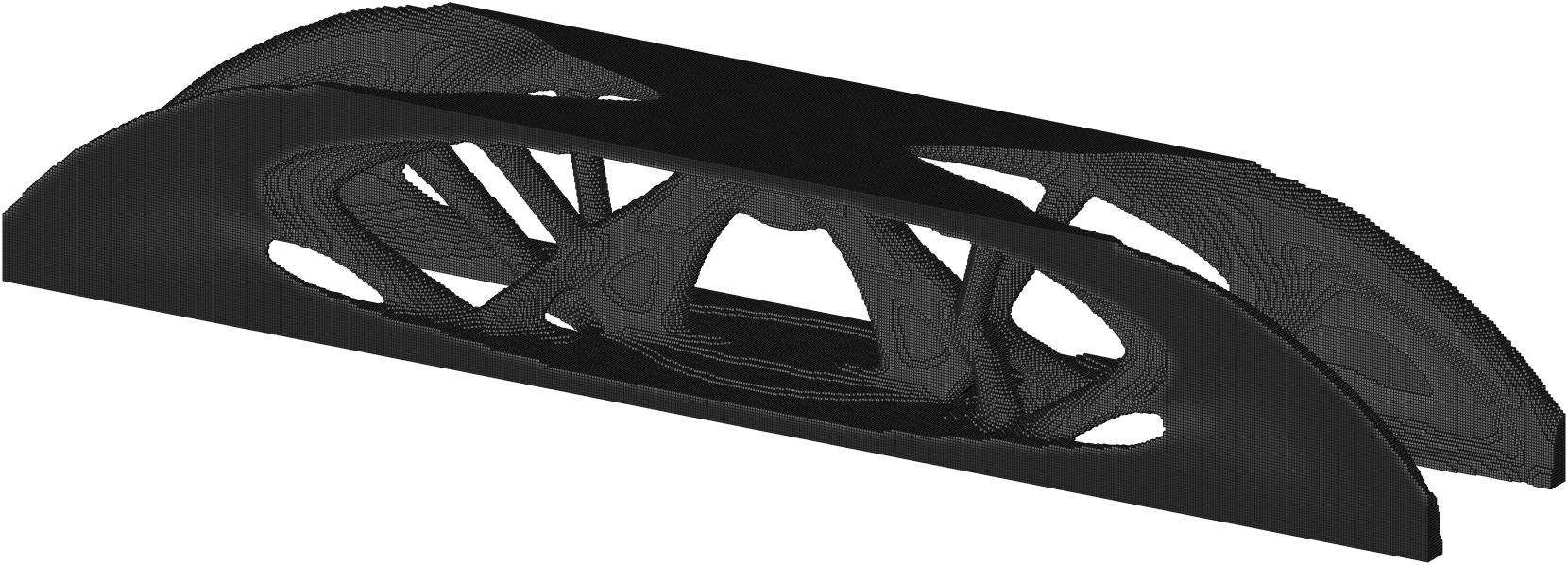}
	\end{subfigure}%
	\begin{subfigure}{0.5\linewidth}
		\centering
		\includegraphics[width=0.9\linewidth]{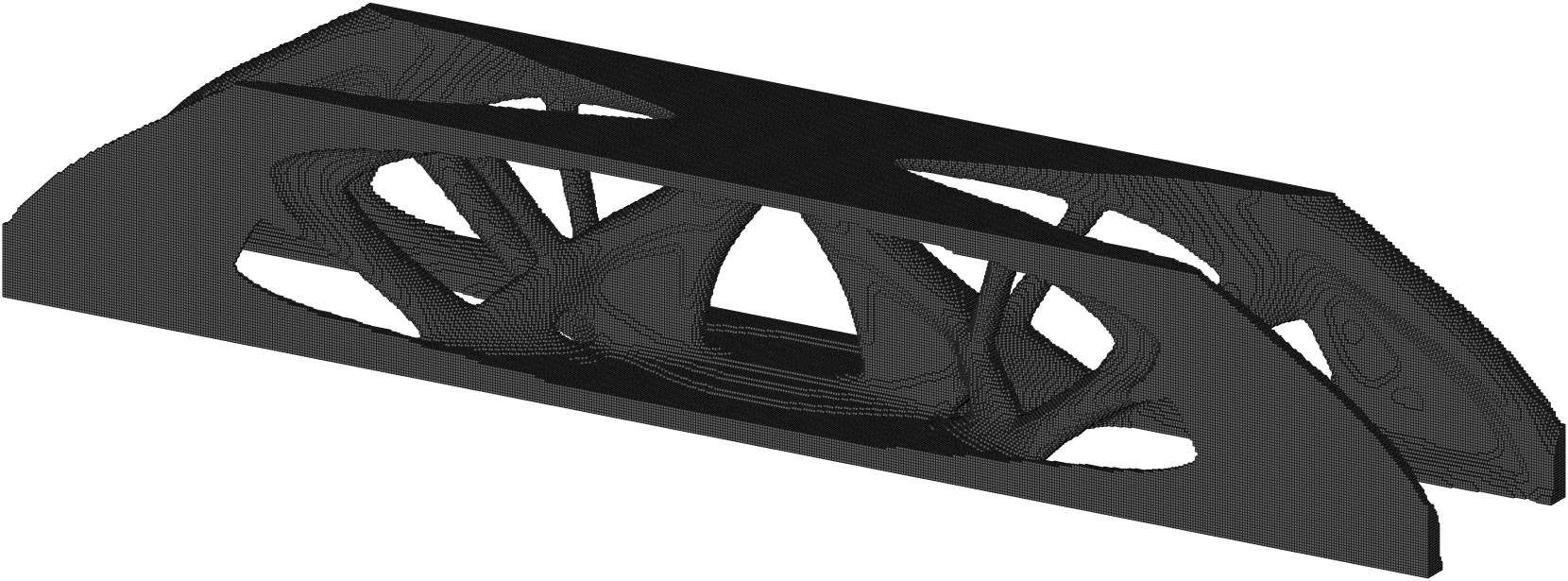}
	\end{subfigure}
	\begin{subfigure}{0.5\linewidth}
		\centering
		\includegraphics[width=0.9\linewidth]{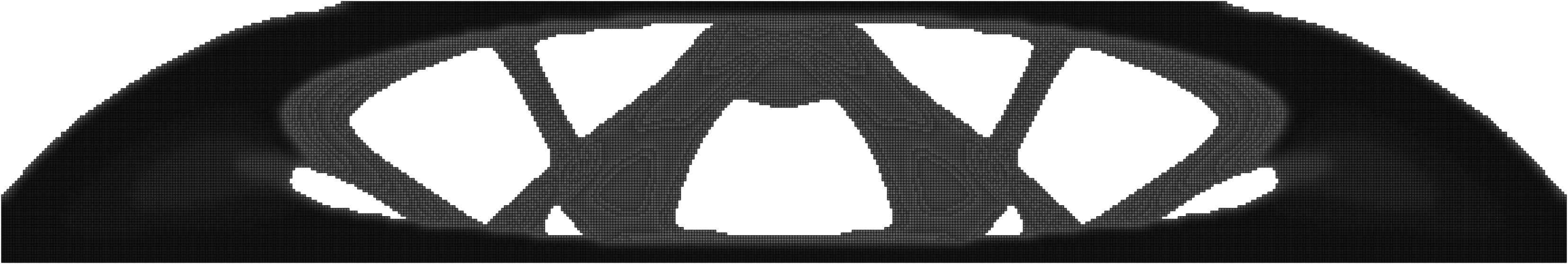}
		\caption{Linear elements ($r_{min} = 1.5$).}
	\end{subfigure}%
	\begin{subfigure}{0.5\linewidth}
		\centering
		\includegraphics[width=0.9\linewidth]{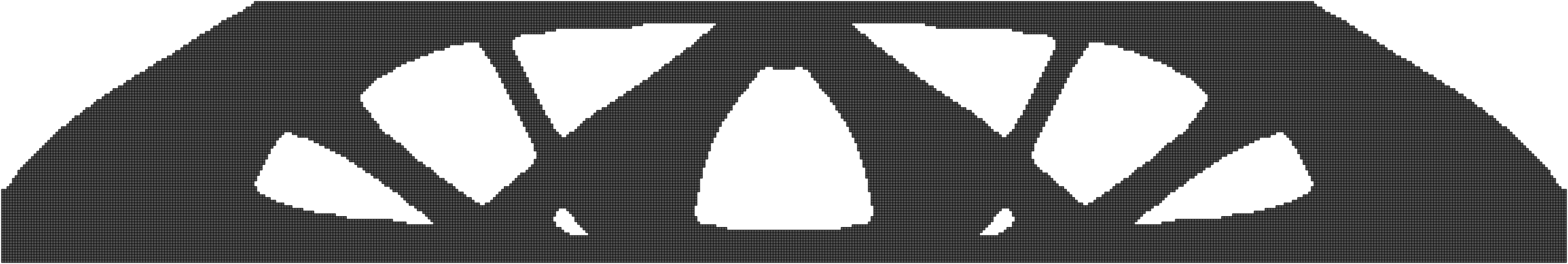}
		\caption{Adaptive strategy with maximum degree $2$ ($r_{min} = 1.5$).}
	\end{subfigure}
	
	\begin{subfigure}{0.5\linewidth}
		\centering
		\includegraphics[width=0.9\linewidth]{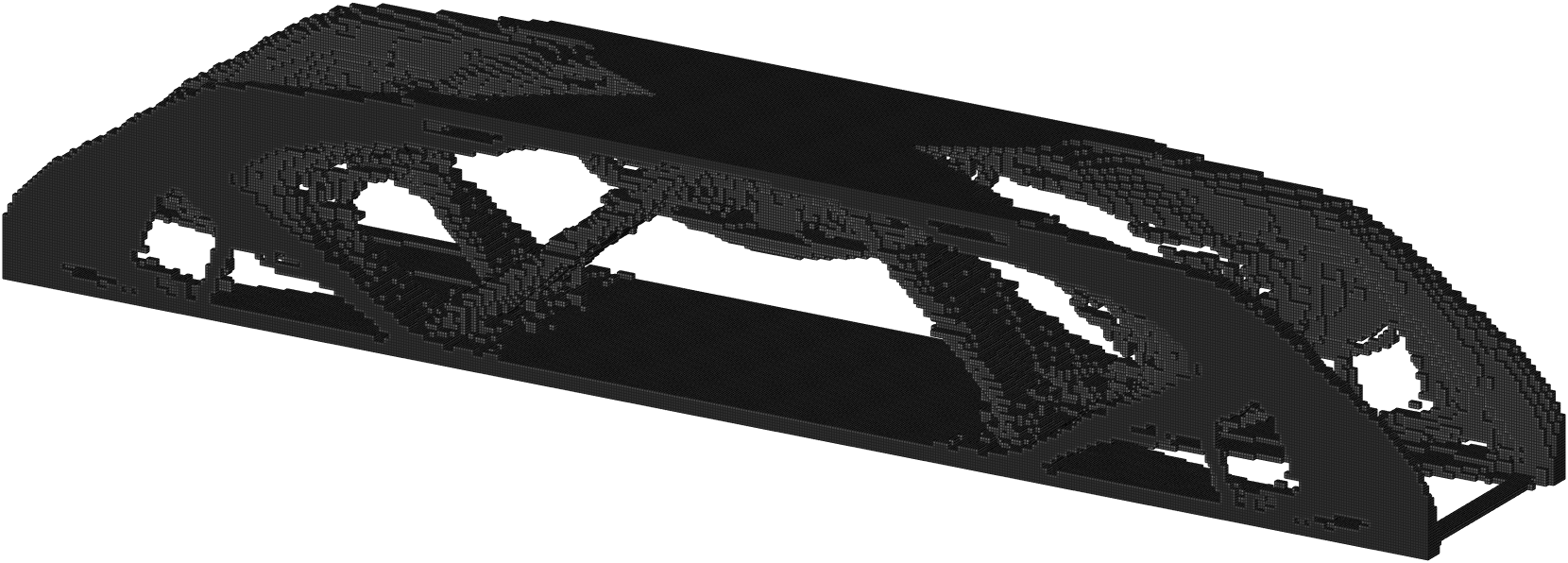}
	\end{subfigure}%
	\begin{subfigure}{0.5\linewidth}
		\centering
		\includegraphics[width=0.9\linewidth]{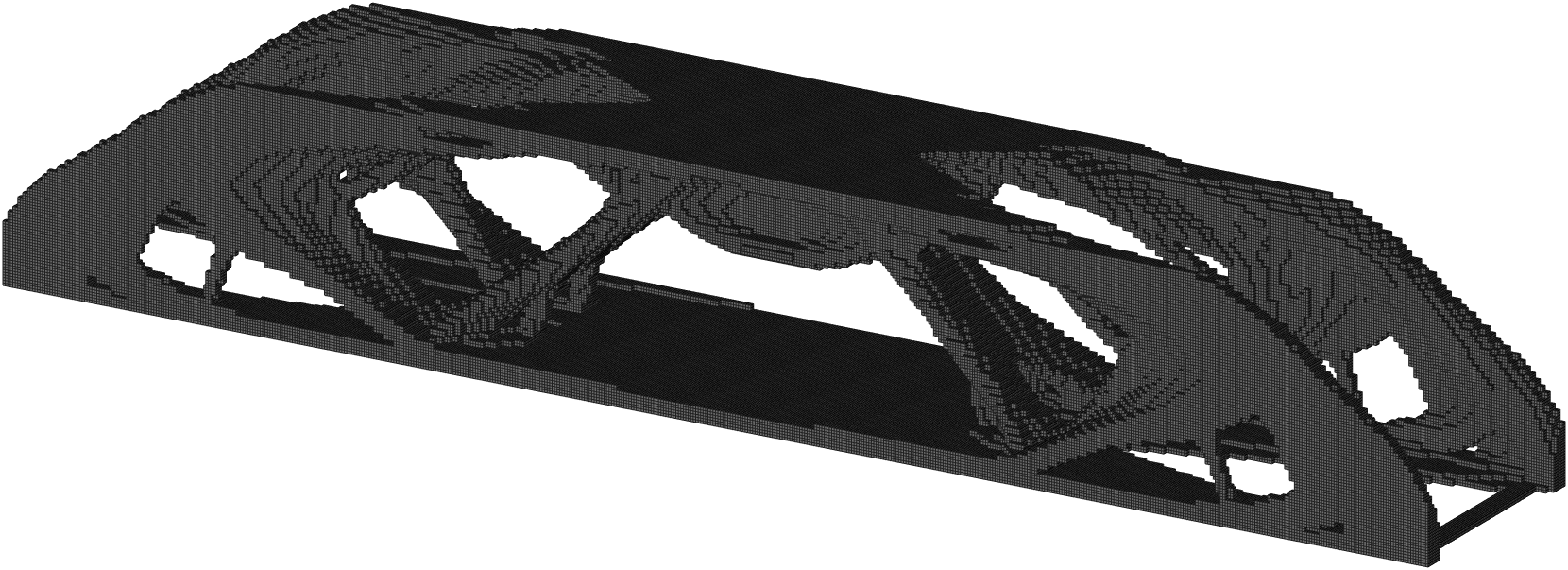}
	\end{subfigure}
	\begin{subfigure}{0.5\linewidth}
		\centering
		\includegraphics[width=0.9\linewidth]{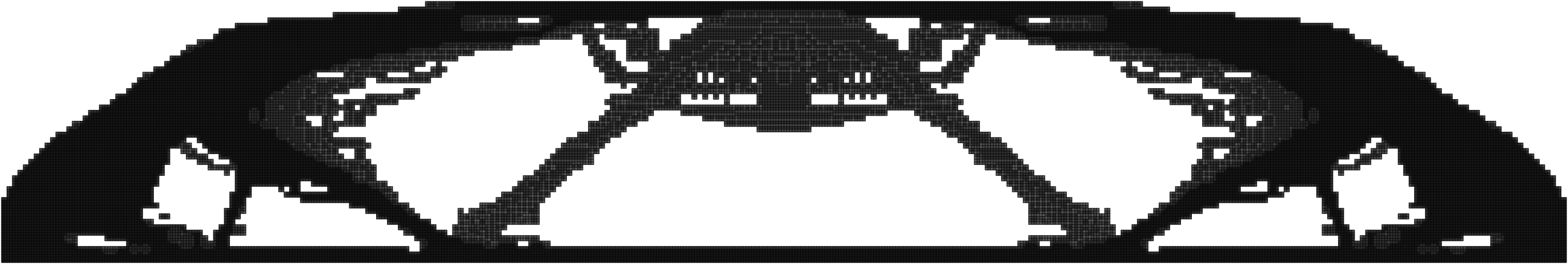}
		\caption{Linear elements ($r_{min} = 0.6$).}
	\end{subfigure}%
	\begin{subfigure}{0.5\linewidth}
		\centering
		\includegraphics[width=0.9\linewidth]{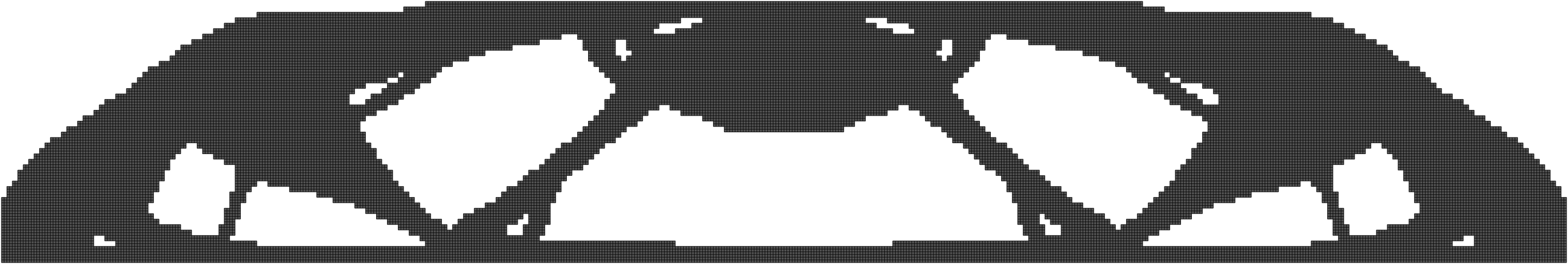}
		\caption{Adaptive strategy with maximum degree $2$ ($r_{min} = 0.6$).}
	\end{subfigure}
	
	\caption{Solutions of the {\tt mbb576x96x96} problem with $15\%$ volume, using MR with $n_{mr} = 4$ and $d_{mr} = 2$.}
	\label{fig: mbb_fix}
\end{figure*}

\begin{figure*}[h]
	\begin{subfigure}{0.25\linewidth}
		\centering
		\includegraphics[width=0.95\linewidth]{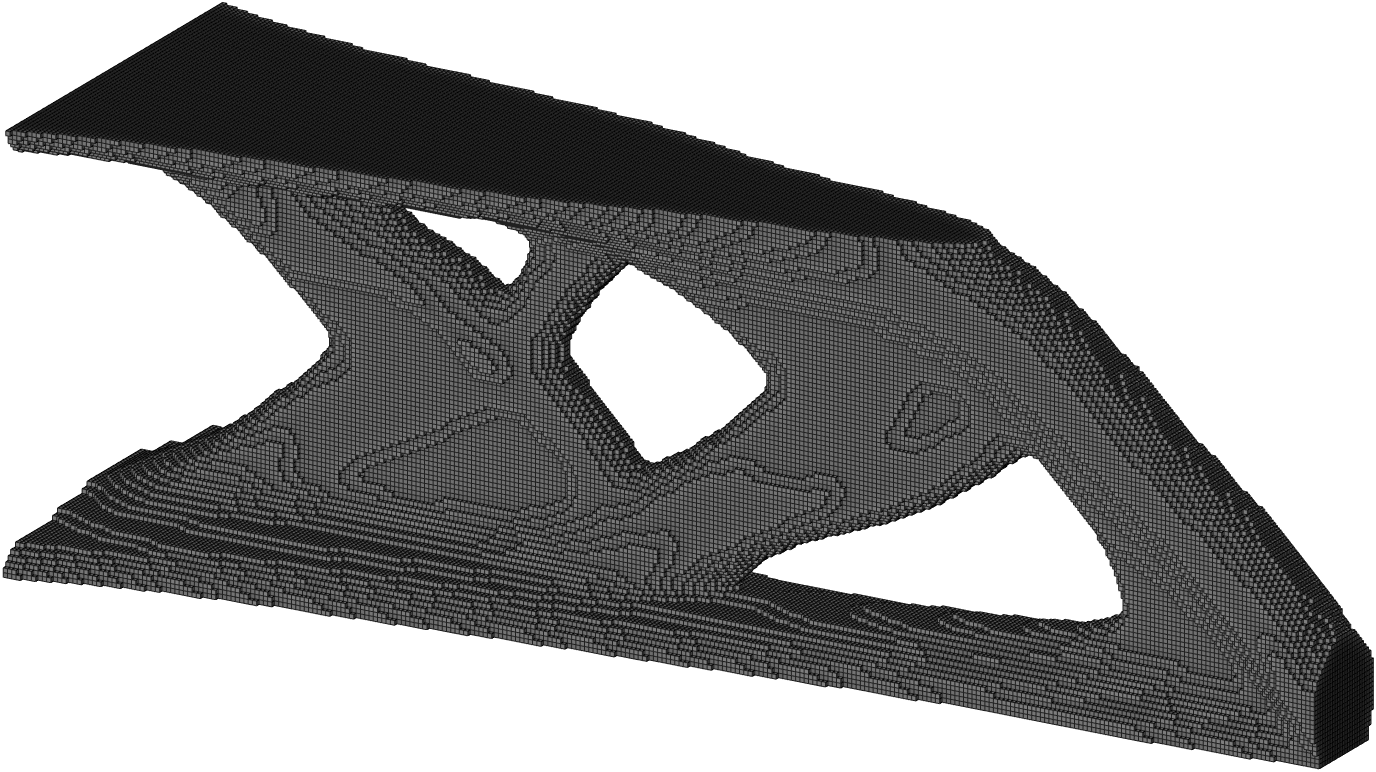}
	\end{subfigure}%
	\begin{subfigure}{0.25\linewidth}
		\centering
		\includegraphics[width=0.95\linewidth]{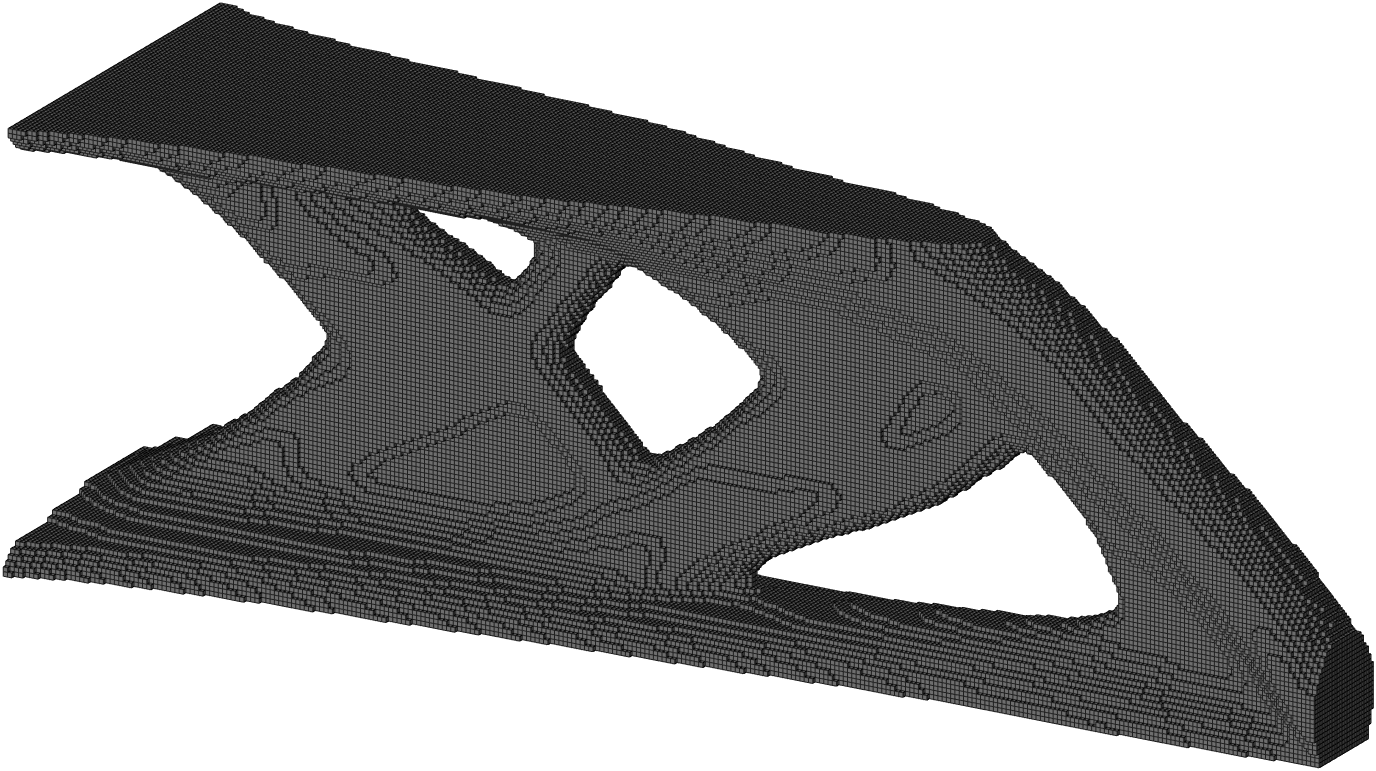}
	\end{subfigure}%
	\begin{subfigure}{0.25\linewidth}
		\centering
		\includegraphics[width=0.95\linewidth]{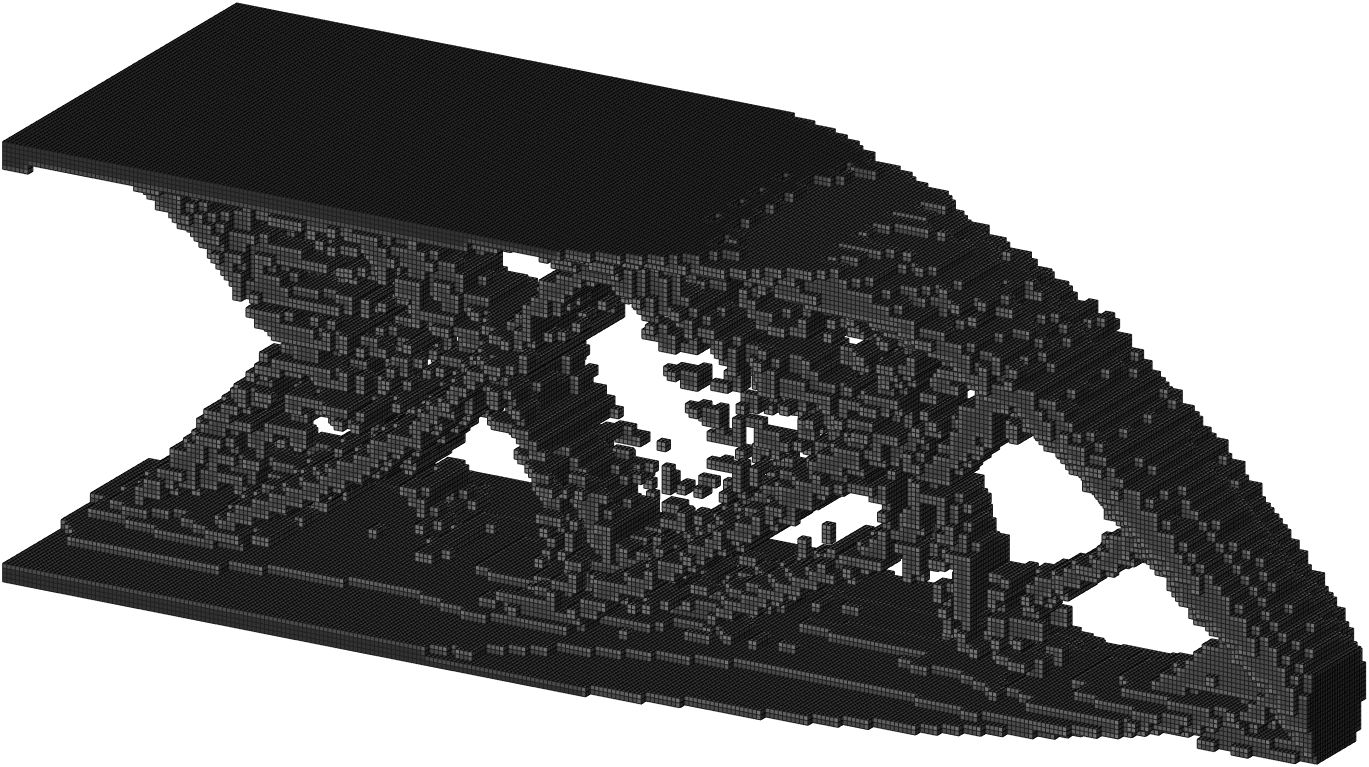}
	\end{subfigure}%
	\begin{subfigure}{0.25\linewidth}
		\centering
		\includegraphics[width=0.95\linewidth]{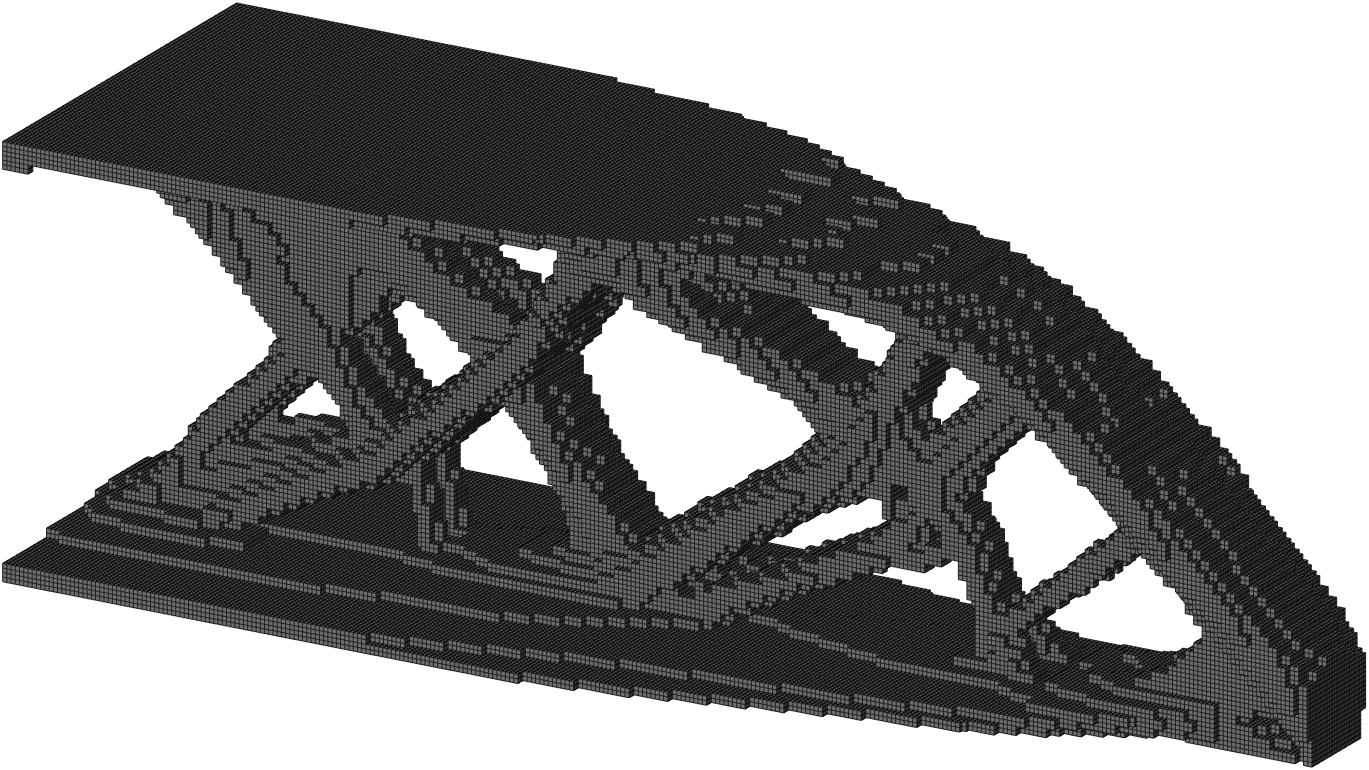}
	\end{subfigure}
	
	\begin{subfigure}{0.25\linewidth}
		\centering
		\includegraphics[width=0.95\linewidth]{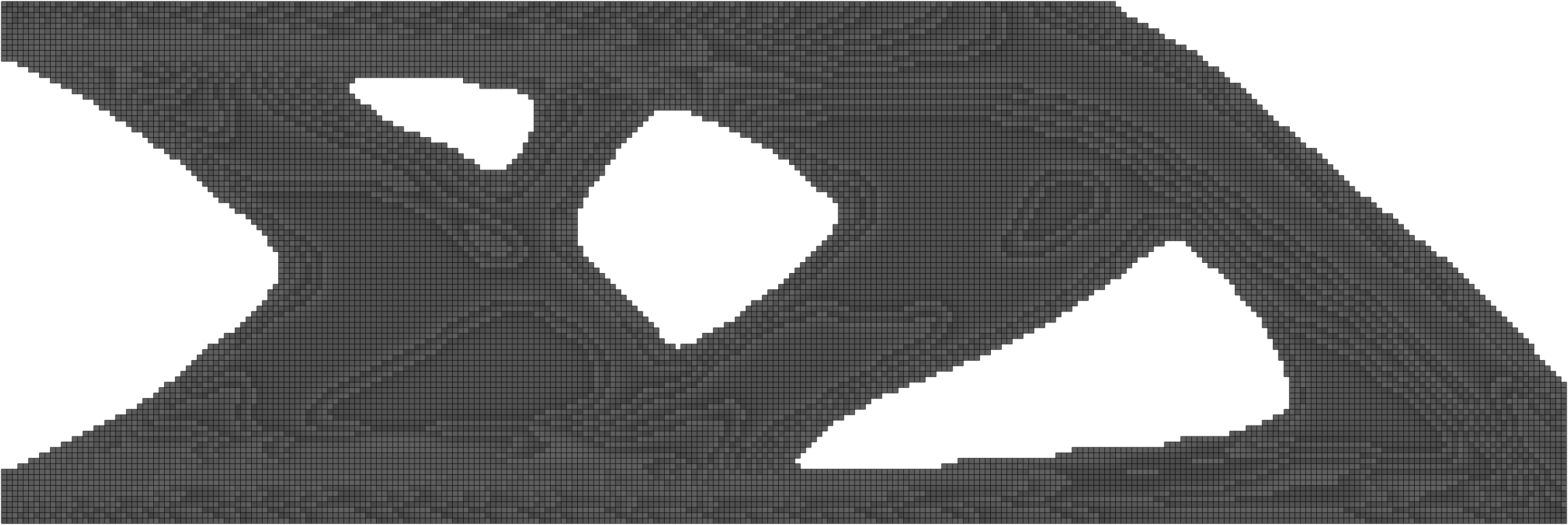}
		\caption{Linear elements \\ ($r_{min} = 1.5$). \centering}
	\end{subfigure}%
	\begin{subfigure}{0.25\linewidth}
		\centering
		\includegraphics[width=0.95\linewidth]{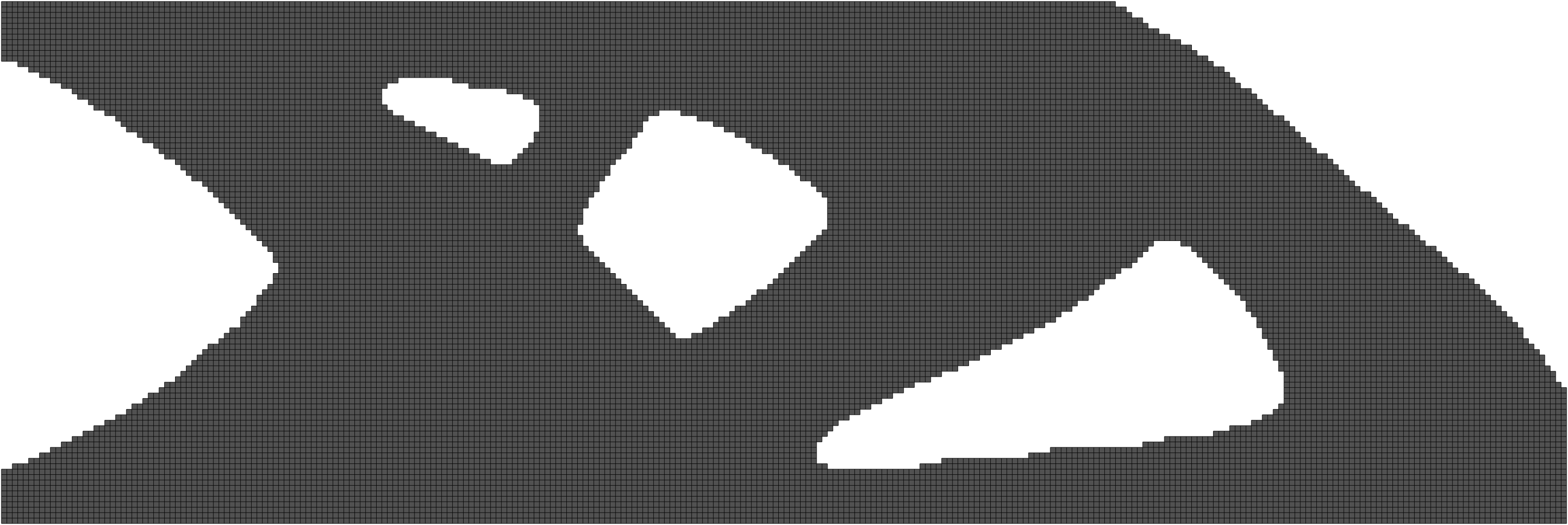}
		\caption{Adaptive strategy with max. degree $2$ ($r_{min} = 1.5$). \centering}
	\end{subfigure}%
	\begin{subfigure}{0.25\linewidth}
		\centering
		\includegraphics[width=0.95\linewidth]{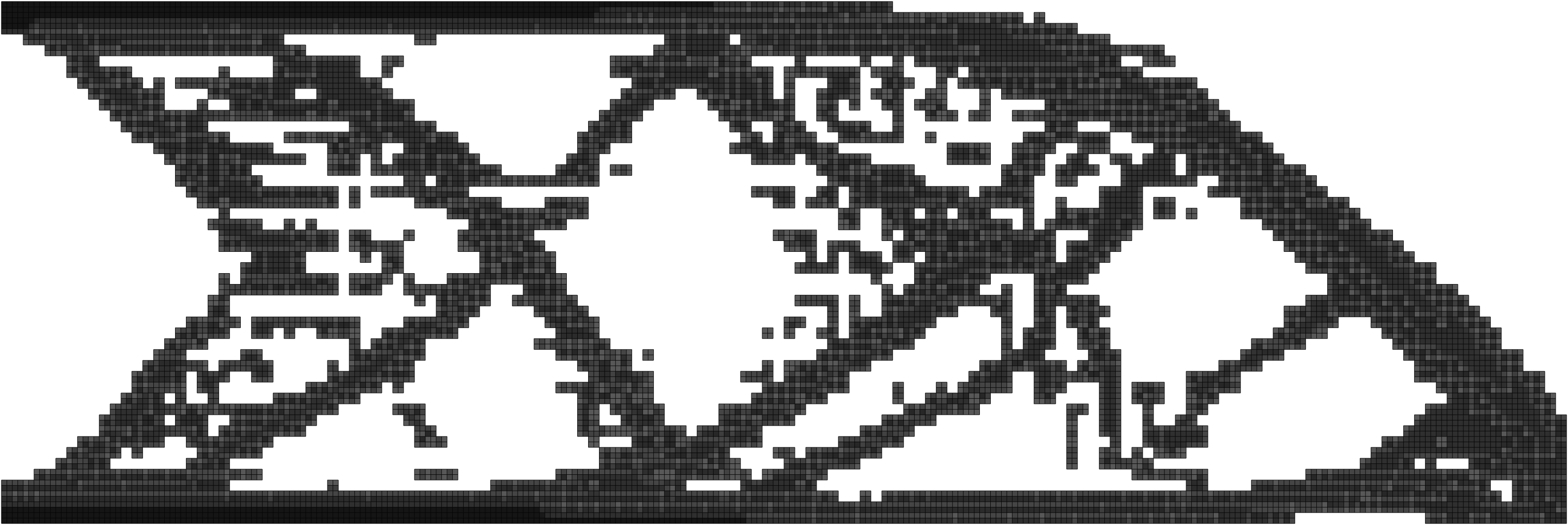}
		\caption{Linear elements \\ ($r_{min} = 0.6$). \centering}
	\end{subfigure}%
	\begin{subfigure}{0.25\linewidth}
		\centering
		\includegraphics[width=0.95\linewidth]{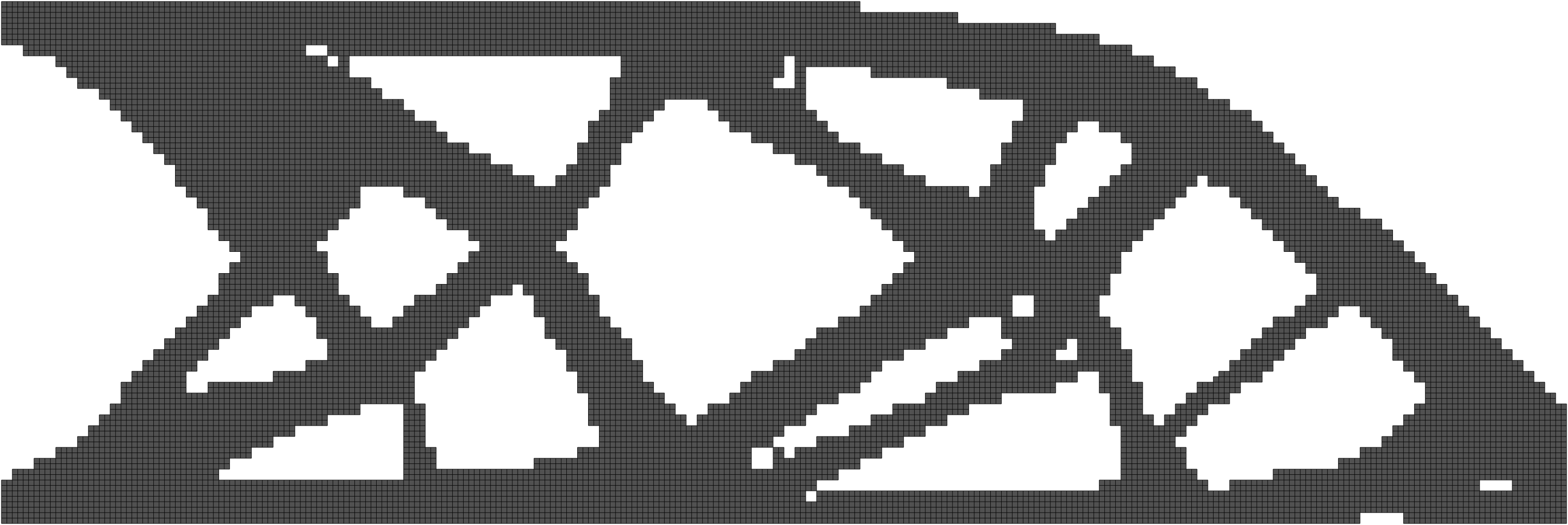}
		\caption{Adaptive strategy with max. degree $2$ ($r_{min} = 0.6$). \centering}
	\end{subfigure}
	
	\caption{Solutions of the {\tt cb192x64x64} problem with $15\%$ volume, using MR with $n_{mr} = 4$ and $d_{mr} = 2$.}
	\label{fig: cb_fix}
\end{figure*}

\begin{figure*}[h]
	\begin{subfigure}{0.25\linewidth}
		\centering
		\includegraphics[width=0.8\linewidth]{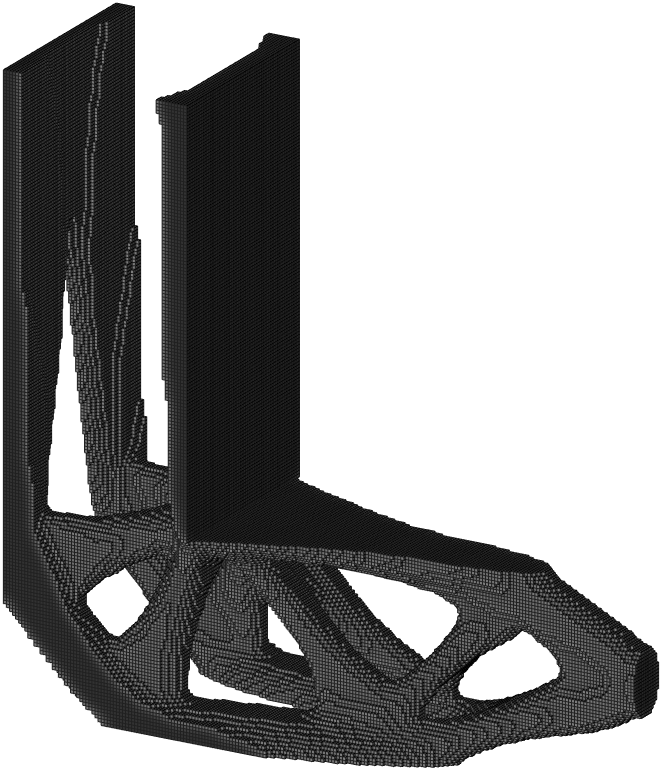}
	\end{subfigure}%
	\begin{subfigure}{0.25\linewidth}
		\centering
		\includegraphics[width=0.8\linewidth]{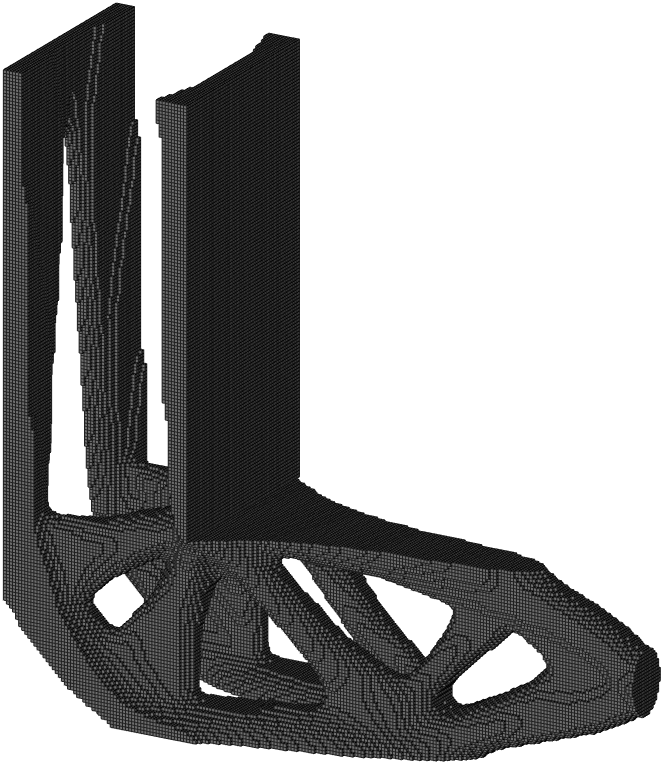}
	\end{subfigure}%
	\begin{subfigure}{0.25\linewidth}
		\centering
		\includegraphics[width=0.8\linewidth]{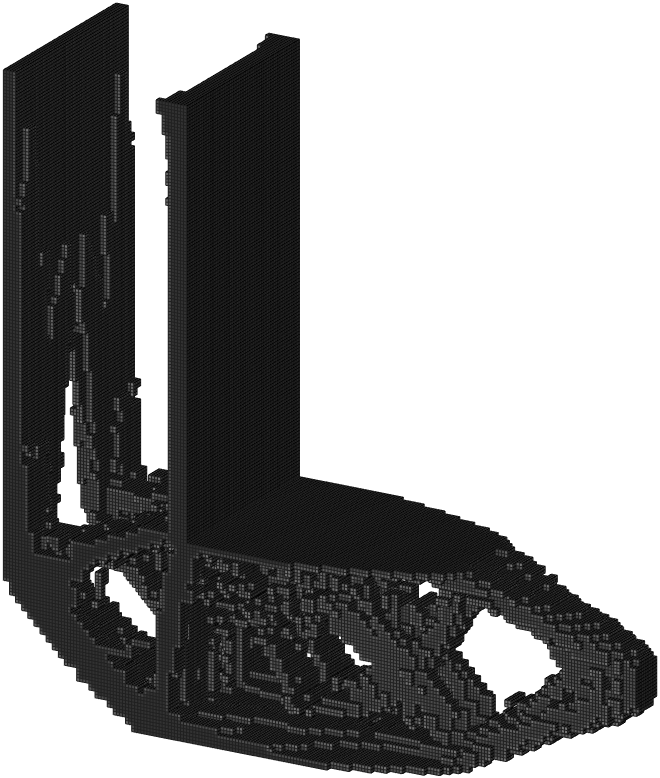}
	\end{subfigure}%
	\begin{subfigure}{0.25\linewidth}
		\centering
		\includegraphics[width=0.8\linewidth]{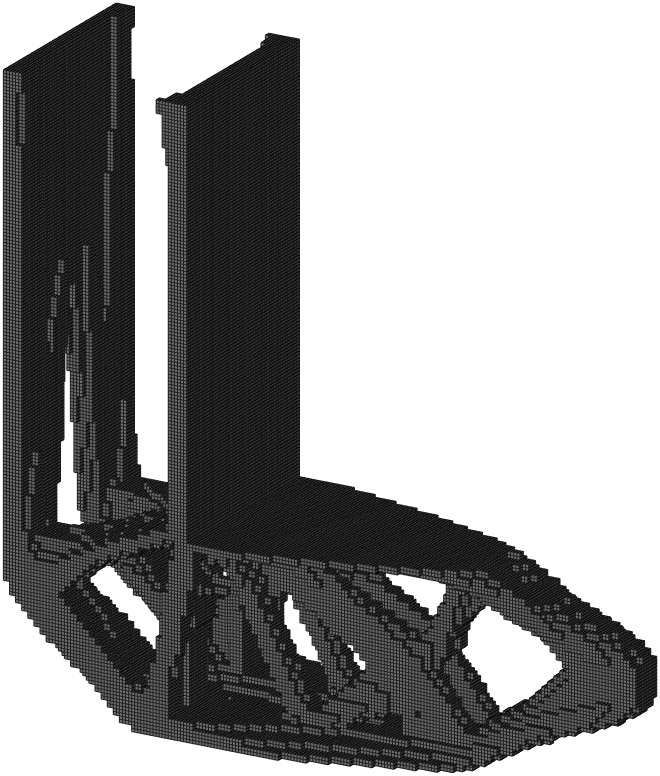}
	\end{subfigure}
	
	\begin{subfigure}{0.25\linewidth}
		\centering
		\includegraphics[width=0.8\linewidth]{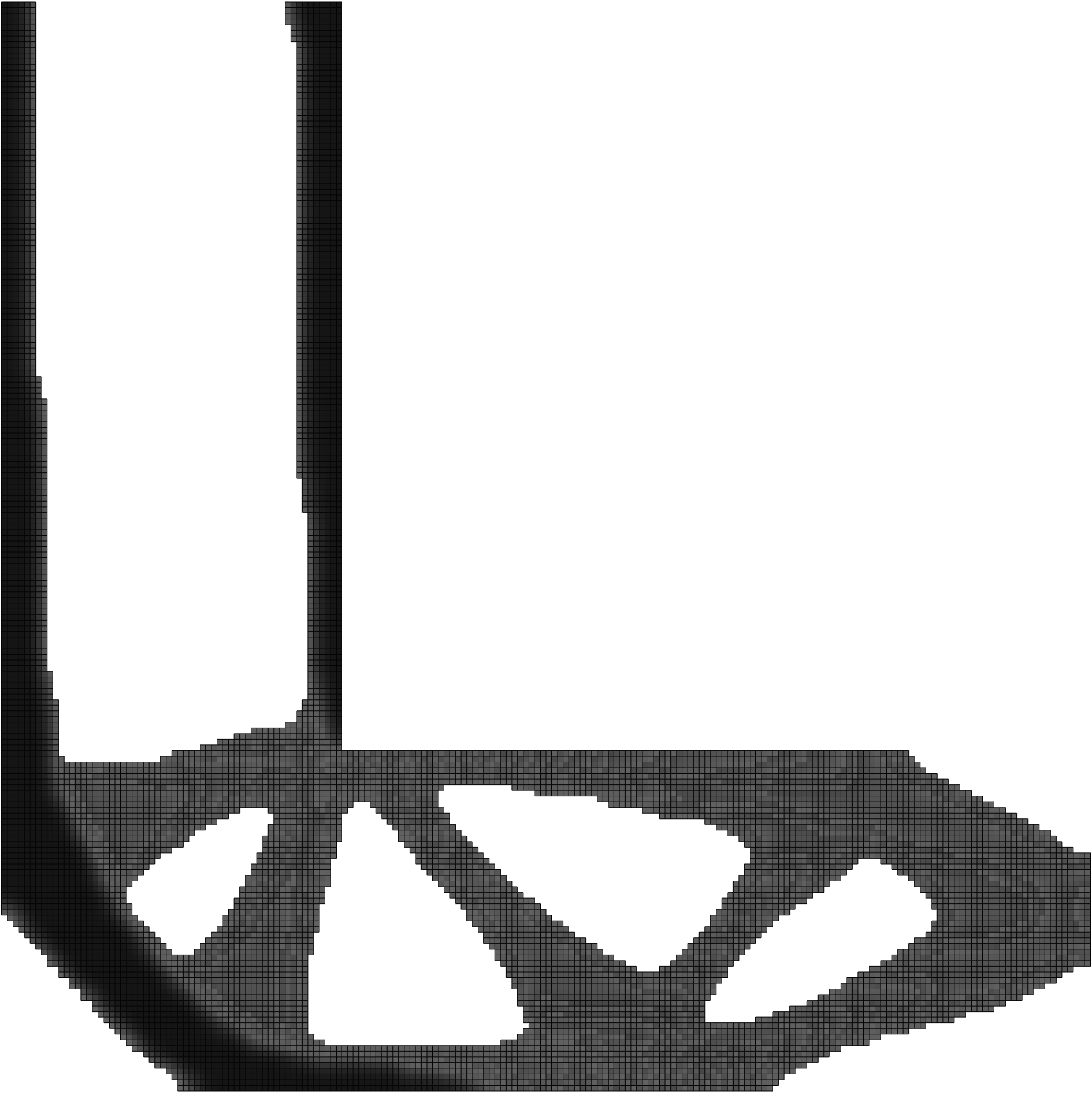}
		\caption{Linear elements \\ ($r_{min} = 1.5$). \centering}
	\end{subfigure}%
	\begin{subfigure}{0.25\linewidth}
		\centering
		\includegraphics[width=0.8\linewidth]{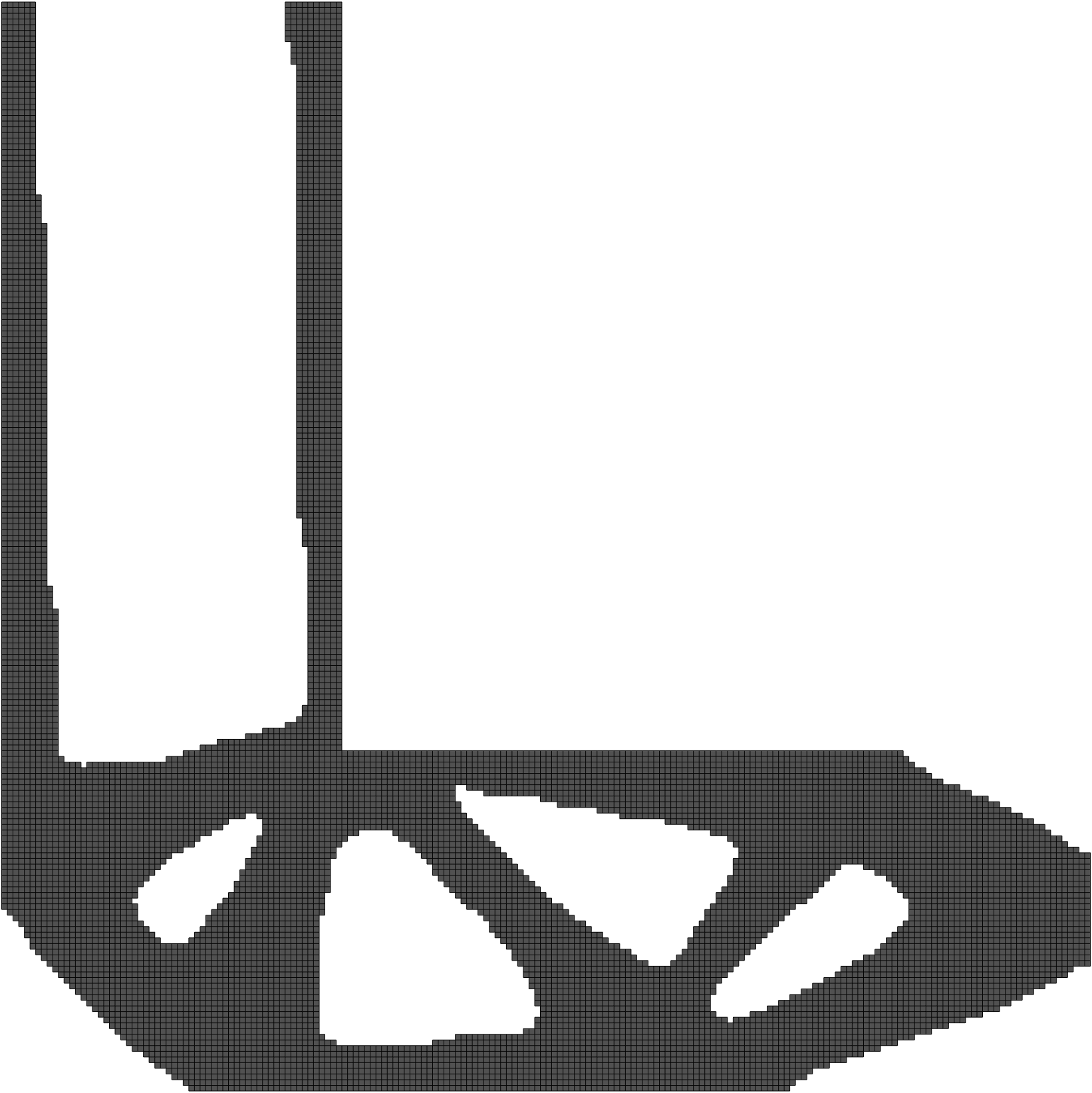}
		\caption{Adaptive strategy with max. degree $2$ ($r_{min} = 1.5$). \centering}
	\end{subfigure}%
	\begin{subfigure}{0.25\linewidth}
		\centering
		\includegraphics[width=0.8\linewidth]{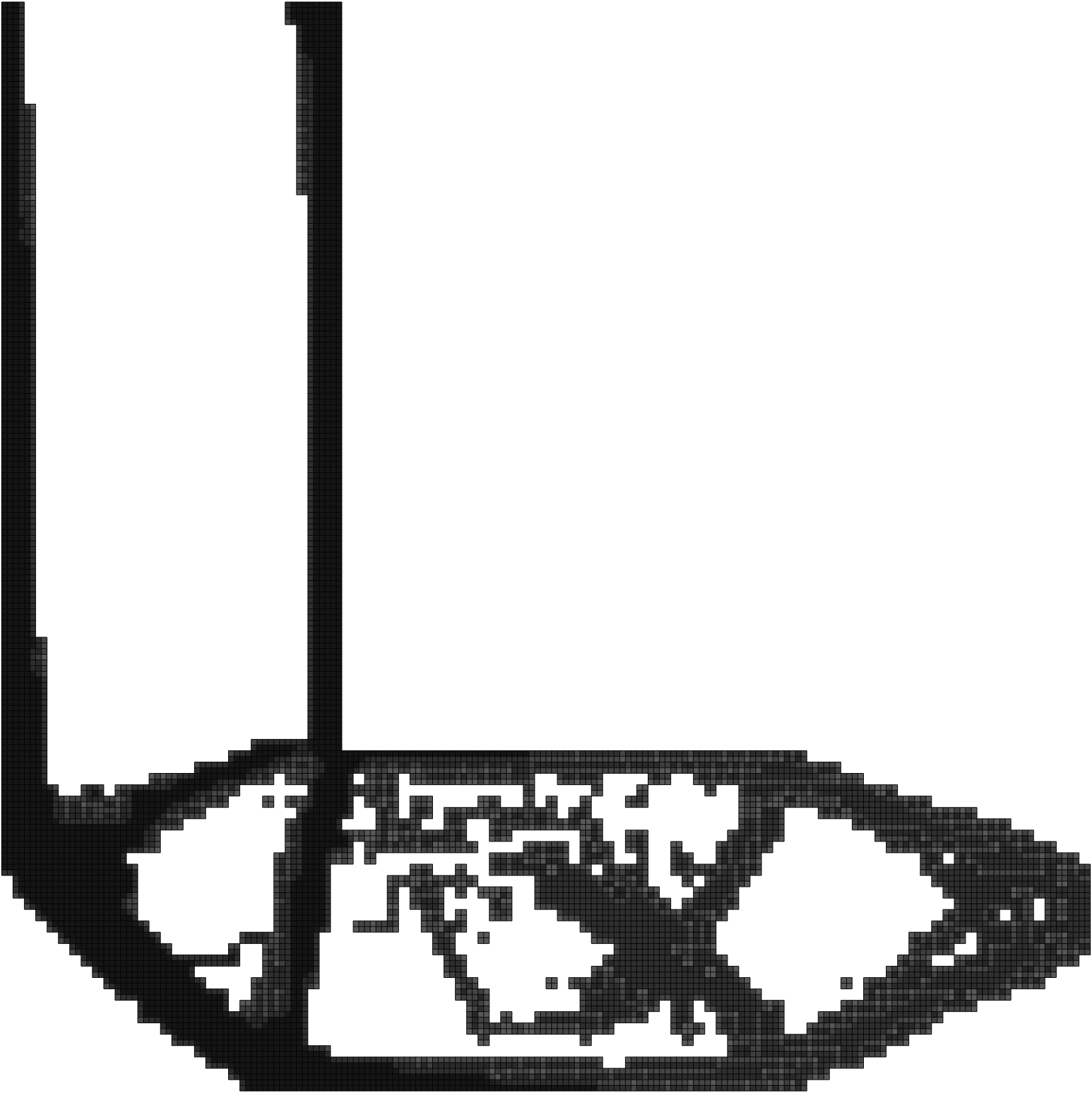}
		\caption{Linear elements \\ ($r_{min} = 0.6$). \centering}
	\end{subfigure}%
	\begin{subfigure}{0.25\linewidth}
		\centering
		\includegraphics[width=0.8\linewidth]{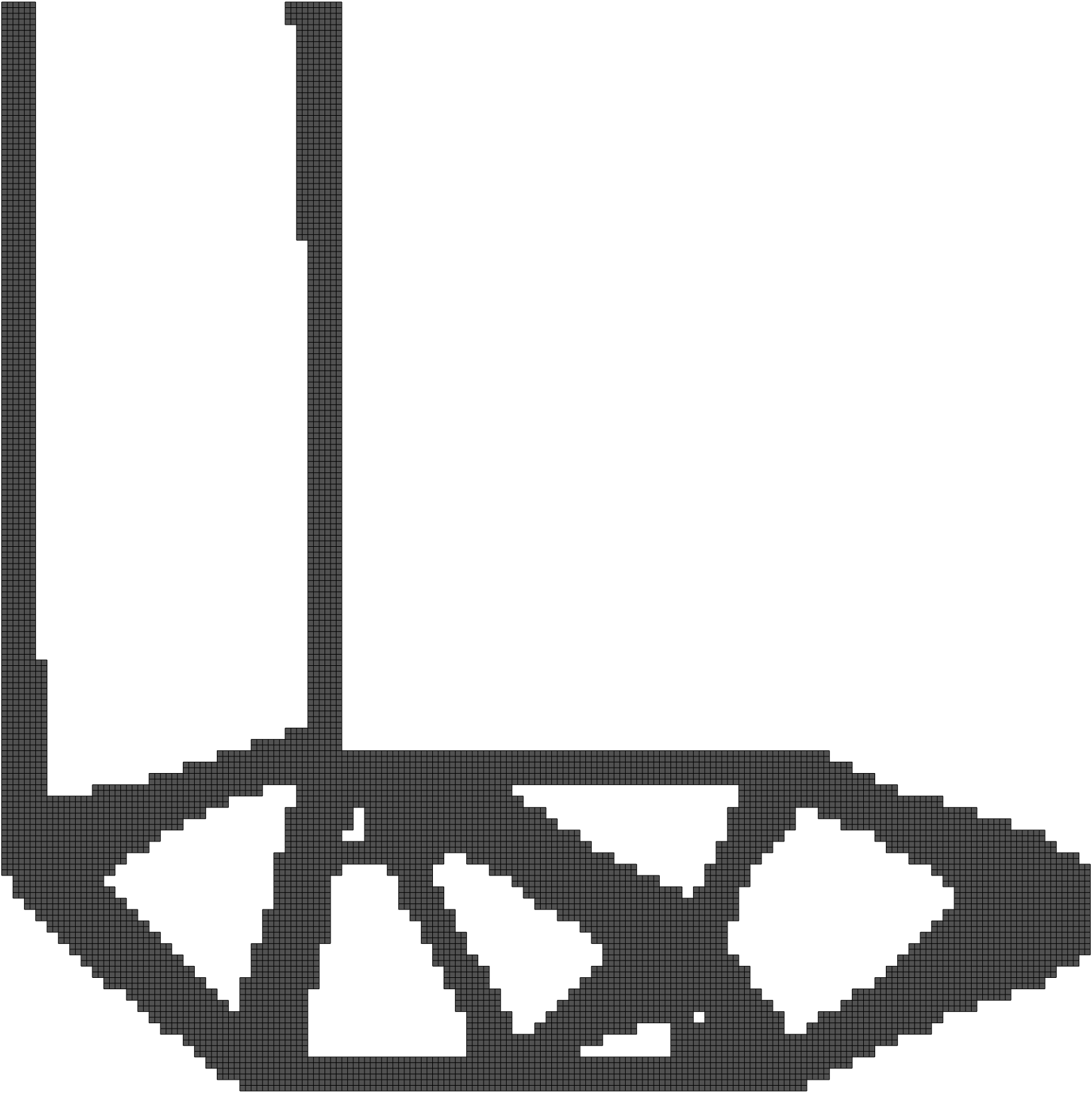}
		\caption{Adaptive strategy with max. degree $2$ ($r_{min} = 0.6$). \centering}
	\end{subfigure}
	
	\caption{Solutions of the {\tt ls192x192x64} problem with $18\%$ volume, using MR with $n_{mr} = 4$ and $d_{mr} = 2$.}
	\label{fig: ls_fix}
\end{figure*}

\begin{figure*}[h]
	\begin{subfigure}{0.5\linewidth}
		\centering
		\includegraphics[width=0.9\linewidth]{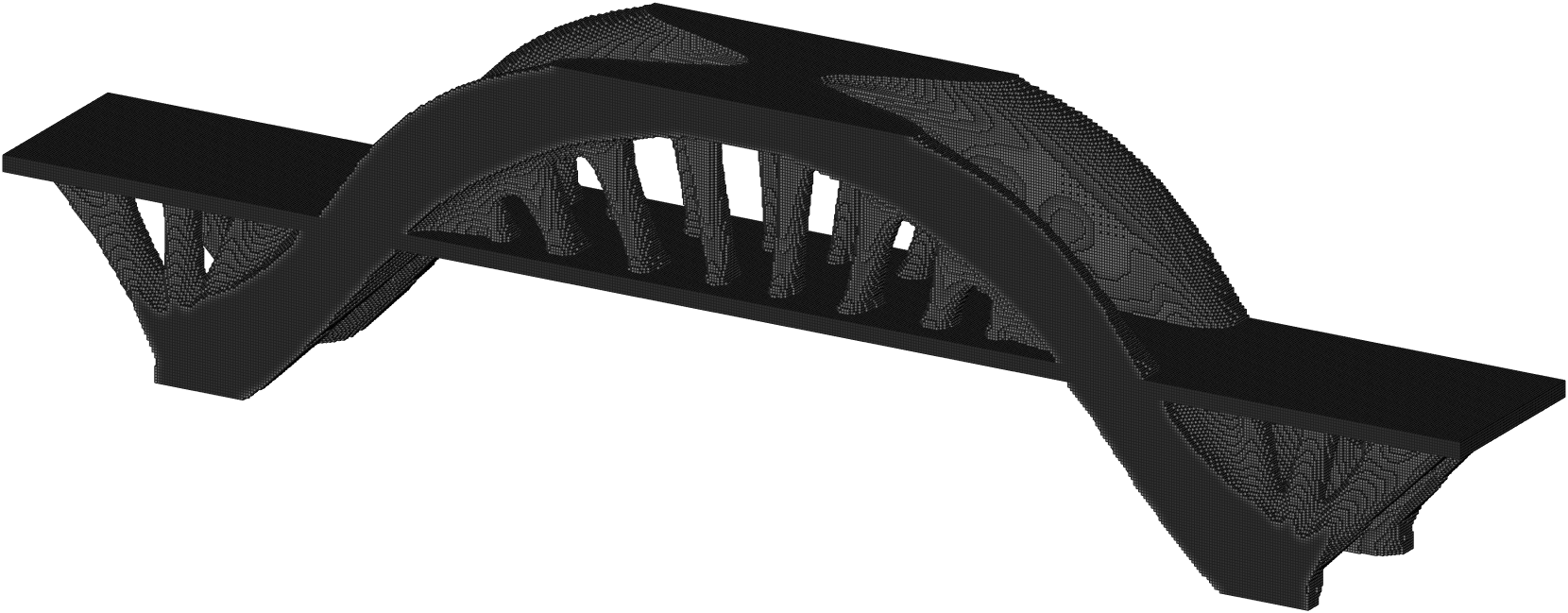}
	\end{subfigure}%
	\begin{subfigure}{0.5\linewidth}
		\centering
		\includegraphics[width=0.9\linewidth]{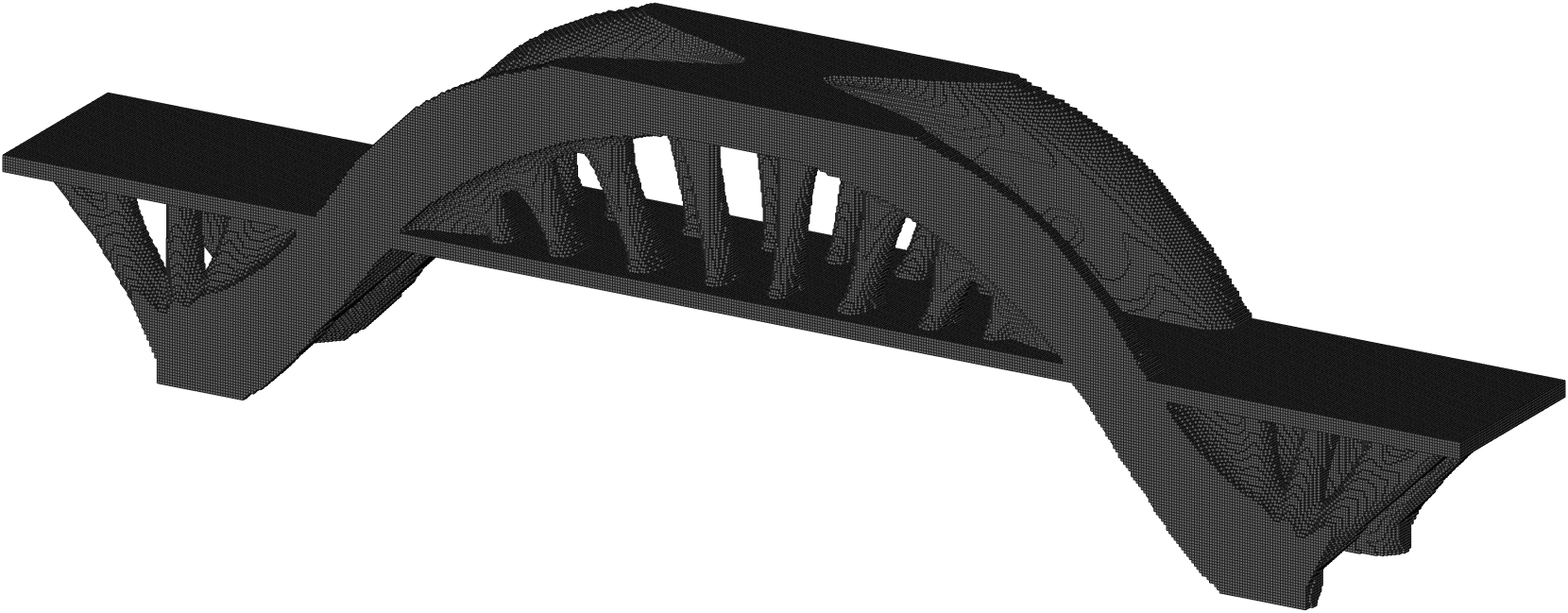}
	\end{subfigure}
	\begin{subfigure}{0.5\linewidth}
		\centering
		\includegraphics[width=0.9\linewidth]{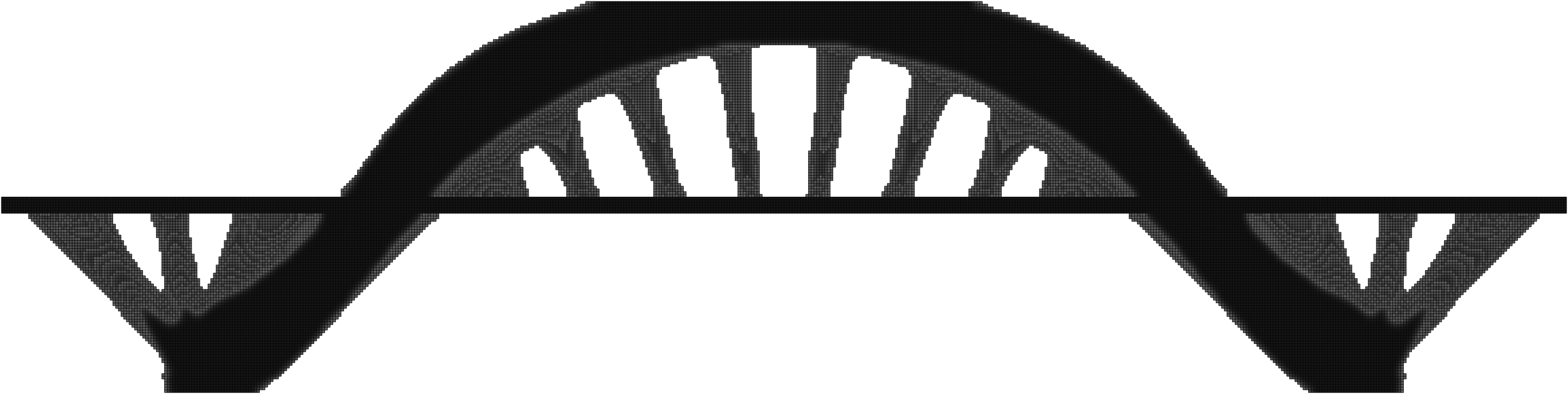}
		\caption{Linear elements ($r_{min} = 1.5$).}
	\end{subfigure}%
	\begin{subfigure}{0.5\linewidth}
		\centering
		\includegraphics[width=0.9\linewidth]{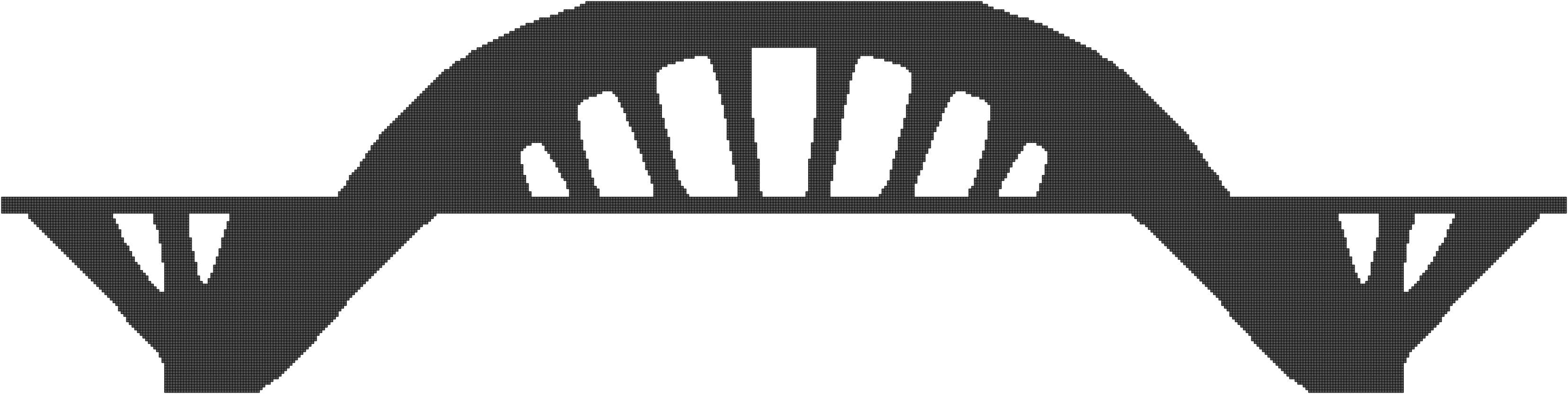}
		\caption{Adaptive strategy with maximum degree $2$ ($r_{min} = 1.5$).}
	\end{subfigure}
	
	\begin{subfigure}{0.5\linewidth}
		\centering
		\includegraphics[width=0.9\linewidth]{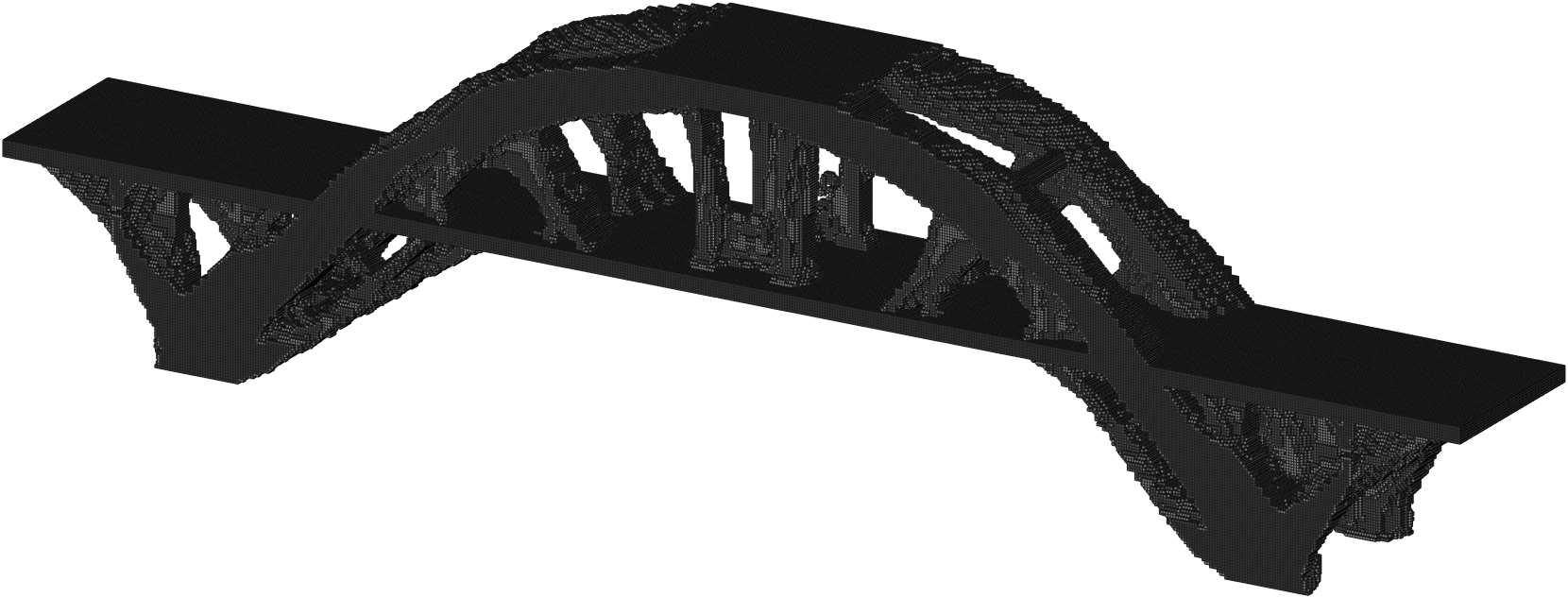}
	\end{subfigure}%
	\begin{subfigure}{0.5\linewidth}
		\centering
		\includegraphics[width=0.9\linewidth]{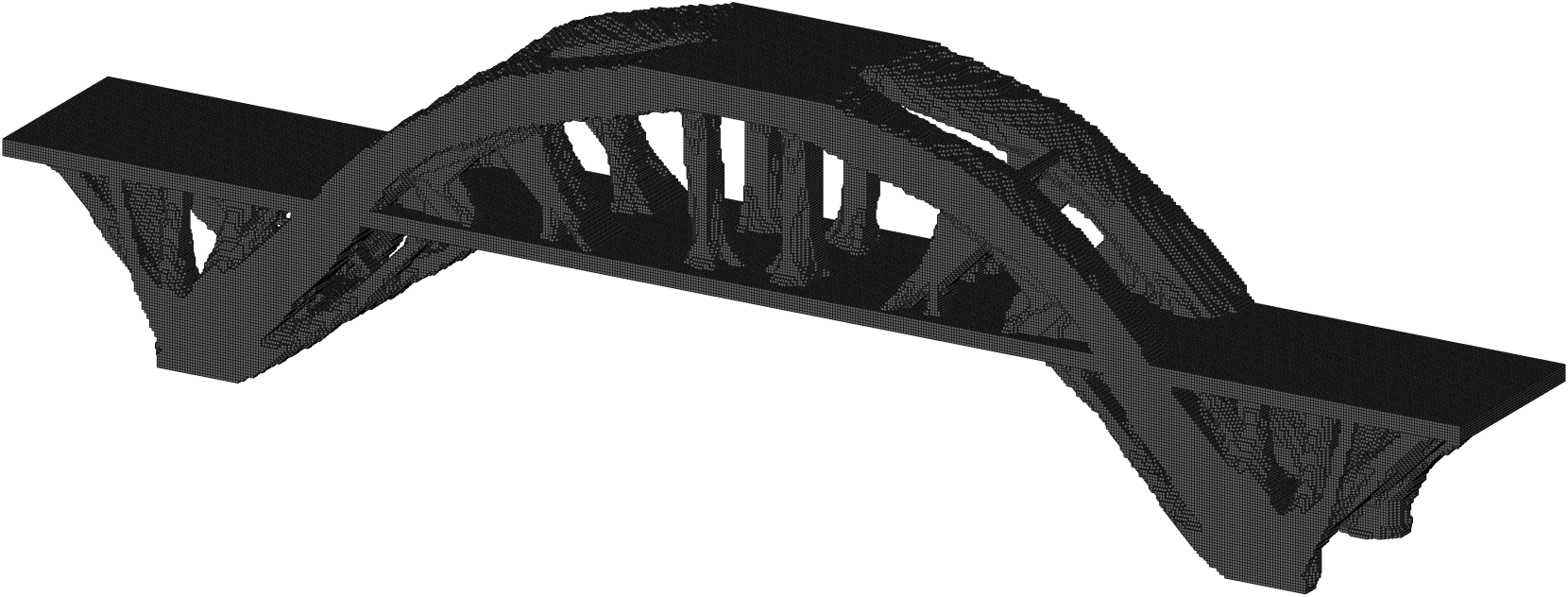}
	\end{subfigure}
	\begin{subfigure}{0.5\linewidth}
		\centering
		\includegraphics[width=0.9\linewidth]{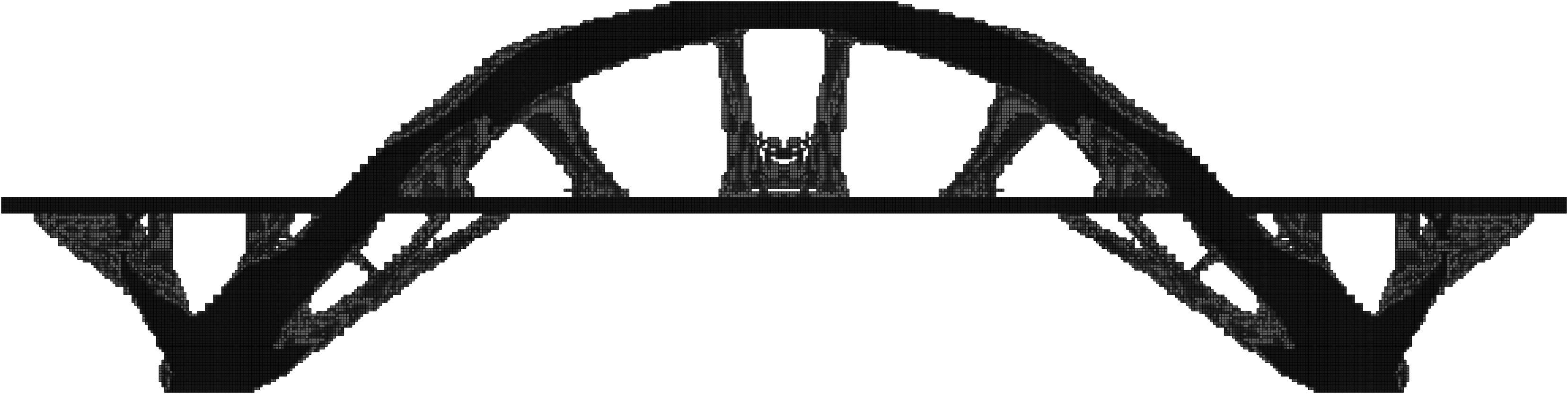}
		\caption{Linear elements ($r_{min} = 0.6$).}
	\end{subfigure}%
	\begin{subfigure}{0.5\linewidth}
		\centering
		\includegraphics[width=0.9\linewidth]{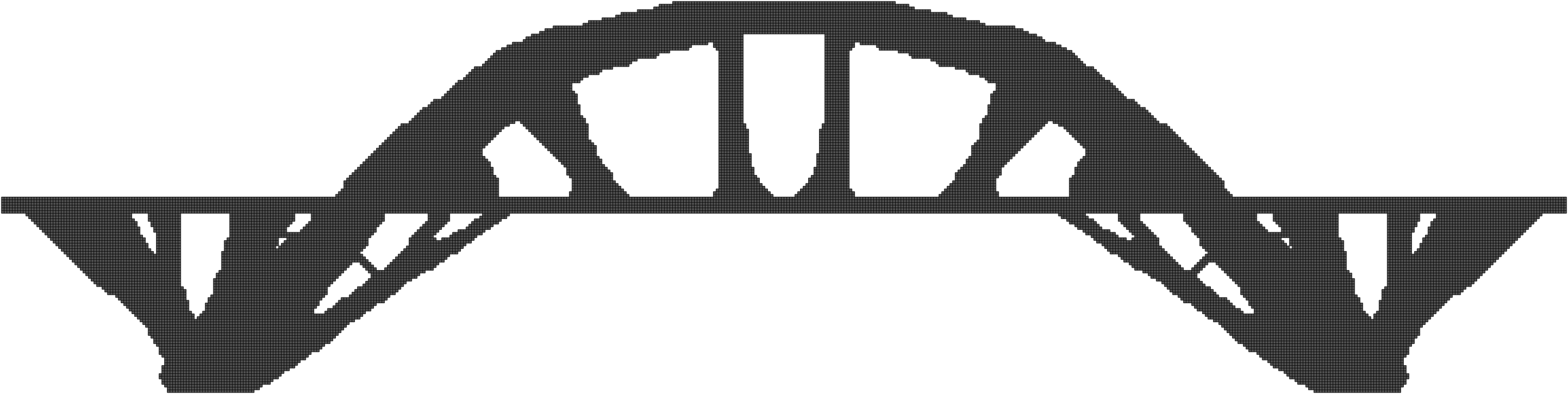}
		\caption{Adaptive strategy with maximum degree $2$ ($r_{min} = 0.6$).}
	\end{subfigure}
	
	\caption{Solutions of the {\tt bd576x144x72} problem with $15\%$ volume, using MR with $n_{mr} = 3$ and $d_{mr} = 2$.}
	\label{fig: bd_fix}
\end{figure*}

In our first experiment, we compare the solution obtained using linear elements with the final solution produced by our adaptive strategy, using quadratic elements but still without fixing variables. Figures \ref{fig: mbb_fix}, \ref{fig: cb_fix}, \ref{fig: ls_fix} and \ref{fig: bd_fix} show the structures obtained for the {\tt mbb576x96x96}, {\tt cb288x96x96}, {\tt ls192x192x64} and {\tt bd576x144x72} problems, respectively, using different values for the filter radius $r_{min}$. 

As expected, more detailed structures are obtained when the filter radius is reduced. However, when the radius is excessively small (generally when $r_{min} < 1.0$), the solutions obtained with linear elements present various artefacts. Fortunately, the adaptive strategy is able to eliminate most of these artefacts, improving the quality of the structure and essentially maintaining the topology obtained in the first step.

Table \ref{tab: fix_F} presents the optimal objective function values for each problem. In this table, $F_0$ corresponds to the compliance for the initial solution (obtained with linear elements), $F$ is the compliance for the solution obtained after solving the problem with quadratic elements and $F_{prj}$ is the compliance of the final solution, obtained after applying the thresholding strategy. To make the figures comparable, all of the values were recalculated in the density mesh with linear elements. The number of density variables of each problem, $n_\rho$, is given only to stress its size.

\begin{table}[h]
	\centering
	\caption{Objective function values obtained with the adaptive strategy.}
	\def\arraystretch{1.2}
	\setlength{\tabcolsep}{8pt}
	\small
	\begin{tabular}{c|c|c|rrr}
		\hline 
		Problem & $n_{\rho}$ & $r_{min}$ & \multicolumn{1}{c}{$F_0$} & \multicolumn{1}{c}{$F$} & \multicolumn{1}{c}{$F_{prj}$} \\ \hline 
		\multirow{2}{*}{\tt mbb576x96x96} & \multirow{2}{*}{1327104} & 1.5 & $2.001 \times 10^{1}$ & $1.981 \times 10^{1}$ & $1.642 \times 10^{1}$ \\
		& & 0.6 & $1.911 \times 10^{1}$ & $1.799 \times 10^{1}$ & $1.706 \times 10^{1}$ \\
		\hline 
		\multirow{2}{*}{\tt cb288x96x96} & \multirow{2}{*}{1327104} &1.5 & $1.957 \times 10^{1}$ & $1.956 \times 10^{1}$ & $1.711 \times 10^{1}$ \\
		& & 0.6 & $2.818 \times 10^{1}$ & $1.848 \times 10^{1}$ & $1.751 \times 10^{1}$ \\
		\hline 
		\multirow{2}{*}{\tt ls192x192x64} & \multirow{2}{*}{1179648} &1.5 & $4.141 \times 10^{1}$ & $4.087 \times 10^{1}$ & $3.016 \times 10^{1}$ \\
		& & 0.6 & $6.089 \times 10^{1}$ & $3.331 \times 10^{1}$ & $3.049 \times 10^{1}$ \\
		\hline 
		\multirow{2}{*}{\tt bd576x144x72} & \multirow{2}{*}{1492992} &1.5 & $9.726 \times 10^{5}$ & $9.714 \times 10^{5}$ & $8.905 \times 10^{5}$ \\
		& & 0.6 & $9.772 \times 10^{5}$ & $9.380 \times 10^{5}$ & $9.067 \times 10^{5}$ \\
		\hline 
	\end{tabular} 
	\label{tab: fix_F}
\end{table}

The objective function values obtained with the adaptive strategy ($F$) are always smaller than those achieved using linear elements ($F_0$). For a filter radius greater than $1.0$, the difference between these two values is small, indicating that the solutions obtained with degree $2$ remain close to the solutions found initially. On the other hand, when $r_{min}<1.0$, an increase on the degree of the elements not only improves the quality of the solution, but also significantly reduces the compliance, especially for the {\tt ls192x192x64} problem. In the following tests, we consider $r_{min} = 0.6$.

Let us now see how the frequency with which we choose the variables to be fixed affects the performance of our algorithm. Henceforth, we will use the notation E0 to designate the original adaptive strategy, where no variable is fixed, and E$d$ to represent each one of the four strategies for fixing variables proposed in Subsection \ref{sec: adapt_fix}, with $n = 5$ for strategies E3 and E4. 
Table \ref{tab: fix_results} presents the numerical results obtained with each strategy for different problems. The number of SLP iterations ($N_{it}$) and the number of rejected steps ($N_{rs}$) with superscripts $1$ and $2$ correspond to those necessary to obtain the initial solution with degree $1$ and the solution with degree $2$, respectively. We observe that, for all problems, just one thresholding attempt was necessary, requiring $N_{it}^{p}$ extra SLP iterations with $N_{rs}^{p}$ rejected steps. 

\begin{table*}[h!]
	\centering
	\caption{Results obtained with the adaptive strategy to increase the degree of the elements.}
	\def\arraystretch{1.2}
	\setlength{\tabcolsep}{2.5pt}
	\small
	\begin{tabular}{c|c|ccc|r|rrrrr|r}
		\hline 
		\!\!\!Problem\!\! & Strat. & $N_{it}^1$ ($N_{rs}^1$) & $N_{it}^2$ ($N_{rs}^2$) & $N_{it}^{p}$ ($N_{rs}^{p}$) & \multicolumn{1}{c|}{$F_{prj}$} & \multicolumn{1}{c}{$T_{sp}$} & \multicolumn{1}{c}{$T_{s}$} & \multicolumn{1}{c}{$T_{g}$} & \multicolumn{1}{c}{$T_{lp}$} & \multicolumn{1}{c|}{$T_{fx}$} & \multicolumn{1}{c}{$T_{total}$}\! \\ \hline 
		\multirow{5}{*}{\begin{tabular}{c}\tt \!\!\!mbb576\!\! \\ \tt \!\!\!x96x96\!\!\end{tabular}} & E0 & 22 (+1) & 33 (+1) & 7 (+0) & 17.064 & 86.6 & 4582.5 & 118.8 & 791.9 & 0.0 & 5917.4\! \\
		& E1 & 22 (+1) & 33 (+1) & 7 (+0) & 17.064 & 43.0 & 2566.0 & 66.0 & 686.4 & 1.0 & 3685.0\! \\
		& E2 & 22 (+1) & 33 (+1) & 7 (+0) & 17.064 & 46.0 & 2668.4 & 65.6 & 688.0 & 18.5 & 3809.7\! \\
		& E3 & 22 (+1) & 33 (+1) & 7 (+0) & 17.064 & 42.9 & 2552.2 & 65.5 & 686.4 & 4.0 & 3674.6\! \\
		& E4 & 22 (+1) & 33 (+1) & 7 (+0) & 17.064 & 42.9 & 2520.0 & 65.9 & 686.4 & 1.7 & 3640.6\! \\
		\hline 
		\multirow{5}{*}{\begin{tabular}{c}\tt cb288\\ \tt \!\!\!x96x96\!\!\end{tabular}} & E0 & 22 (+1) & 79 (+7) & 3 (+0) & 17.513 & 182.0 & 9777.3 & 216.8 & 1260.7 & 0.0 & 12101.5\! \\
		& E1 & 22 (+1) & 70 (+1) & 19 (+1) & 17.496 & 99.0 & 5506.0 & 130.2 & 999.3 & 0.8 & 7387.6\! \\
		& E2 & 22 (+1) & 70 (+1) & 19 (+1) & 17.496 & 101.0 & 5198.7 & 123.0 & 991.3 & 38.7 & 7104.9\! \\
		& E3 & 22 (+1) & 70 (+1) & 19 (+1) & 17.496 & 93.3 & 5208.9 & 123.1 & 989.9 & 8.0 & 7073.7\! \\
		& E4 & 22 (+1) & 70 (+1) & 19 (+1) & 17.496 & 95.0 & 5292.7 & 125.7 & 994.9 & 1.7 & 7162.1\! \\
		\hline 
		\multirow{5}{*}{\begin{tabular}{c}\tt ls192\\ \tt \!\!\!x192x64\!\!\end{tabular}} & E0 & 27 (+0) & 32 (+1) & 3 (+2) & 30.491 & 71.4 & 3574.6 & 52.7 & 208.0 & 0.0 & 4196.8\! \\
		& E1 & 27 (+0) & 32 (+1) & 3 (+2) & 30.491 & 30.7 & 1882.1 & 37.5 & 193.1 & 0.9 & 2421.9\! \\
		& E2 & 27 (+0) & 32 (+1) & 3 (+2) & 30.491 & 33.1 & 1901.8 & 37.0 & 189.0 & 16.2 & 2455.5\! \\
		& E3 & 27 (+0) & 32 (+1) & 3 (+2) & 30.491 & 30.8 & 1861.9 & 37.2 & 188.7 & 3.5 & 2399.7\! \\
		& E4 & 27 (+0) & 32 (+1) & 3 (+2) & 30.491 & 30.5 & 1846.7 & 37.1 & 188.4 & 1.3 & 2381.2\! \\
		\hline 
	\end{tabular} 
	\label{tab: fix_results}
\end{table*}

In general, the number of iterations, as well as the objective function value ($F_{prj}$), was the same for all strategies. The only exception is the {\tt cb288x96x96} problem, for which the algorithm took fewer iterations when the variables were fixed, showing also a slightly smaller final compliance in this case. The structures obtained with all of the strategies were practically identical. 

Table \ref{tab: fix_results} also shows the time (in seconds) required by each strategy to obtain the final solution, which includes the time spent both for solving the initial problem (with linear elements) and for solving the problem again with quadratic elements, as well as the time consumed applying the thresholding strategy. For comparison purposes, in addition to the overall time ($T_{total}$), the table includes the time spent setting up the preconditioner ($T_{sp}$), solving linear systems ($T_s$), computing gradients ($T_g$), solving LP subproblems ($T_{lp}$) and selecting the variables to be fixed ($T_{fx}$). 

Fixing some displacements when solving the problem with elements of degree $2$ reduces the size of the linear systems, so that $T_{sp}$ and $T_s$ are also reduced. Furthermore, since some design variables are also fixed, we see a reduction in $T_g$ and $T_{lp}$ as well. Since the time spent choosing the variables to be fixed ($T_{fx}$) was generally small in relation to the total time, we conclude that all of the proposed strategies for fixing variables are effective in reducing the time spent to solve the problems. 

Comparing the different strategies, we notice that E4 is slightly better for problems {\tt mbb576x96x96} and {\tt ls192x192x64}, and E3 is slightly better for {\tt cb288x96x96}. We remark that strategy E1 should generally be avoided, since fixing the variables based only in the initial solution can prevent the SLP algorithm from converging to a feasible KKT solution, though it did not happen in our tests. On the other hand, strategy E2, in which we select the variables to be fixed at every iteration, can lead to an increase in $T_{sp}$, $T_{s}$ and $T_{fx}$ and, by extension, the total time. Based on these preliminary results, we believe that E3 and E4 are the best strategies. However, one can also try other approaches, such as taking more or less SLP outer iterations to choose the variables to be fixed, or even fixing variables according to other indicator. 

\begin{figure}[h]
  \centering
	\begin{subfigure}{0.45\linewidth}
		\centering
		\includegraphics[width=\linewidth]{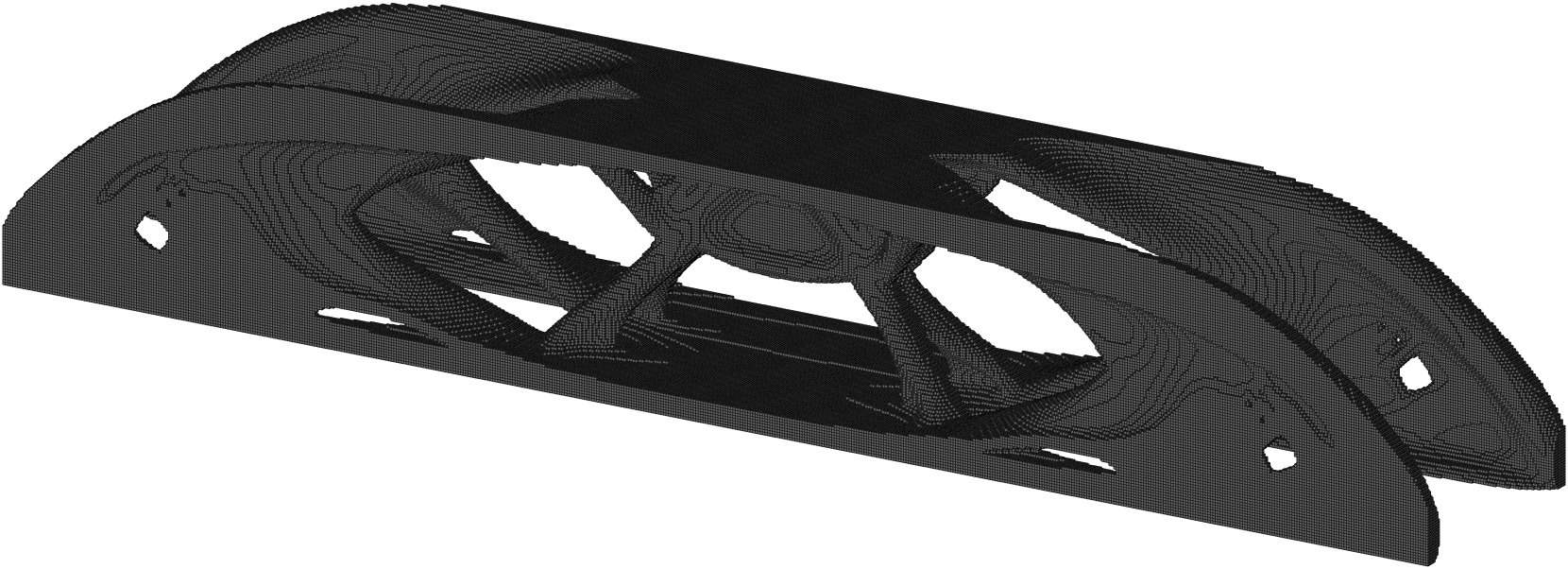}
		\caption{Traditional method.}
	\end{subfigure} \hspace{0.5cm}
	\begin{subfigure}{0.45\linewidth}
		\centering
		\includegraphics[width=\linewidth]{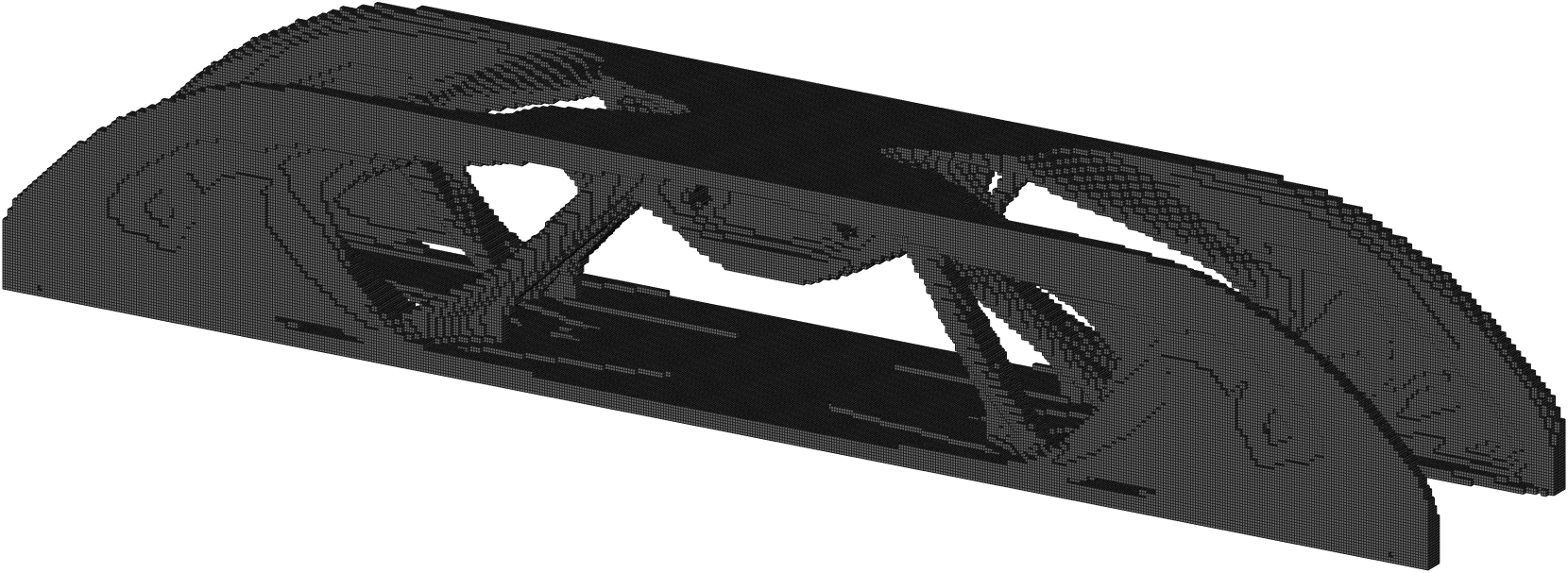}
		\caption{Multiresolution with quadratic elements.}
	\end{subfigure}
	\begin{subfigure}{0.45\linewidth}
		\centering
		\includegraphics[width=\linewidth]{mbb144x24x24quarter_mrn4d2_vol15_r0,6_deg1a2_prj}
		\caption{Proposed adaptive strategy.}
	\end{subfigure}
	\caption{Structures obtained for the {\tt mbb576x96x96} problem using different methods ($r_{min} = 0.6$).}
	\label{fig: fix_mbb_compare}
\end{figure}

\begin{figure}[h]
  \centering
	\begin{subfigure}{0.25\linewidth}
		\centering
		\includegraphics[width=0.95\linewidth]{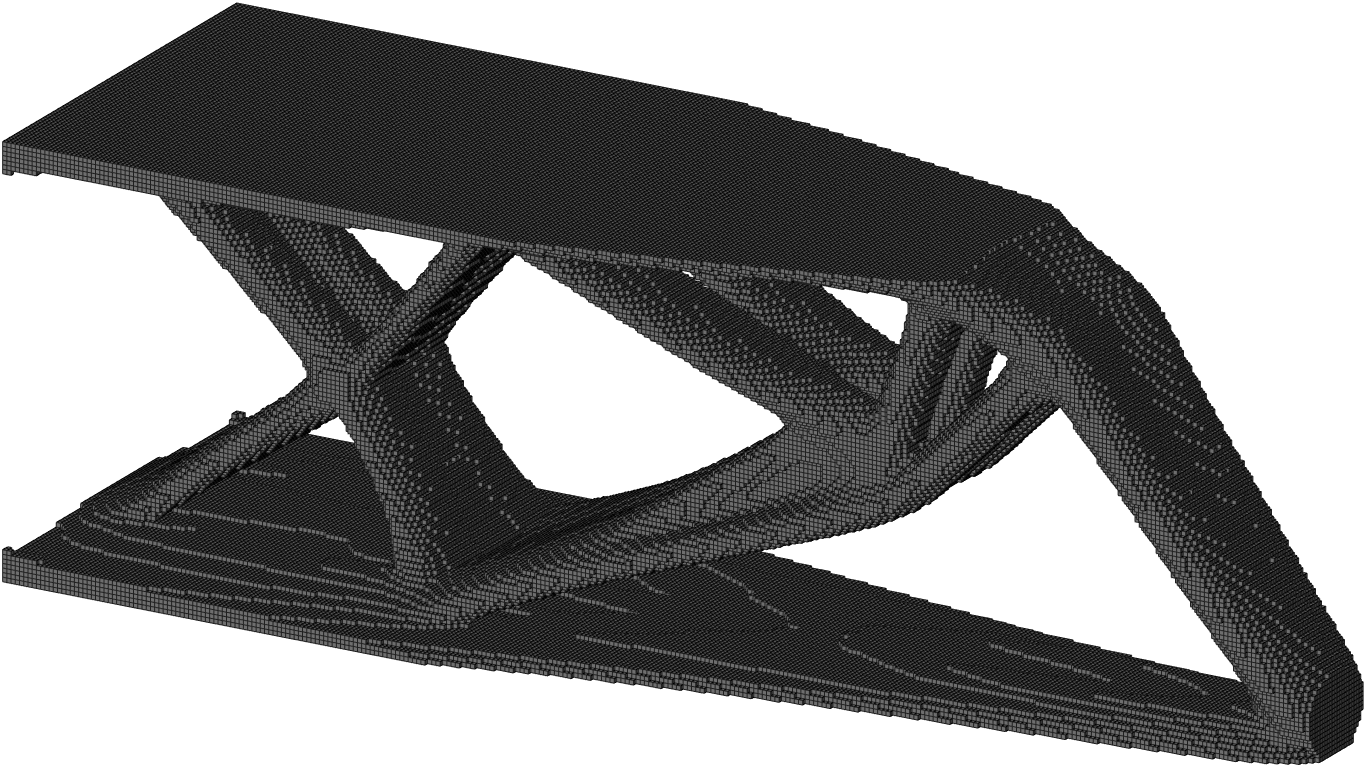}
		\caption{Traditional \\ method. \centering}
	\end{subfigure}%
	\hspace{0.3cm}
	\begin{subfigure}{0.25\linewidth}
		\centering
		\includegraphics[width=0.95\linewidth]{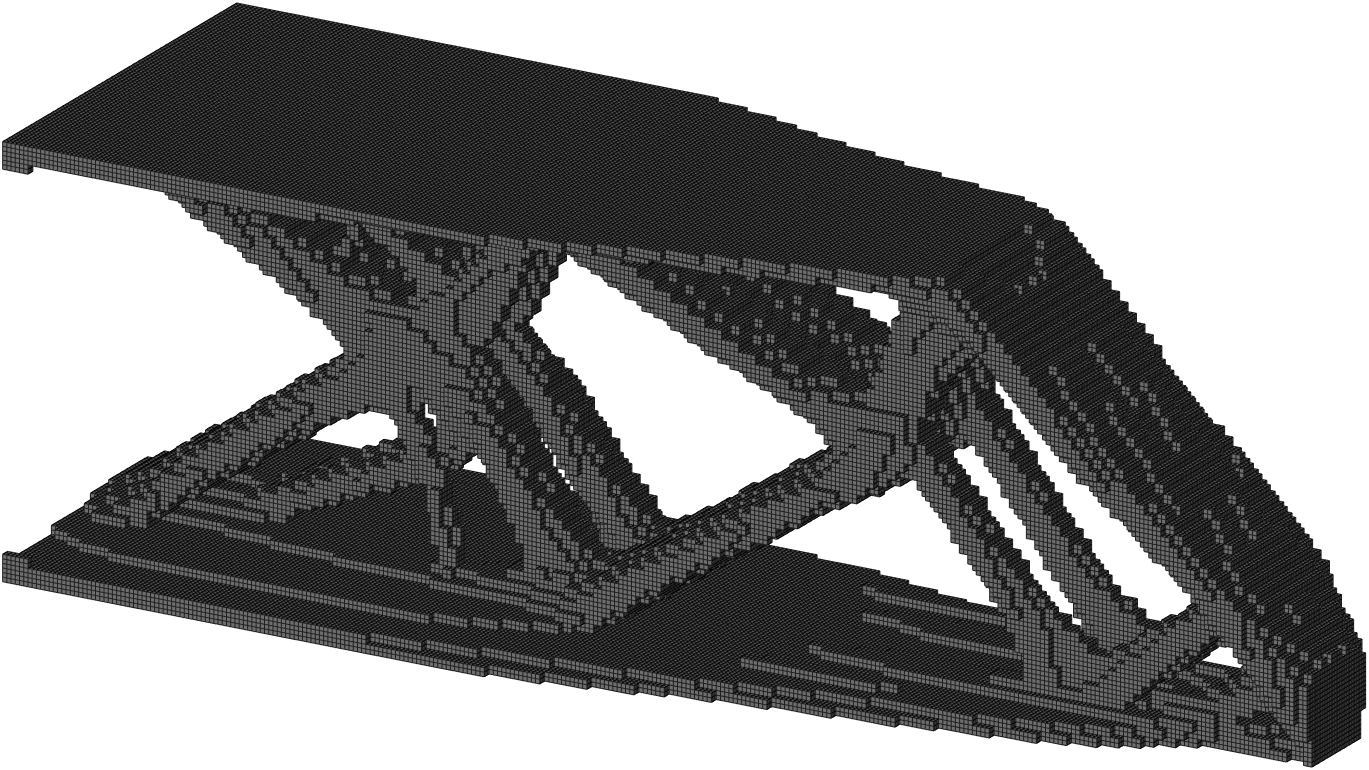}
		\caption{Multiresolution with \\ quadratic elements. \centering}
	\end{subfigure}
	\hspace{0.3cm}
	\begin{subfigure}{0.25\linewidth}
		\centering
		\includegraphics[width=0.95\linewidth]{cb72x24x24half_mrn4d2_vol15_r0,6_deg1a2_prj}
		\caption{Proposed adaptive \\ strategy. \centering}
	\end{subfigure}
	\caption{Structures obtained for the {\tt cb288x96x96} problem using different methods ($r_{min} = 0.6$).}
	\label{fig: fix_cb_compare}
\end{figure}

\begin{figure}[h]
	\centering
	\begin{subfigure}{0.25\linewidth}
		\centering
		\includegraphics[width=0.8\linewidth]{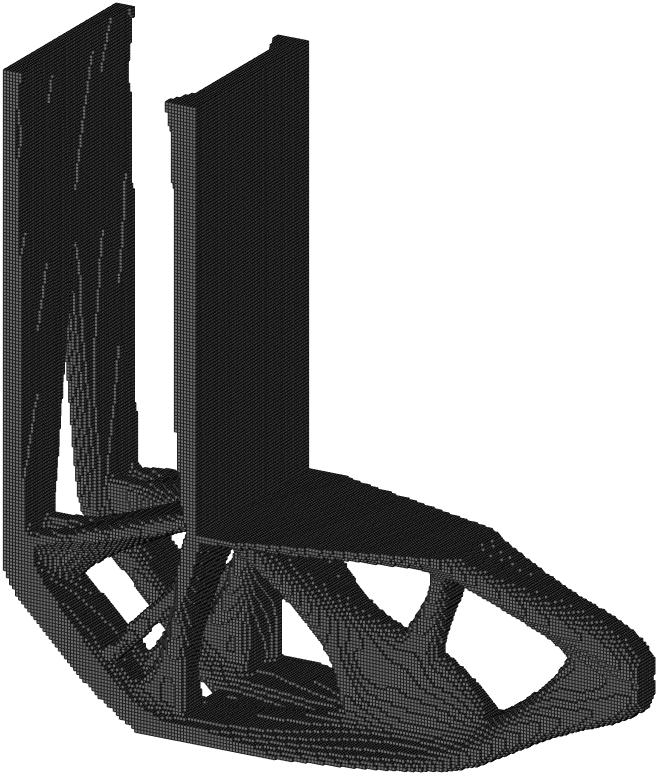}
		\caption{Traditional \\ method. \centering}
	\end{subfigure}%
	\hspace{0.3cm}
	\begin{subfigure}{0.25\linewidth}
		\centering
		\includegraphics[width=0.8\linewidth]{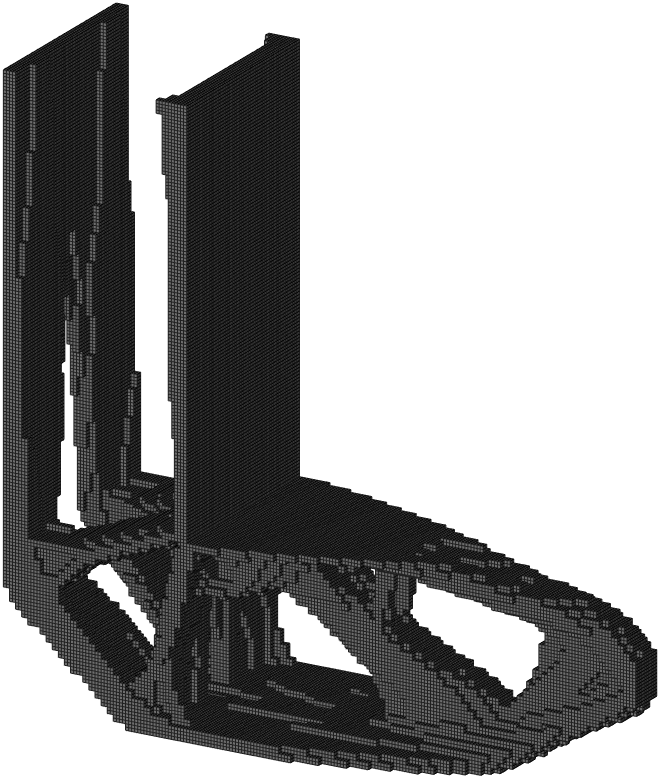}
		\caption{Multiresolution with \\ quadratic elements. \centering}
	\end{subfigure}
	\hspace{0.3cm}
	\begin{subfigure}{0.25\linewidth}
		\centering
		\includegraphics[width=0.8\linewidth]{ls48x48x16half_mrn4d2_r0,6_deg1a2_prj}
		\caption{Proposed adaptive \\ strategy. \centering}
	\end{subfigure}
	\caption{Structures obtained for the {\tt ls192x192x64} problem using different methods ($r_{min} = 0.6$).}
	\label{fig: fix_ls_compare}
\end{figure}

To evaluate the effectiveness of the adaptive strategy, we also compare our results with those obtained using the traditional method (with linear elements and without relying on multiresolution), as well as to those achieved by multiresolution using quadratic elements from the beginning. In our comparison, we adopt strategy E4 to choose the variables that should be fixed. 
Figures \ref{fig: fix_mbb_compare}, \ref{fig: fix_cb_compare} and \ref{fig: fix_ls_compare} show the structures obtained with the three methods for problems {\tt mbb576x96x96}, {\tt cb288x96x96} and {\tt ls192x192x64}, respectively. 

We notice that each method reached a different local minimizer, so that the material distributions are not the same. Besides, the structures obtained with the adaptive strategy are more detailed, while those obtained with multiresolution with degree $2$ from the beginning have the less detailed topologies.

Table \ref{tab: fix_compare} presents the number of iterations and the total time spent (in seconds) by the methods. For the adaptive strategy, the figures in the table include the steps in which the elements have degree $1$ and $2$. As expected, the traditional method usually requires less iterations, but it always takes longer to solve the problems because there are many more variables to consider when solving linear systems and LP subproblems. 
We also notice that our adaptive strategy managed to both reduce the total time and find high quality solutions. However, the compliance values obtained with our strategy are not always the smallest, except for the {\tt cb288x96x96} problem. In short, the proposed strategy offers an interesting trade-off between design performance and computational cost.

\begin{table*}[h!]
	\centering
	\caption{Comparison of the traditional method, the multiresolution using quadratic elements and the adaptive strategy.}
	\def\arraystretch{1.2}
	\setlength{\tabcolsep}{8pt}
	\small
	\begin{tabular}{cc|crrr}
		\hline 
		Problem & Method & \multicolumn{1}{c}{$N_{it}$ ($N_{rs}$)} & \multicolumn{1}{c}{$F$} & \multicolumn{1}{c}{$F_{prj}$} & \multicolumn{1}{c}{$T_{total}$} \\ \hline 
		\multirow{3}{*}{\tt mbb576x96x96} & Traditional & 35 (+1) & 17.861 & 16.616 & 30013.4 \\ 
		& MR with degree $2$ & 23 (+1) & 18.160 & 17.158 & 7199.2 \\
		& Adaptive strategy & 55 (+2) & 17.987 & 17.064 & 3640.6 \\
		\hline 
		\multirow{3}{*}{\tt cb288x96x96} & Traditional & 23 (+1) & 19.018 & 17.534 & 32605.1 \\
		& MR with degree $2$ & 56 (+2) & 18.695 & 17.694 & 10694.4 \\
		& Adaptive strategy & 92 (+2) & 18.460 & 17.496 & 7162.1 \\
		\hline 
		\multirow{3}{*}{\tt ls192x192x64} & Traditional & 28 (+1) & 29.905 & 26.498 & 17337.4 \\ 
		& MR with degree $2$ & 72 (+3) & 33.933 & 31.116 & 11305.0 \\ 
		& Adaptive strategy & 59 (+1) & 33.305 & 30.491 & 2381.2 \\ 
		\hline 
	\end{tabular} 
	\label{tab: fix_compare}
\end{table*}

\begin{figure}[h]
	\begin{subfigure}{0.5\linewidth}
		\centering
		\includegraphics[width=0.9\linewidth]{mbb144x24x24quarter_mrn4d2_vol15_r0,6_deg1a2_prj_L}
		\caption{Maximum degree $2$. \centering}
	\end{subfigure}
	\begin{subfigure}{0.5\linewidth}
		\centering
		\includegraphics[width=0.9\linewidth]{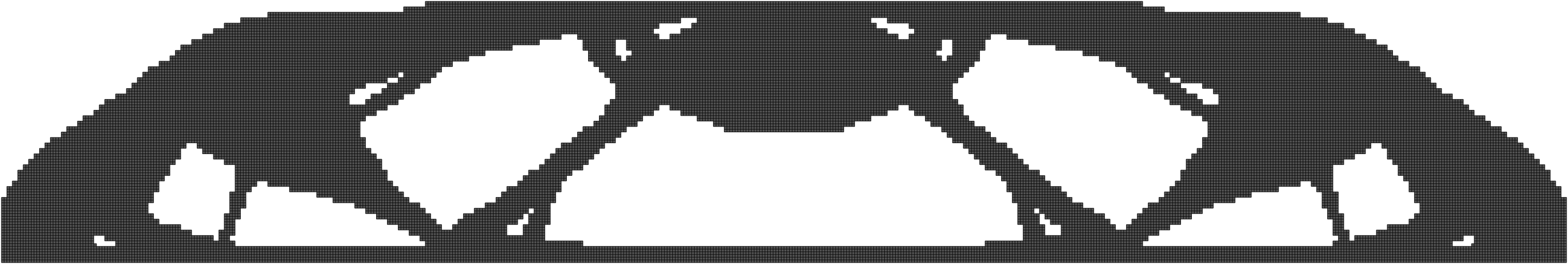}
		\caption{Maximum degree $3$. \centering}
	\end{subfigure}
	
	\caption{Structures obtained for the {\tt mbb576x96x96} problem with $15\%$ volume, using MR with $n_{mr} = 4$, $d_{mr} = 2$, $r_{min} = 0.6$ and the adaptive strategy to increase the degree of the elements.}
	\label{fig: mbb_fix_deg3}
\end{figure}

\begin{figure}[!h]
  \centering
	\begin{subfigure}{0.25\linewidth}
		\centering
		\includegraphics[width=0.8\linewidth]{ls48x48x16half_mrn4d2_r0,6_deg1a2_prj_L}
		\caption{Maximum degree $2$. \centering}
	\end{subfigure}%
	\hspace{0.5cm}
	\begin{subfigure}{0.25\linewidth}
		\centering
		\includegraphics[width=0.8\linewidth]{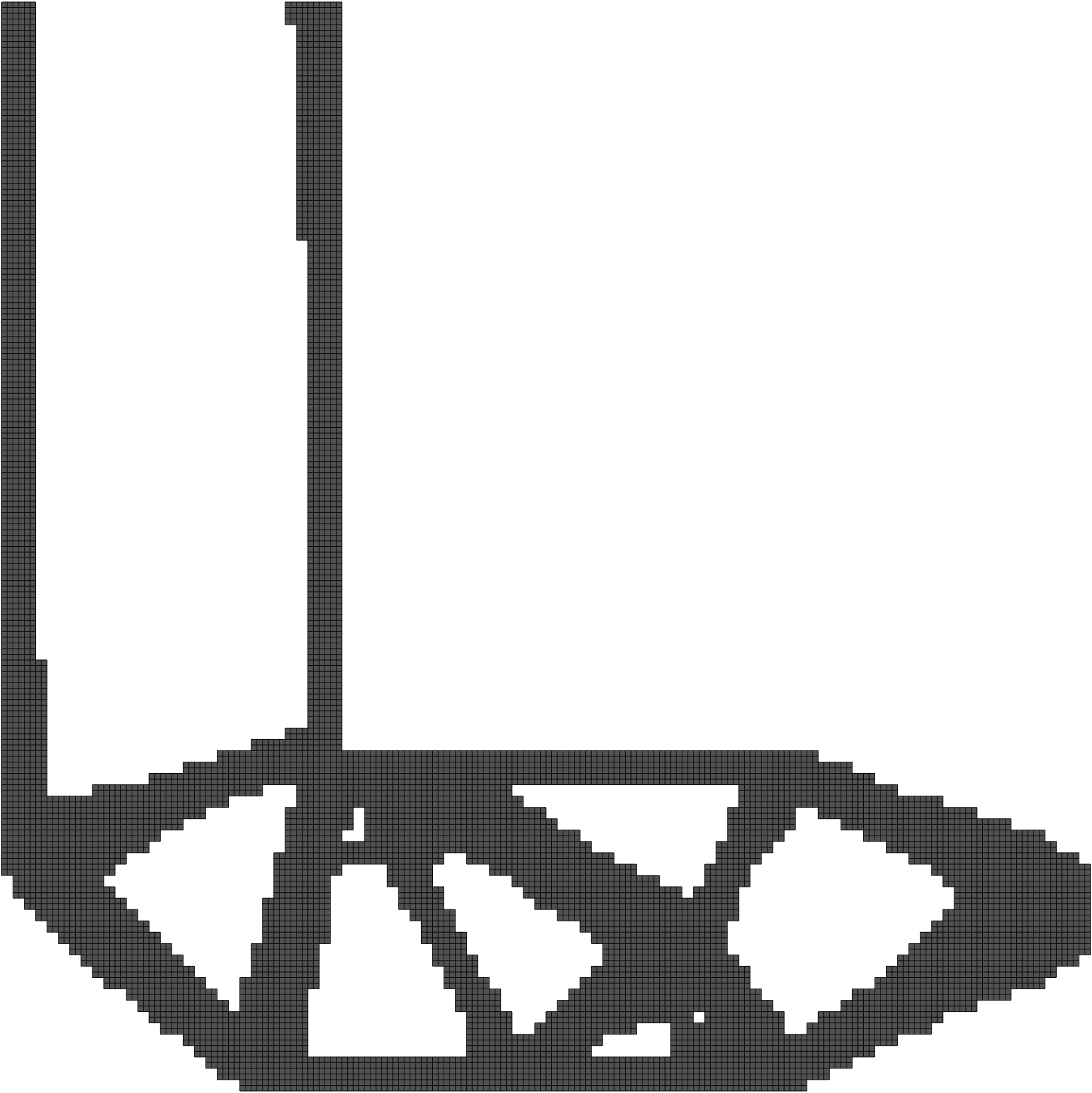}
		\caption{Maximum degree $3$. \centering}
	\end{subfigure}
	
	\caption{Structures obtained for the {\tt ls192x192x64} problem with $18\%$ volume, using MR with $n_{mr} = 4$, $d_{mr} = 2$, $r_{min} = 0.6$ and the adaptive strategy to increase the degree of the elements.}
	\label{fig: ls_fix_deg3}
\end{figure}

So far, we have tested the adaptive strategy increasing the degree of the elements from $1$ to $2$. To complete this section, we show now some results obtained solving the MBB beam and the L-shaped beam problems with elements of degree up to $3$.  
Figures \ref{fig: mbb_fix_deg3} and \ref{fig: ls_fix_deg3} compare the side views of the structures obtained using elements with degree up to $2$ and $3$. Table \ref{tab: fix_deg3} presents the numerical results for these tests.

\begin{table}[!h]
	\centering
	\caption{Results obtained with the adaptive strategy to increase the degree of the elements, considering different maximum degrees.}
	\def\arraystretch{1.2}
	\setlength{\tabcolsep}{8pt}
	\small
	\begin{tabular}{c|c|r|r}
		\hline 
		Problem & Max. degree & \multicolumn{1}{c|}{$F_{prj}$} & \multicolumn{1}{c}{$T_{total}$} \\ \hline 
		\multirow{2}{*}{\tt mbb576x96x96} & $2$ & 17.064 & 3640.6 \\
		& $3$ & 17.089 & 11560.9 \\
		\hline 
		\multirow{2}{*}{\tt ls192x192x64} & $2$ & 30.491 & 2381.2 \\
		& $3$ & 30.482 & 15197.6 \\ 
		\hline 
	\end{tabular} 
	\label{tab: fix_deg3}
\end{table}

In general, the topologies obtained by both strategies are almost indistinguishable. In addition, the variation on the compliance was also negligible. 
On the other hand, increasing the degree of the finite elements to $3$ has a considerable impact on the time spent to solve the problems. Even so, the time spent by the multiresolution algorithm did not exceed that taken by the traditional method, indicating that this approach can be used in some cases. However, it is clear to us that that elements of degree higher than 3 should be avoided, since their use may undermine the efficiency of the multiresolution algorithm.

% ---------- %

\section{Conclusions}

In this work, we combine several techniques to efficiently solve three-dimensional structural topology optimization problems. Our algorithm, tested through a Matlab implementation, was able to solve a variety of problems at a low computational cost. 

To solve the optimization problem, we employ a sequential linear programming method with a stopping criterion based on first order optimality conditions. The results show that, besides being more mathematically sound than its usual counterparts, such as the optimality criteria method and the MMA, our algorithm managed to solve the problems spending few iterations, still obtaining local optimal topologies that meet the purpose of each structure from an engineering point of view. 

We successfully applied both the geometric and the algebraic multigrid as preconditioners for the conjugate gradient method, used to solve the linear equilibrium systems. Geometric multigrid showed better results, but it is restricted to problems with regular meshes and specific discretizations. The algebraic version is more robust and can be applied on a wider variety of problems, although its setup cost is high. 

We propose a new thresholding strategy based on the gradient of the Lagrangian of the problem, with which we substantially reduce the number of intermediate densities, obtaining binary solutions that are closer to local minimizers. 

Combining the SLP method with a multiresolution scheme using three different meshes, we managed to obtain high-resolution structures while reducing the time spent on all steps of the topology optimization process. Using the same parameters of the algorithm, the results obtained with the multiresolution were close to those obtained with the traditional method. 

However, if we raise the number of design variables and density elements but keep fixed the radius of application of the filter (used for the projection of design variables to densities), we get solutions that are mesh independent, so that the increase in resolution did not provide structures with a larger number of members. On the other hand, reducing the filter radius to obtain more detailed structures, we observe a major drawback of multiresolution, that is the occurrence of artefacts. 

According to our tests, when the geometric multigrid method is employed, it is more advantageous to use elements of the Lagrange family with degree at most $2$. For elements of degree greater than $2$, the serendipity family is preferable. Besides that, the solutions obtained with both types of elements were equally good and the use of elements of degree $2$ or $3$ is enough to soften the occurrence of artefacts in three-dimensional structures.

To mitigate the increase in time when problems are solved using quadratic or cubic elements, we propose an adaptive strategy in which we increase the degree of the elements and suppress variables from the problem. This strategy proved to be quite efficient, as we successfully obtained accurate and detailed structures without visible artefacts consuming less time. 

In fact, for large-scale problems, the new algorithm was able to obtain topologies similar to those provided by the traditional method while reducing the solution time by a mean factor of 6.7. Moreover, it allows us to solve even larger problems, for which memory consumption makes the use of the traditional method prohibitive. And although the acceleration strategies proposed here may be combined with other nonlinear optimization methods, the remarkably low number of iterations required by the SLP method leads us to conclude that it is a very competitive algorithm for solving three-dimensional topology optimization problems.

As future work, we intend to extend the methods and strategies proposed here to the solution of more challenging problems, such as those related to compliant mechanisms \cite{Sigmund1997,Zhu} as well as structures with manufacturing constrains \cite{ZegardPaulino,Zhu2021}. For solving these difficult problems, we plan to code the program in C++, so we can use faster libraries for the solution of linear systems and linear programming problems.

\section*{Declarations}

\def\ackname{Ethics approval and consent to participate}
\begin{acknowledgements}
Not applicable.
\end{acknowledgements}

\def\ackname{Consent for publication}
\begin{acknowledgements}
Not applicable.
\end{acknowledgements}

\def\ackname{Funding}
\begin{acknowledgements}
This work was supported by Funda\c c\~ao de Amparo \`a Pesquisa do Estado de S\~ao Paulo (grant 2018/24293-0); and Coordena\c c\~ao de Aperfei\c coamento de Pessoal de N\'\i vel Superior.
\end{acknowledgements}

\def\ackname{Availability of data and materials}
\begin{acknowledgements}
The code and the dataset used in this study are openly available at {\tt https://github.com/AlfredoVitorino/TopOpt-SLP-3D}.
\end{acknowledgements}

\def\ackname{Competing interests}
\begin{acknowledgements}
The authors declare that they have no competing interests.
\end{acknowledgements}

\def\ackname{Authors' contributions}
\begin{acknowledgements}
A.V. and F.G. contributed to the design and implementation of the research, to the analysis of the results and to the writing of the manuscript.
\end{acknowledgements}

% ---------- %

% References 

% ---------- %

\end{document}